\documentclass[reqno, 10pt, centertags,draft]{amsart}
\usepackage{amsmath,amsthm,amscd,amssymb,latexsym,upref,enumerate}

\makeatletter
\def\theequation{\@arabic\c@equation}

\newcommand{\bbN}{{\mathbb{N}}}
\newcommand{\bbR}{{\mathbb{R}}}

\newcommand{\bbZ}{{\mathbb{Z}}}
\newcommand{\bbC}{{\mathbb{C}}}

\newcommand{\cH}{{\mathcal H}}

\newcommand{\cM}{{\mathcal M}}
\newcommand{\cN}{{\mathcal N}}
\newcommand{\cO}{{\mathcal O}}

\newcommand{\no}{\nonumber}
\newcommand{\lb}{\label}
\newcommand{\f}{\frac}

\newcommand{\ol}{\overline}

\newcommand{\wti}{\widetilde}
\newcommand{\ga}{\gamma}
\newcommand{\la}{\lambda}
\newcommand{\al}{\alpha}
\newcommand{\si}{\sigma}
\newcommand{\ve}{\varepsilon}
\newcommand{\Oh}{O}
\newcommand{\oh}{o}

\newcommand{\abs}[1]{\left\lvert#1\right\rvert}

\newcommand{\loc}{\text{\rm{loc}}}

\newcommand{\rank}{\text{\rm{rank}}}

\newcommand{\dom}{\text{\rm{dom}}}

\newcommand{\slim}{\text{\rm{s-lim}}}
\newcommand{\slimes}{\text{\rm{l.i.m.}}}

\newcommand{\supp}{\text{\rm{supp}}}

\newcommand{\bi}{\bibitem}
\newcommand{\hatt}{\widehat}

\newcommand{\essran}{\text{\rm ess.ran}}
\newcommand{\CmR}{\bbC\backslash\bbR}

\renewcommand{\Re}{\text{\rm Re}}
\renewcommand{\Im}{\text{\rm Im}}
\renewcommand{\ln}{\text{\rm ln}}

\numberwithin{equation}{section}

\newtheorem{theorem}{Theorem}[section]

\newtheorem{lemma}[theorem]{Lemma}
\newtheorem{corollary}[theorem]{Corollary}
\newtheorem{hypothesis}[theorem]{Hypothesis}
\newtheorem{example}[theorem]{Example}
\theoremstyle{definition}
\newtheorem{definition}[theorem]{Definition}
\newtheorem{remark}[theorem]{Remark}

\begin{document}

\title[Spectral Theory and Strongly Singular Potentials]{On Spectral
Theory for Schr\"odinger Operators with Strongly Singular Potentials}
\author[F.\ Gesztesy, and M.\ Zinchenko]{Fritz
Gesztesy and Maxim Zinchenko}
\address{Department of Mathematics,
University of Missouri, Columbia, MO 65211, USA}
\email{fritz@math.missouri.edu}
\urladdr{http://www.math.missouri.edu/personnel/faculty/gesztesyf.html}
\address{Department of Mathematics,
University of Missouri, Columbia, MO 65211, USA}
\email{maxim@math.missouri.edu}
\thanks{Based upon work supported by the National Science
Foundation under Grant No.\ DMS-0405526.}
\thanks{{\it Math.\ Nachrichten} {\bf 279}, 1041--1082 (2006).}

\begin{abstract}
We examine two kinds of spectral theoretic situations:  First, we
recall the case of  self-adjoint half-line Schr\"odinger operators
on $[a,\infty)$, $a\in\bbR$, with a regular finite end point $a$ and the
case of Schr\"odinger operators on the real line with locally integrable
potentials, which naturally lead to Herglotz functions and $2\times 2$
matrix-valued Herglotz functions representing the associated
Weyl--Titchmarsh coefficients. Second, we contrast this with the case
of self-adjoint half-line Schr\"odinger operators on $(a,\infty)$ with
a potential strongly singular at the end point $a$. We focus on situations
where the potential is so singular that the associated maximally
defined Schr\"odinger operator is self-adjoint (equivalently, the
associated minimally defined Schr\"odinger operator is essentially
self-adjoint) and hence no boundary condition is required at the finite
end point $a$. For this case we show that the Weyl--Titchmarsh coefficient
in this strongly singular context still determines the associated spectral
function, but ceases to posses the Herglotz property. However, as
will be shown, Herglotz function techniques continue to play a decisive
role in the spectral theory for strongly singular Schr\"odinger operators.
\end{abstract}

\maketitle

\section{Introduction} \lb{s1}

The principal goal of this paper is to study singular Schr\"odinger
operators on a half-line $[a,\infty)$, $a\in\bbR$, with strongly singular
potentials at the finite end point $a$ in the sense that
\begin{equation}
V\in L^1_{\loc}((a,\infty);dx), \, \text{ $V$ real-valued, }
V\notin L^1([a,b];dx), \quad b>a.  \lb{1.1}
\end{equation}
For previous studies of strongly singular Schr\"odinger operators we refer,
for instance, to \cite{AG74}, \cite{At82}--\cite{AF99}, \cite{BG85}, \cite{BS95},
\cite{EK05}--\cite{GP84}, \cite{Me64},
\cite{Na74}, \cite{Pe74}, \cite{Ry05}--\cite{SS03} and the references
therein. (Many of these references treat, in fact, a discrete set of
singularities on
$\bbR$ or on
$(a,b)$, $-\infty\leq a<b\leq\infty$.) Quite recently, singular potentials
became again a popular object of study from various points of views: Some
groups study singular interactions in connections with scales of Hilbert
spaces (see, e.g.,
\cite{DHD03}, \cite{HM03}--\cite{HM04b} and the
references therein), while other groups study strongly singular
interactions in the context of Pontryagin spaces (we refer, e.g., to
\cite{BT04}, \cite{DH03}, \cite{DLS04}, \cite{DS00}, \cite{RS04} and the
references therein).

Our point of departure in connection with strongly singular potentials
is quite different: We focus on the derivation of the spectral function
for strongly singular half-line Schr\"odinger operators starting from
the resolvent (and hence the Green's function). In stark contrast to the
standard situation of Schr\"odinger operators on a half-line
$[a,\infty)$, $a\in\bbR$, with a regular end point $a$, where the
associated spectral function generates the measure in the Herglotz
representation of the Weyl--Titchmarsh coefficient, we show that
half-line Schr\"odinger operators with strongly singular potentials at
the endpoint $a$ lead to spectral functions which are related to the
analog of a Weyl--Titchmarsh coefficient which ceases to be a Herglotz
function. In fact, the strongly singular potentials studied in this
paper are so singular at $a$ that the associated maximally defined
Schr\"odinger operator is self-adjoint (equivalently, the associated
minimal Schr\"odinger operator is essentially self-adjoint) and hence no
boundary condition is required at the finite endpoint $a$.

In Section \ref{s2} we recall the essential ingredients of
standard spectral theory for self-adjoint Schr\"odinger operators
on a half-line $[a,\infty)$, $a\in\bbR$, with a regular end point
$a$ and problems on the real line with locally integrable
potentials. In either case the notion of a spectral function or
$2\times 2$ matrix spectral function is intimately connected with
Herglotz functions and $2\times 2$ Herglotz matrices representing
the celebrated Weyl--Titchmarsh coefficients. This section is, in
part, of an expository nature. In stark contrast to the half-line
case with a regular finite endpoint $a$ in Section \ref{s2}, we
will show in Section \ref{s3} in the case of strongly singular
potentials $V$ on $(a,\infty)$ with singularity concentrated at
the endpoint $a$, that the corresponding spectral functions are no
longer derived from associated Herglotz functions (although,
certain Herglotz functions still play an important role in this
context). We present and contrast two approaches in Section
\ref{s3}: First we discuss the case where the reference point
$x_0$ coincides with the singular endpoint $a$, leading to a
scalar Weyl--Titchmarsh coefficient and a scalar spectral
function. Alternatively, we treat the case where the
reference point $x_0$ belongs to the interior of the interval
$(a,\infty)$, leading to a $2\times 2$ matrix-valued
Weyl--Titchmarsh and spectral function. Finally, in Section
\ref{s4} we provide a detailed discussion of the explicitly
solvable example $V(x)=[\gamma^2-(1/4)]x^{-2}$, $x\in(0,\infty)$,
$\gamma\in[1,\infty)$. Again we illustrate the two approaches with a
choice of reference point $x_0=0$ and $x_0\in(0,\infty)$.

\section{Spectral Theory and Herglotz Functions}
\lb{s2}

In this section we separately recall basic spectral theory for the
case of half-line Schr\"odinger operators with a regular left endpoint
and the case of full-line Schr\"odinger operators with locally
integrable potentials and their relationship to Herglotz functions and
matrices. The material of this section is standard and various parts of
it can be found, for instance, in \cite{BE05}, \cite[Ch.\ 9]{CL85},
\cite[Sect.\ XIII.5]{DS88}, \cite[Ch.\ 2]{EK82},
\cite{Ev04}, \cite[Ch.\ 10]{Hi69}, \cite{HS98}, \cite{Ko49}, \cite{Le51},
\cite[Ch.\ 2]{LG75}, \cite[Ch.\ VI]{Na68}, \cite[Ch.\ 6]{Pe88},
\cite[Chs.\ II, III]{Ti62}, \cite[Sects.\ 7--10]{We87}.

Starting with the half-line case (with a regular left endpoint) we
introduce the following main assumption:

\begin{hypothesis} \lb{h2.1}
$(i)$ Let $a\in\bbR$ and assume that
\begin{equation}
V\in L^1 ([a,c];dx) \, \text{ for all $c\in(a,\infty)$, } \, V \,
\text{real-valued.} \lb{2.1}
\end{equation}
$(ii)$ Introducing the differential expression $\tau_+$ given by
\begin{equation}
\tau_+=-\f{d^2}{dx^2} + V(x), \quad x\in (a,\infty), \lb{2.2}
\end{equation}
we assume $\tau_+$ to be in the limit point case at $+\infty$.
\end{hypothesis}

Associated with the differential expression $\tau_+$ one introduces the
self-adjoint Schr\"odinger operator $H_{+,\alpha}$ in
$L^2([a,\infty);dx)$ by
\begin{align}
&H_{+,\alpha}f=\tau_+ f, \quad \alpha\in [0,\pi), \no \\
&f\in \dom(H_{+,\alpha})=\big\{g\in L^2([a,\infty);dx)\,\big|\, g, g'
\in AC([a,c]) \text{ for all $c\in(a,\infty)$;}    \lb{2.3} \\
& \hspace*{3.1cm}
\sin(\alpha) g'(a_+)+\cos(\alpha) g(a_+)=0; \,
\tau_+ g\in L^2([a,\infty); dx)\big\}. \no
\end{align}
Here (and in the remainder of this manuscript) $\prime$ denotes $d/dx$
and $AC([c,d])$ denotes the class of absolutely continuous functions on
the closed interval $[c,d]$.

\begin{remark} \lb{r2.2}
For simplicity we chose the half-line $[a,\infty)$ rather than a finite
interval $[a,b)$, $a<b<\infty$. Moreover, we chose the limit point
hypothesis of $\tau_+$ at the right end point to avoid having to consider
any boundary conditions at that point. Both limitations can be removed.
\end{remark}

Next, we introduce the standard fundamental system of solutions
$\phi_\alpha(z,\cdot)$ and $\theta_\alpha(z,\cdot)$, $z\in\bbC$,
of
\begin{equation}
(\tau_+ \psi)(z,x) = z \psi(z,x), \quad x\in [a,\infty), \lb{2.4}
\end{equation}
satisfying the initial conditions at the point $x=a$,
\begin{equation}
\phi_\alpha(z,a)=-\theta'_\alpha(z,a)=-\sin(\alpha), \;\;
\phi'_\alpha(z,a)=\theta_\alpha(z,a)=\cos(\alpha), \;\;
\alpha \in [0,\pi). \lb{2.5}
\end{equation}
For future purpose we emphasize that for any fixed $x\in [a,\infty)$,
$\phi_\alpha(z,x)$ and $\theta_\alpha(z,x)$ are entire with
respect to $z$ and that
\begin{equation}
W(\theta_\alpha(z,\cdot),\phi_\alpha(z,\cdot))(x)=1, \quad z\in\bbC,
\lb{2.6}
\end{equation}
where
\begin{equation}
W(f,g)(x)=f(x)g'(x)-f'(x)g(x) \lb{2.7}
\end{equation}
denotes the Wronskian of $f$ and $g$.

A particularly important special solution of \eqref{2.4} is the
{\it Weyl--Titchmarsh solution} $\psi_{+,\alpha}(z,\cdot)$,
$z\in\bbC\backslash\bbR$, uniquely characterized by
\begin{equation}
\psi_{+,\alpha}(z,\cdot)\in L^2([a,\infty);dx), \quad
\sin(\alpha)\psi'_{+,\alpha}(z,a)+\cos(\alpha)\psi_{+,\alpha}(z,a)=1,
\quad z\in\bbC\backslash\bbR.  \lb{2.8}
\end{equation}
The second condition in \eqref{2.8} just determines the normalization
of $\psi_{+,\alpha}(z,\cdot)$ and defines it uniquely. The crucial
condition in \eqref{2.8} is the $L^2$-property which uniquely determines
$\psi_{+,\alpha}(z,\cdot)$ up to constant multiples by the limit point
hypothesis of $\tau_+$ at $\infty$. In particular, for
$\alpha, \beta \in [0,\pi)$,
\begin{equation}
\psi_{+,\alpha}(z,\cdot)=C(z,\alpha,\beta)\psi_{+,\beta}(z,\cdot) \,
\text{ for some coefficient $C(z,\alpha,\beta)\in\bbC$.}  \lb{2.9}
\end{equation}
The normalization in \eqref{2.8} shows that $\psi_{+,\alpha}(z,\cdot)$
is of the type
\begin{equation}
\psi_{+,\alpha}(z,x)=\theta_\alpha(z,x)+m_{+,\alpha}(z)\phi_\alpha(z,x),
\quad z\in\bbC\backslash\bbR, \; x\in [a,\infty) \lb{2.10}
\end{equation}
for some coefficient $m_{+,\alpha}(z)$, the {\it Weyl--Titchmarsh
$m$-function} associated with $\tau_+$ and $\alpha$.

Next, we recall the fundamental identity
\begin{equation}
\int_a^{\infty} dx\,\psi_{+,\alpha}(z_{1},x)
\psi_{+,\alpha}(z_{2},x) = \frac{m_{+,\alpha}(z_{1})-
m_{+, \alpha}(z_{2})}{z_{1}-z_{2}}, \quad z_1, z_2
\in\bbC\backslash\bbR, \; z_1 \neq z_2. \lb{2.11}
\end{equation}
It is a consequence of the elementary fact
\begin{equation}
\f{d}{dx}
W(\psi(z_1,\cdot),\psi(z_2,\cdot))(x)=(z_1-z_2)\psi(z_1,x)\psi(z_2,x)
\lb{2.12}
\end{equation}
for solutions $\psi(z_j,\cdot)$, $j=1,2$, of \eqref{2.4}, and the fact
that $\tau_+$ is assumed to be in the limit point case at $\infty$ which
implies
\begin{equation}
\lim_{x\uparrow\infty} W(\psi_{+,\alpha}(z_{1},\cdot),
\psi_{+,\alpha}(z_{2},\cdot))(x)=0. \lb{2.13}
\end{equation}
Moreover, since $\ol{\psi_{+,\alpha}(z,\cdot)}$ is the unique solution
of $\tau_+ \psi(\ol z,x) = \ol z \psi(\ol z,x)$, $x\in [a,\infty)$,
satisfying
\begin{equation}
\ol{\psi_{+,\alpha}(z,\cdot)} \in L^2([a,\infty);dx), \quad
\sin(\alpha)\ol{\psi'_{+,\alpha}(z,a)}
+\cos(\alpha)\ol{\psi_{+,\alpha}(z,a)}=1, \lb{2.14}
\end{equation}
and since
\begin{equation}
\ol{\phi_\alpha(z,x)} = \phi_\alpha(\ol z,x), \quad
\ol{\theta_\alpha(z,x)} = \theta_\alpha(\ol z,x), \quad
z\in\bbC, \; x\in [a,\infty),   \lb{2.15}
\end{equation}
one concludes that $\ol{\psi_{+,\alpha}(z,\cdot)}$ is the
Weyl--Titchmarsh solution of $\tau_+ \psi(\ol z,x) = \ol z \psi(\ol
z,x)$, $x\geq a$, and hence
\begin{equation}
\ol{m_{+,\alpha}(z)} = m_{+,\alpha}(\ol z), \quad
z\in\bbC\backslash\bbR.
\lb{2.16}
\end{equation}
Thus, choosing $z_1=z$, $z_2=\ol z$ in \eqref{2.11}, one infers
\begin{equation}
\int_a^{\infty} dx\,|\psi_{+,\alpha}(z,x)|^2
                               = \frac{\Im(m_{+,\alpha}(z))}{\Im(z)}, \quad z
\in\bbC\backslash\bbR. \lb{2.17}
\end{equation}

Before we turn to the proper interpretation of formulas \eqref{2.16}
and \eqref{2.17}, we briefly take a look at the Green's function
$G_{+,\alpha}(z,x,x')$ of $H_{+,\alpha}$. Using
\eqref{2.5}, \eqref{2.6}, and \eqref{2.8} one obtains,
\begin{equation}
G_{+,\alpha}(z,x,x') = \begin{cases}
\phi_{\alpha}(z,x)\psi_{+,\alpha}(z,x'), & a\leq x \leq x', \\
\phi_{\alpha}(z,x')\psi_{+,\alpha}(z,x), & a\leq x' \leq x,
\end{cases} \quad z\in\bbC\backslash\bbR \lb{2.18}
\end{equation}
and thus,
\begin{align}
&((H_{+,\alpha}-zI)^{-1}f)(x)
=\int^{\infty}_{a} dx' \, G_{+,\alpha}
(z,x,x')f(x'), \lb{2.19} \\
& \hspace*{1.6cm} z\in\bbC\backslash\bbR, \; x\in[a,\infty), \;
f\in L^{2}([a,\infty);dx). \no
\end{align}

Next we mention the following analyticity result (for the notion of
Herglotz functions we refer to Appendix \ref{sA}). Here and in the
remainder of this manuscript, $\chi_\cM$ denotes the characteristic
function of a set $\cM\subset \bbR$.

\begin{lemma} \lb{l2.3}
Assume Hypothesis \ref{h2.1} and let $\alpha\in [0,\pi)$. Then
$m_{+,\alpha}$ is analytic on $\bbC\backslash\sigma(H_{+,\alpha})$,
moreover, $m_{+,\alpha}$ is a Herglotz function. In addition, for each
$x\in [a,\infty)$, $\psi_{+,\alpha}(\cdot,x)$ and
$\psi'_{+,\alpha}(\cdot,x)$ are analytic on
$\bbC\backslash\sigma(H_{+,\alpha})$.
\end{lemma}
\begin{proof}
Pick real numbers $c$ and $d$ such that $a\leq c<d<\infty$. Then, using
\eqref{2.18} and \eqref{2.19} one computes
\begin{align}
& \int_{\sigma(H_{+,\alpha})}
\f{d \|E_{H_{+,\alpha}}(\lambda)\chi_{[c,d]}\|_{L^2([a,\infty);dx)}^2}{\lambda-z} =
\big(\chi_{[c,d]},(H_{+,\alpha}-zI)^{-1}
\chi_{[c,d]}\big)_{L^2([a,\infty);dx)} \no \\
& \quad = \int_c^d dx \int_c^x dx'\,\theta_\alpha(z,x)
\phi_\alpha(z,x') + \int_c^d dx\int_x^d dx'\,
\phi_\alpha(z,x)\theta_\alpha(z,x')  \lb{2.20} \\
& \qquad + m_{+,\alpha}(z) \bigg[\int_c^d dx \,
\phi_\alpha(z,x)\bigg]^2, \quad
z\in\bbC\backslash\sigma(H_{+,\alpha}).  \no
\end{align}
Since the left-hand side of \eqref{2.20} is analytic with respect to $z$
on $\bbC\backslash\sigma(H_{+,\alpha})$ and since $\phi_\alpha(\cdot,x)$
and $\theta_\alpha(\cdot,x)$ are entire for fixed $x\in
[a,\infty)$ with $\phi_\alpha(z,\cdot)$, $\theta_\alpha(z,\cdot)$, and
their first $x$-derivatives being absolutely continuous on each
interval $[a,b]$, $b>a$, one concludes that
$m_{+,\alpha}$ is analytic in a sufficiently small open
neighborhood $\cN_{z_0}$ of a given point $z_0\in
\bbC\backslash\sigma(H_{+,\alpha})$, as long as we can guarantee
the existence of $c(z_0), d(z_0) \in [a,\infty)$ such that
\begin{equation}
\int_{c(z_0)}^{d(z_0)} dx \,\phi_\alpha(z,x) \neq 0, \quad z\in\cN_{z_0}.
\lb{2.20a}
\end{equation}
The latter is shown as follows: First, pick $z_0\in
\bbC\backslash\sigma(H_{+,\alpha})$. Then since $\phi_{\alpha}(z_0,\cdot)$
does not vanish identically, one can find $c(z_0), d(z_0)\in[a,\infty)$
such that
\begin{equation}
\int_{c(z_0)}^{d(z_0)} dx\, \phi_\alpha(z_0,x) \neq 0. \lb{2.20b}
\end{equation}
Since
\begin{equation}
\int_{c(z_0)}^{d(z_0)} dx\, \phi_\alpha(z,x)
\end{equation}
is entire with respect to $z$, \eqref{2.20b} guarantees the
existence of an open neighborhood $\cN_{z_0}$ of $z_0$ such that
\eqref{2.20a} holds. Since $z_0\in\bbC\backslash\si(H_{+,\al})$
was chosen arbitrary, $m_{+,\al}$ is analytic on
$\bbC\backslash\si(H_{+,\al})$. Together with \eqref{2.16} and
\eqref{2.17} this proves that $m_{+,\al}$ is a Herglotz function.
By \eqref{2.10} (and its $x$-derivative), $\psi_{+,\al}(\cdot,x)$
and $\psi'_{+,\al}(\cdot,x)$ are analytic on
$\bbC\backslash\si(H_{+,\al})$ for each $x\in[a,\infty)$.
\end{proof}

\begin{remark} \lb{r2.4}
Traditionally, one proves analyticity of $m_{+,\alpha}$ on
$\bbC\backslash\bbR$ by first restricting $H_{+,\alpha}$ to the
interval $[a,b]$ (introducing a self-adjoint boundary condition at the
endpoint $b$) and then controls the uniform limit of a sequence of
meromorphic Weyl--Titchmarsh coefficients analytic on
$\bbC\backslash\bbR$ as $b\uparrow\infty$. We chose the somewhat
roundabout proof of Lemma \ref{l2.3} based on the fundamental identity
\eqref{2.20} in view of Section \ref{s3}, in which we consider strongly
singular potentials at $x=a$, where the traditional approach leading to a
Weyl--Titchmarsh coefficient $m_+$ possessing the Herglotz property is not
applicable, but the current method of proof relying on the family of
spectral projections $\{E_{H_{+,\alpha}}\}_{\lambda\in\bbR}$, the Green's
function $G_{+,\alpha}(z,x,x')$ of $H_{+,\alpha}$, and identity
\eqref{2.20}, remains in effect.
\end{remark}

Moreover, we recall the following well-known facts on $m_{+,\alpha}$:
\begin{align}
&\lim\limits_{\epsilon\downarrow 0}\, i\epsilon\, m_{+,\alpha}
(\lambda+i\epsilon) =\begin{cases}
0, &\phi_{\alpha}(\lambda,\,\cdot\,)\notin L^{2}([a,\infty);dx), \\
-\|\phi_{\alpha}(\lambda, \,\cdot\,)\|^{-2}_{L^2([a,\infty);dx)},
&\phi_{\alpha}
(\lambda, \,\cdot\,)\in L^{2}([a,\infty);dx), \end{cases} \lb{2.21}  \\
& \hspace*{7.67cm} \lambda\in\bbR, \; \alpha\in [0,\pi), \no \\
&m_{+,\alpha_{1}}(z) =\frac{-\sin(\alpha_{1}-\alpha_{2}) +
\cos(\alpha_{1}-\alpha_{2}) m_{+, \alpha_{2}}(z)}
{\cos(\alpha_{1}-\alpha_{2}) +\sin(\alpha_{1}-\alpha_{2})
m_{+,\alpha_{2}}(z)}, \quad \alpha_1, \alpha_2 \in [0,\pi), \lb{2.22} \\
&m_{+,\alpha}(z) \underset{z\to i\infty}{=} \begin{cases}
\cot(\alpha)+\frac{i}{\sin^{2}(\alpha)}\, z^{-1/2} - \frac
{\cos(\alpha)}{\sin^{3}(\alpha)}\,z^{-1} + o(z^{-1}), & \alpha\in
(0,\pi), \\ iz^{1/2} + o(1), & \alpha=0. \end{cases} \lb{2.23}
\end{align}
The asymptotic behavior \eqref{2.23} then implies the
Herglotz representation of $m_{+,\alpha}$ (cf.\ Theorem
\ref{tA.2}\,$(iii)$),
\begin{equation}
m_{+,\alpha}(z) = \begin{cases} c_{+,\alpha}
+ \int_{\bbR}  d\rho_{+,\alpha}(\lambda) \Big[
\frac{1}{\lambda-z} -\frac{\lambda}{1+\lambda^2}\Big], &
\alpha\in [0,\pi),  \\[3mm]
\cot(\alpha) + \int_{\bbR} d\rho_{+,\alpha}(\lambda)\,
(\lambda-z)^{-1}, & \alpha\in (0,\pi), \end{cases}
\quad z\in\bbC\backslash\bbR  \lb{2.25}
\end{equation}
with
\begin{equation}
\int_{\bbR} \frac{d\rho_{+,\alpha}(\lambda)}
{1+|\lambda|} \, \begin{cases} < \infty, & \alpha\in (0,\pi), \\
=\infty, & \alpha =0, \end{cases} \quad
\int_{\bbR} \frac{d\rho_{+,0}(\lambda)}
{1+\lambda^2} < \infty. \lb{2.26}
\end{equation}

We note that in formulas \eqref{2.10}--\eqref{2.25} one
can of course replace $z\in\bbC\backslash\bbR$ by
$z\in\bbC\backslash\sigma(H_{+,\alpha})$.

For future purposes we also note the following result, a version of
Stone's formula in the weak sense (cf., e.g., \cite[p.\ 1203]{DS88}).

\begin{lemma} \lb{l2.4a}
Let $T$ be a self-adjoint operator in a complex separable Hilbert space
$\cH$ $($with scalar product denoted by $(\cdot,\cdot)_\cH$, linear in the
second factor$)$ and denote by $\{E_T(\lambda)\}_{\lambda\in\bbR}$ the
family of self-adjoint right-continuous spectral projections associated
with $T$, that is, $E_T(\lambda)=\chi_{(-\infty,\lambda]}(T)$,
$\lambda\in\bbR$. Moreover, let $f,g \in\cH$, $\lambda_1,\lambda_2\in\bbR$,
$\lambda_1<\lambda_2$, and $F\in C(\bbR)$. Then,
\begin{align}
&(f,F(T)E_{T}((\lambda_1,\lambda_2])g)_{\cH} \no \\
& \quad = \lim_{\delta\downarrow 0}\lim_{\varepsilon\downarrow 0}
\frac{1}{2\pi i}
\int_{\lambda_1+\delta}^{\lambda_2+\delta} d\lambda \, F(\lambda)
\big[\big(f,(T-(\lambda+i\varepsilon) I_{\cH})^{-1}g\big)_{\cH}  \no \\
& \hspace*{4.9cm} - \big(f,(T-(\lambda-i\varepsilon)I_{\cH})^{-1}
g\big)_{\cH}\big]. \lb{2.26a}
\end{align}
\end{lemma}
\begin{proof}
First, assume $F\geq 0$. Then
\begin{align}
& \big(F(T)^{1/2}E_{T}((\lambda_1,\lambda_2])f, (T-z I_{\cH})^{-1}
F(T)^{1/2}E_{T}((\lambda_1,\lambda_2])f\big)_{\cH} \no \\
& \quad = \int_\bbR d\big(f,E_T(\lambda)f\big)_{\cH} \, F(\lambda)
\chi_{(\lambda_1,\lambda_2]}(\lambda)(\lambda-z)^{-1} \no \\
& \quad = \int_\bbR \f{d\big(F(T)^{1/2}
\chi_{(\lambda_1,\lambda_2]}(T)f,E_T(\lambda)
F(T)^{1/2}\chi_{(\lambda_1,\lambda_2]}(T)f\big)_{\cH}}{(\lambda-z)},
\quad z\in\bbC_+  \lb{2.26b}
\end{align}
is a Herglotz function and hence \eqref{2.26a} for $g=f$ follows from
\eqref{A.4}. If $F$ is not nonnegative, one decomposes $F$ as
$F=(F_1-F_2)+i(F_3-F_4)$ with $F_j\geq 0$, $1\leq j \leq 4$ and
applies \eqref{2.26b} to each $j\in\{1,2,3,4\}$. The general case
$g\neq f$ then follows from the case $g=f$ by polarization.
\end{proof}

Next, we relate the family of spectral projections,
$\{E_{H_{+,\alpha}}(\lambda)\}_{\lambda\in\bbR}$, of the self-adjoint
operator $H_{+,\alpha}$ and the spectral function
$\rho_{+,\alpha}(\lambda)$, $\lambda\in\bbR$, which generates the
measure in the Herglotz representation \eqref{2.25} of $m_{+,\alpha}$.

We first note that for $F\in C(\bbR)$,
\begin{align}
&(f,F(H_{+,\alpha})g)_{L^2([a,\infty);dx)}= \int_{\bbR}
d (f,E_{H_{+,\alpha}}(\lambda)g)_{L^2([a,\infty);dx)}\,
F(\lambda), \no \\
& f, g \in\dom(F(H_{+,\alpha})) \lb{2.27} \\
& \quad =\bigg\{h\in L^2([a,\infty);dx)
\,\bigg|\,
\int_{\bbR} d \|E_{H_{+,\alpha}}(\lambda)h\|_{L^2([a,\infty);dx)}^2
\, |F(\lambda)|^2 < \infty\bigg\}. \no
\end{align}
Equation \eqref{2.27} extends to measurable functions $F$ and holds
also in the strong sense, but the displayed weak version will suffice
for our purpose.

In the following, $C_0^\infty((c,d))$, $-\infty \leq c<d\leq \infty$,
denotes the usual space of infinitely differentiable functions of
compact support contained in $(c,d)$.

\begin{theorem} \lb{t2.5}
Let $\alpha\in [0,\pi)$, $f,g \in C^\infty_0((a,\infty))$,
$F\in C(\bbR)$, and $\lambda_1, \lambda_2 \in\bbR$,
$\lambda_1<\lambda_2$. Then,
\begin{equation}
(f,F(H_{+,\alpha})E_{H_{+,\alpha}}((\lambda_1,\lambda_2])
g)_{L^2([a,\infty);dx)} =  \big(\hatt
f_{+,\alpha},M_FM_{\chi_{(\lambda_1,\lambda_2]}} \hatt
g_{+,\alpha}\big)_{L^2(\bbR;d\rho_{+,\alpha})}, \lb{2.28}
\end{equation}
where we introduced the notation
\begin{equation}
\hatt h_{+,\alpha}(\lambda)=\int_a^\infty dx \,
\phi_\alpha(\lambda,x)h(x), \quad \lambda \in\bbR, \;
h\in C^\infty_0((a,\infty)), \lb{2.29}
\end{equation}
and $M_G$ denotes the maximally defined operator of multiplication by the
$d\rho_{+,\alpha}$-measurable function $G$ in
the Hilbert space $L^2(\bbR;d\rho_{+,\alpha})$,
\begin{align}
\begin{split}
& \big(M_G\hatt h\big)(\lambda)=G(\lambda)\hatt h(\lambda)
\, \text{ for $d\rho_{+,\alpha}$-a.e.\ $\lambda\in\bbR$}, \lb{2.30} \\
& \hatt h\in\dom(M_G)=\big\{\hatt k \in L^2(\bbR;d\rho_{+,\alpha}) \,\big|\,
G\hatt k \in L^2(\bbR;d\rho_{+,\alpha})\big\}.
\end{split}
\end{align}
Here $d\rho_{+,\alpha}$ is the measure in the Herglotz representation of
the Weyl--Titchmarsh function $m_{+,\alpha}$ $($cf.\ \eqref{2.25}$)$.
\end{theorem}
\begin{proof}
The point of departure for deriving \eqref{2.28} is Stone's formula
\eqref{2.26a} applied to $T=H_{+,\alpha}$,
\begin{align}
&(f,F(H_{+,\alpha})E_{H_{+,\alpha}}((\lambda_1,\lambda_2])
g)_{L^2([a,\infty);dx)} \no \\ & \quad =
\lim_{\delta\downarrow 0}\lim_{\varepsilon\downarrow 0}
\frac{1}{2\pi i} \int_{\lambda_1+\delta}^{\lambda_2+\delta}
d\lambda \, F(\lambda)
\big[\big(f,(H_{+,\alpha}-(\lambda+i\varepsilon)
I)^{-1}g\big)_{L^2([a,\infty);dx)}
\no \\
& \hspace*{4.9cm} - \big(f,(H_{+,\alpha}-(\lambda-i\varepsilon)
I)^{-1}g\big)_{L^2([a,\infty);dx)}\big]. \lb{2.31}
\end{align}
Insertion of \eqref{2.18} and \eqref{2.19} into \eqref{2.31} then
yields the following:
\begin{align}
&(f,F(H_{+,\alpha})E_{H_{+,\alpha}}((\lambda_1,\lambda_2])
g)_{L^2([a,\infty);dx)} = \lim_{\delta\downarrow
0}\lim_{\varepsilon\downarrow 0} \frac{1}{2\pi i}
\int_{\lambda_1+\delta}^{\lambda_2+\delta} d\lambda \, F(\lambda) \no \\
& \quad \times \int_a^\infty dx
\bigg\{\bigg[\ol{f(x)} \psi_{+,\alpha}(\lambda
+i\varepsilon,x) \int_a^x dx'\,
\phi_\alpha(\lambda+i\varepsilon,x')g(x') \no \\
& \hspace*{2.4cm} + \ol{f(x)} \phi_\alpha(\lambda+i\varepsilon,x)
\int_x^\infty dx'\, \psi_{+,\alpha}(\lambda+i\varepsilon,x') g(x')
\bigg]  \no \\
& \hspace*{1.9cm}  -\bigg[\ol{f(x)}
\psi_{+,\alpha}(\lambda -i\varepsilon,x) \int_a^x dx'\,
\phi_\alpha(\lambda-i\varepsilon,x')g(x') \no \\
& \hspace*{2.4cm}  + \ol{f(x)}
\phi_\alpha(\lambda-i\varepsilon,x)
\int_x^\infty dx' \,\psi_{+,\alpha}(\lambda-i\varepsilon,x')
g(x')\bigg]\bigg\}.  \lb{2.32}
\end{align}
Freely interchanging the $dx$ and $dx'$ integrals with
the limits and the $d\lambda$ integral (since all integration
domains are finite and all integrands are continuous), and inserting
expression \eqref{2.10} for $\psi_{+,\alpha}(z,x)$ into \eqref{2.32}, one
obtains
\begin{align}
&(f,F(H_{+,\alpha})E_{H_{+,\alpha}}((\lambda_1,\lambda_2])
g)_{L^2([a,\infty);dx)}
=\int_a^\infty dx\, \ol{f(x)}\bigg\{\int_a^x dx' \, g(x') \no \\
& \qquad \times \lim_{\delta\downarrow
0}\lim_{\varepsilon\downarrow 0} \frac{1}{2\pi i}
\int_{\lambda_1+\delta}^{\lambda_2+\delta} d\lambda \, F(\lambda)
\Big[\big[\theta_\alpha(\lambda,x) +
m_{+,\alpha}(\lambda+i\varepsilon)\phi_\alpha(\lambda,x)\big]
\phi_\alpha(\lambda,x') \no \\
& \hspace*{4cm} -\big[\theta_\alpha(\lambda,x) +
m_{+,\alpha}(\lambda-i\varepsilon)\phi_\alpha(\lambda,x)\big]
\phi_\alpha(\lambda,x')\Big] \no \\
& \quad +\int_x^\infty dx'\, g(x') \lim_{\delta\downarrow 0}
\lim_{\varepsilon\downarrow 0} \frac{1}{2\pi i}
\int_{\lambda_1+\delta}^{\lambda_2+\delta} d\lambda \, F(\lambda)
\lb{2.33} \\
& \hspace*{3.5cm} \times \Big[
                   \phi_\alpha(\lambda,x) \big[\theta_\alpha(\lambda,x') +
m_{+,\alpha}(\lambda+i\varepsilon)\phi_\alpha(\lambda,x')\big]  \no \\
& \hspace{4cm} -\phi_\alpha(\lambda,x)
\big[\theta_\alpha(\lambda,x') +
m_{+,\alpha}(\lambda-i\varepsilon)\phi_\alpha(\lambda,x')\big]\Big]\bigg\}.
\no
\end{align}
Here we employed the fact that for fixed $x\in [a,\infty)$,
$\theta_\alpha(z,x)$ and $\phi_\alpha(z,x)$ are entire with
respect to  $z$, that $\theta_\alpha(\la,x)$ and
$\phi_\alpha(\la,x)$ are real-valued for $\la\in\bbR$, that
$\theta_\alpha(z,\cdot), \phi_\alpha(z,\cdot) \in AC([a,c])$ for
all $c>a$, and hence that
\begin{align}
\begin{split}
\theta_\alpha(\lambda\pm i\varepsilon,x)
&\underset{\varepsilon\downarrow 0}{=}
\theta_\alpha(\lambda,x) \pm
i\varepsilon(d/dz)\theta_\alpha(z,x)|_{z=\lambda} + \Oh(\ve^2),
\lb{2.33A} \\
\phi_\alpha(\lambda\pm i\varepsilon,x)
&\underset{\varepsilon\downarrow 0}{=} \phi_\alpha(\lambda,x)
\pm i\varepsilon(d/dz)\phi_\alpha(z,x)|_{z=\lambda} + \Oh(\ve^2)
\end{split}
\end{align}
with $\Oh(\varepsilon^2)$ being uniform with respect to
$(\lambda,x)$ as long as $\lambda$ and $x$ vary in
compact subsets of $\bbR\times [a,\infty)$. Moreover, we used that
\begin{align}
\begin{split}
&\varepsilon|m_{+,\alpha}(\lambda+i\varepsilon)|\leq
C(\lambda_1,\lambda_2,\varepsilon_0) \, \text{ for } \, \lambda\in
[\lambda_1,\lambda_2], \; 0<\varepsilon\leq\varepsilon_0,
\lb{2.33a} \\
&\varepsilon |\Re(m_{+,\alpha}(\lambda+i\varepsilon))|
\underset{\varepsilon\downarrow 0}{=}\oh(1), \quad \lambda\in \bbR.
\end{split}
\end{align}
In particular, utilizing \eqref{2.33A} and \eqref{2.33a},
$\phi_\alpha(\lambda\pm i\varepsilon,x)$ and
$\theta_\alpha(\lambda\pm i\varepsilon,x)$ have been replaced by
$\phi_\alpha(\lambda,x)$ and $\theta_\alpha(\lambda,x)$ under the
$d\lambda$ integrals in \eqref{2.33}. Cancelling appropriate terms in
\eqref{2.33},  simplifying the remaining terms, and using \eqref{2.16}
then yield
\begin{align}
&(f,F(H_{+,\alpha})E_{H_{+,\alpha}}((\lambda_1,\lambda_2])
g)_{L^2([a,\infty);dx)}
=\int_a^\infty dx\, \ol{f(x)}\int_a^\infty dx' \, g(x') \no \\
& \quad \times \lim_{\delta\downarrow 0}\lim_{\varepsilon\downarrow 0}
\frac{1}{\pi}
\int_{\lambda_1+\delta}^{\lambda_2+\delta} d\lambda \, F(\lambda)
\phi_\alpha(\lambda,x)\phi_\alpha(\lambda,x')
\Im(m_{+,\alpha}(\lambda+i\varepsilon)). \lb{2.34}
\end{align}
Using the fact that by \eqref{A.4}
\begin{equation}
\int_{(\lambda_1,\lambda_2]} d\rho_{+,\alpha}(\lambda) =
\rho_{+,\alpha}((\lambda_1,\lambda_2])=
\lim_{\delta\downarrow 0}\lim_{\varepsilon\downarrow 0}
\frac{1}{\pi}\int_{\lambda_1+\delta}^{\lambda_2+\delta} d\lambda \,
\Im(m_{+,\alpha}(\lambda+i\varepsilon)), \lb{2.35}
\end{equation}
and hence that
\begin{align}
\int_{\bbR} d\rho_{+,\alpha}(\lambda)\, h(\lambda) &=
\lim_{\varepsilon\downarrow 0} \frac{1}{\pi}\int_{\bbR} d\lambda \,
\Im(m_{+,\alpha}(\lambda+i\varepsilon))\, h(\lambda), \quad
h\in C_0(\bbR), \lb{2.36} \\
\int_{(\lambda_1,\lambda_2]} d\rho_{+,\alpha}(\lambda)\, k(\lambda) &=
\lim_{\delta\downarrow 0} \lim_{\varepsilon\downarrow 0} \frac{1}{\pi}
\int_{\lambda_1+\delta}^{\lambda_2+\delta} d\lambda \,
\Im(m_{+,\alpha}(\lambda+i\varepsilon))\, k(\lambda), \quad
k\in C(\bbR), \lb{2.37}
\end{align}
(with $C_0(\bbR)$ the space of continuous compactly supported functions
on $\bbR$) one concludes
\begin{align}
&(f,F(H_{+,\alpha})E_{H_{+,\alpha}}((\lambda_1,\lambda_2])
g)_{L^2([a,\infty);dx)} \no \\ & \quad =\int_a^\infty dx\,
\ol{f(x)}\int_a^\infty dx' \, g(x') \int_{(\lambda_1,\lambda_2]}
d\rho_{+,\alpha}(\lambda) F(\lambda)
\phi_\alpha(\lambda,x)\phi_\alpha(\lambda,x') \no \\
& \quad = \int_{(\lambda_1,\lambda_2]} d\rho_{+,\alpha}(\lambda)
F(\lambda) \, \ol{\hatt f_{+,\alpha}(\lambda)}\, \hatt
g_{+,\alpha}(\lambda), \lb{2.38}
\end{align}
using \eqref{2.29} and interchanging the $dx$, $dx'$ and
$d\rho_{+,\alpha}$ integrals once more.
\end{proof}

\begin{remark} \lb{r2.6}
Theorem \ref{t2.5} is of course well-known. We presented a detailed proof
since this proof will serve as the model for generalizations to
strongly singular potentials and hence pave the way into somewhat
unchartered territory in Section \ref{s3}. In this context it is
worthwhile to examine the principal ingredients entering the proof of
Theorem \ref{t2.5}: Let $\lambda_j \in\bbR$, $j=0,1,2$, $\lambda_1 <
\lambda_2$, and $\varepsilon_0 >0$. Then the following items played a
crucial role in the proof of Theorem \ref{t2.5}:
\begin{align}
& (i) \;\;\, \text{For each $x\in [a,\infty)$, $\theta_\alpha(z,x)$ and
$\phi_\alpha(z,x)$ are entire with respect to $z$} \no \\
& \qquad \text{and real-valued for $z\in\bbR$.} \no \\
& (ii) \;\; \text{$m_{+,\alpha}$ is analytic on
$\bbC\backslash\bbR$.}
\no \\
& (iii) \; \ol{m_{+,\alpha}(z)}=m_{+,\alpha}(\ol z), \quad z\in\bbC_+.
\no \\
& (iv) \;\, \varepsilon |m_{+,\alpha}(\lambda+i\varepsilon)|\leq
C, \quad \lambda \in [\lambda_1,\lambda_2], \; 0<\varepsilon\leq
\varepsilon_0. \lb{2.38a} \\
& (v) \;\;\; \varepsilon |\Re(m_{+,\alpha}(\lambda+i\varepsilon))|
\underset{\varepsilon\downarrow 0}{=} \oh(1), \quad \lambda\in\bbR. \no
\\
& (vi) \;\;\, \rho_{+,\alpha}(\lambda)-\rho_{+,\alpha}(\lambda_0) =
\lim_{\delta\downarrow 0}\lim_{\varepsilon\downarrow 0} \f{1}{\pi}
\int_{\lambda_0+\delta}^{\lambda+\delta} d\mu \,
\Im(m_{+,\alpha}(\mu+i\varepsilon)). \no \\
& \qquad \; \text{defines a nondecreasing function $\rho_{+,\alpha}$ and
hence a measure on $\bbR$.} \no
\end{align}
Of course, properties $(ii)$--$(vi)$ are satisfied by any Herglotz
function. However, as we will see in Sections \ref{s3} and \ref{s4},
properties $(ii)$--$(vi)$ (possibly restricting $z$ to a sufficiently
small neighborhood of $\bbR$) are also crucial in connection with a class
of strongly singular potentials at $x=a$, where the analog of the
coefficient $m_{+,\alpha}$ will necessarily turn out to be a non-Herglotz
function. In particular, one can (and we will in Section
\ref{s3}) use an analog of \eqref{2.20} to prove items $(ii)$--$(vi)$ in
\eqref{2.38a} (for $|\Im(z)|$ sufficiently small) without ever invoking the
Herglotz property of $m_{+,\alpha}$, by just using the fact that the
left-hand side of \eqref{2.20} is a Herglotz function
whether or not the potential $V$ is strongly singular at the endpoint $a$.
Thus, the mere existence of the family of spectral projections
$\{E_{H_{+,\alpha}}(\lambda)\}_{\lambda\in\bbR}$ of the
self-adjoint operator $H_{+,\alpha}$ implies properties of the type
$(ii)$--$(vi)$.
\end{remark}

\begin{remark} \lb{r2.6a}
The effortless derivation of the link between the family of spectral
projections $E_{H_{+,\alpha}}(\cdot)$ and the spectral function
$\rho_{+,\alpha}(\cdot)$ of $H_{+,\alpha}$ in Theorem \ref{t2.5} applies
equally well to half-line Dirac-type operators and Hamiltonian systems (see
the extensive literature cited, e.g., in \cite{CG02}) and to half-lattice
Jacobi- (cf.\ \cite{Be68}) and CMV operators (i.e., semi-infinite five-diagonal unitary matrices which are related to orthogonal polynomials on the unit circle in the manner that half-lattice tri-diagonal (Jacobi) matrices are related to orthogonal polynomials on the real line as discussed in detail in \cite{Si05}; cf.\ \cite{GZ05} for an application of Theorem \ref{t2.5} to CMV operators). After
circulating a first draft of this manuscript, it was kindly pointed out to
us by Don Hinton that the idea of linking the family of spectral
projections and the spectral function using Stone's formula as the
starting point can already be found in a paper by Hinton and Schneider
\cite{HS98} published in 1998.
\end{remark}

Actually, one can improve on Theorem \ref{t2.5} and remove the
compact support restrictions on $f$ and $g$ in the usual way. To this
end one considers the map
\begin{equation}
\widetilde U_{+,\alpha}\colon \begin{cases} C_0^\infty((a,\infty))\to
L^2(\bbR; d\rho_{+,\alpha}) \\[1mm]
\hspace*{1.6cm} h \mapsto \hatt h_{+,\alpha}(\cdot)=
\int_a^\infty dx\, \phi_\alpha(\cdot,x) h(x). \end{cases} \lb{2.39}
\end{equation}
Taking $f=g$, $F=1$, $\lambda_1\downarrow -\infty$, and
$\lambda_2\uparrow \infty$ in \eqref{2.28} then shows that
$\widetilde U_{+,\alpha}$ is a densely defined isometry in
$L^2([a,\infty);dx)$, which extends by continuity to an isometry on
$L^2([a,\infty);dx)$. The latter is denoted by $U_{+,\alpha}$ and given by
\begin{equation}
U_{+,\alpha}\colon 
\begin{cases}
L^2([a,\infty);dx)\to L^2(\bbR;d\rho_{+,\alpha}) \\[1mm]
\hspace*{1.95cm}  h \mapsto \hatt h_{+,\alpha}(\cdot)=
\slimes_{b\uparrow\infty}\int_a^b dx\, \phi_\alpha(\cdot,x) h(x),
\end{cases}  \lb{2.40}
\end{equation}
where $\slimes$ refers to the $L^2(\bbR;d\rho_{+,\alpha})$-limit.

The calculation in \eqref{2.38} also yields
\begin{equation}
(E_{H_{+,\alpha}}((\lambda_1,\lambda_2])g)(x)
=\int_{(\lambda_1,\lambda_2]}
d\rho_{+,\alpha}(\lambda)\, \phi_\alpha(\lambda,x)
\hatt g_{+,\alpha}(\lambda), \quad g\in C_0^\infty((a,\infty))
\lb{2.41}
\end{equation}
and subsequently, \eqref{2.41} extends to all $g\in L^2([a,\infty);dx)$
by continuity. Moreover, taking $\lambda_1\downarrow -\infty$ and
$\lambda_2\uparrow \infty$ in \eqref{2.41} using
\begin{equation}
\slim_{\lambda\downarrow -\infty} E_{H_{+,\alpha}}(\lambda)=0, \quad
\slim_{\lambda\uparrow \infty}
E_{H_{+,\alpha}}(\lambda)=I_{L^2([a,\infty);dx)},
\lb{2.42}
\end{equation}
where
\begin{equation}
E_{H_{+,\alpha}}(\lambda)=E_{H_{+,\alpha}}((-\infty,\lambda]), \quad
\lambda\in\bbR, \lb{2.43}
\end{equation}
then yields
\begin{equation}
g(\cdot)=\slimes_{\mu_1\downarrow -\infty, \mu_2\uparrow\infty}
\int_{\mu_1}^{\mu_2} d\rho_{+,\alpha}(\lambda)\,
\phi_\alpha(\lambda,\cdot) \hatt g_{+,\alpha}(\lambda), \quad g\in
L^2([a,\infty);dx),  \lb{2.44}
\end{equation}
where $\slimes$ refers to the $L^2([a,\infty); dx)$-limit.

In addition, one can show
that the map $U_{+,\alpha}$ in \eqref{2.40} is onto and hence that
$U_{+,\alpha}$ is unitary (i.e., $U_{+,\alpha}$ and
$U_{+,\alpha}^{-1}$ are isometric isomorphisms between
$L^2([a,\infty);dx)$ and $L^2(\bbR;d\rho_{+,\alpha})$) with
\begin{equation}
U_{+,\alpha}^{-1} \colon \begin{cases} L^2(\bbR;d\rho_{+,\alpha}) \to
L^2([a,\infty);dx)  \\[1mm]
\hspace*{1.75cm}
\hatt h \mapsto \slimes_{\mu_1\downarrow -\infty, \mu_2\uparrow\infty}
\int_{\mu_1}^{\mu_2} d\rho_{+,\alpha}(\lambda)\,
\phi_\alpha(\lambda,\cdot) \hatt h(\lambda). \end{cases} \lb{2.45}
\end{equation}
Indeed, using \eqref{2.40} and \eqref{2.28} with $\la_1\downarrow-\infty$ and $\la_2\uparrow+\infty$, one has for all $f,g\in L^2([a,\infty);dx)$ and all bounded $F,G\in C(\bbR)$,
\begin{align}
(F(H_{+,\al})f,G(H_{+,\al})g)_{L^2([a,\infty);dx)} &= (f,F(H_{+,\al})^*G(H_{+,\al})g)_{L^2([a,\infty);dx)}   \no \\
&= (\hatt f_{+,\al},\ol{F}G\hatt g_{+,\al})_{L^2(\bbR;d\rho_{+,\alpha})}  \no \\
&= (F\hatt f_{+,\al},G\hatt g_{+,\al})_{L^2(\bbR;d\rho_{+,\alpha})}.   \lb{2.45a}
\end{align}
Then $U_{+,\al}(F(H_{+,\al})f)=F\hatt f_{+,\al}$ $d\rho_{\al,b}$-a.e.\ since for $h=F(H_{+,\al})f$ it follows from \eqref{2.45a} that
\begin{align}
& \|F\hatt f_{+,\al}-\hatt h_{+,\al}\|_{L^2(\bbR;d\rho_{+,\alpha})}^2 = (F\hatt f_{+,\al},F\hatt f_{+,\al})_{L^2(\bbR;d\rho_{+,\alpha})} - (F\hatt f_{+,\al},\hatt h_{+,\al})_{L^2(\bbR;d\rho_{+,\alpha})} \no
\\
&\qquad - (\hatt h_{+,\al},F\hatt f_{+,\al})_{L^2(\bbR;d\rho_{+,\alpha})} + (\hatt h_{+,\al},\hatt h_{+,\al})_{L^2(\bbR;d\rho_{+,\alpha})} \no
\\
& \quad = (F(H_{+,\al})f,F(H_{+,\al})f)_{L^2([a,\infty);dx)} - (F(H_{+,\al})f,h)_{L^2([a,\infty);dx)} \no
\\
&\qquad - (h,F(H_{+,\al})f)_{L^2([a,\infty);dx)} + (h,h)_{L^2([a,\infty);dx)}   \no \\
& \quad = 0. \lb{2.45b}
\end{align}
Next, one notes that for every $\la_0\in\bbR$, the range of $U_{+,\al}$ contains a continuous function $\hatt f_{+,\al}(\la)$ nonvanishing in a neighborhood of $\la_0$. For example, the image of
$f(\, \cdot \,)=\chi_{[c,d]} (\, \cdot \,) \, \phi_\al(\la_0, \, \cdot \,)$, $a<c<d<\infty$, has the above property. It thus follows from $U_{+,\al}(F(H_{+,\al})f)=F\hatt f_{+,\al}$ that the range of $U_{+,\al}$ contains all continuous functions and hence $U_{+,\al}$ is onto. 

We sum up these considerations in a variant of the
spectral theorem for (functions of) $H_{+,\alpha}$.

\begin{theorem} \lb{t2.6}
Let $\alpha\in [0,\pi)$ and $F\in C(\bbR)$. Then,
\begin{equation}
U_{+,\alpha} F(H_{+,\alpha})U_{+,\alpha}^{-1} = M_F \lb{2.46}
\end{equation}
in $L^2(\bbR;d\rho_{+,\alpha})$ $($cf.\ \eqref{2.30}$)$. Moreover,
\begin{align}
& \sigma(F(H_{+,\alpha}))= \essran_{d\rho_{+,\alpha}}(F), \lb{2.46a} \\
& \sigma(H_{+,\alpha})=\supp\,(d\rho_{+,\alpha}),  \lb{2.46b}
\end{align}
and the spectrum of $H_{+,\alpha}$ is simple.
\end{theorem}

Here the essential range of $F$ with respect to a measure
$d\mu$ is defined by
\begin{equation}
\essran_{d\mu}(F)=\{z\in\bbC\,|\, \text{for all
$\varepsilon>0$,} \, \mu(\{\lambda\in\bbR \,|\,
|F(\lambda)-z|<\varepsilon\})>0\}.   \lb{2.46c}
\end{equation}

We conclude the half-line case by recalling the following elementary
example of the  Fourier-sine transform.

\begin{example} \lb{e2.7}
Let $\alpha=0$ and $V(x)=0$ for a.e.\ $x\in(0,\infty)$. Then,
\begin{align}
& \phi_0(\lambda,x)=\frac{\sin(\lambda^{1/2}x)}{\lambda^{1/2}}, \quad
\lambda > 0, \; x\in (0,\infty), \no \\
& m_{+,0}(z)=iz^{1/2}, \quad  z\in\bbC\backslash [0,\infty), \lb{2.47} \\
& d\rho_{+,0}(\lambda)=\pi^{-1}\chi_{[0,\infty)}(\lambda)
\lambda^{1/2}d\lambda, \quad \lambda\in\bbR, \no
\end{align}
and hence,
\begin{align}
& \hatt h_{+,0}(\lambda)=\slimes_{y\uparrow\infty} \int_0^y dx\,
\frac{\sin(\lambda^{1/2}x)}{\lambda^{1/2}} h(x), \quad
h\in L^2([0,\infty);dx), \lb{2.48} \\
& h(x)=\slimes_{\mu\uparrow\infty} \frac{1}{\pi}
\int_{0}^{\mu} \lambda^{1/2}d\lambda \,
\frac{\sin(\lambda^{1/2}x)}{\lambda^{1/2}} \,
\hatt h_{+,0}(\lambda), \quad \hatt h_{+,0}
\in L^2([0,\infty); \pi^{-1}\lambda^{1/2}d\lambda).  \no
\end{align}
Introducing the change of variables
\begin{equation}
p=\lambda^{1/2}>0, \quad \hatt
H(p)=\bigg(\f{2\lambda}{\pi}\bigg)^{1/2}\hatt h_{+,0}(\lambda),
\lb{2.49}
\end{equation}
the pair of equations in \eqref{2.48} takes on the usual symmetric
form of the Fourier-sine transform,
\begin{align}
\begin{split}
                \hatt H(p)&=\slimes_{y\uparrow\infty}
\bigg(\frac{2}{\pi}\bigg)^{1/2} \int_0^y dx\,
\sin(px) h(x), \quad h\in L^2([0,\infty);dx), \lb{2.50} \\
                h(x)&=\slimes_{q \uparrow\infty}
\bigg(\frac{2}{\pi}\bigg)^{1/2}
\int_{0}^{q} dp \, \sin(px) \hatt H(p), \quad \hatt H
\in L^2([0,\infty); dp).
\end{split}
\end{align}
\end{example}

Next, we turn to the case of the entire real line and make the
following basic assumption.

\begin{hypothesis} \lb{h2.8}
$(i)$ Assume that
\begin{equation}
V\in L^1_{\loc} (\bbR;dx), \quad V \, \text{real-valued.}
\lb{2.51}
\end{equation}
$(ii)$ Introducing the differential expression $\tau$ given by
\begin{equation}
\tau=-\f{d^2}{dx^2} + V(x), \quad x\in\bbR, \lb{2.52}
\end{equation}
we assume $\tau$ to be in the limit point case at $+\infty$ and at
$-\infty$.
\end{hypothesis}

Associated with the differential expression $\tau$ one introduces the
self-adjoint Schr\"odinger operator $H$ in $L^2(\bbR;dx)$ by
\begin{align}
\begin{split}
&Hf=\tau f,   \\
&f\in \dom(H) = \big\{g\in L^2(\bbR;dx)\,\big|\, g, g' \in
AC_{\loc}(\bbR); \, \tau g\in L^2(\bbR; dx)\big\}. \lb{2.53}
\end{split}
\end{align}
Here $AC_{\loc}(\bbR)$ denotes the class of locally absolutely
continuous functions on $\bbR$.

As in the half-line context we introduce the usual fundamental
system of solutions $\phi_\alpha(z,\cdot,x_0)$ and
$\theta_\alpha(z,\cdot,x_0)$, $z\in\bbC$, of
\begin{equation}
(\tau \psi)(z,x) = z \psi(z,x), \quad x\in \bbR \lb{2.54}
\end{equation}
with respect to a fixed reference point $x_0\in\bbR$, satisfying the
initial conditions at the point $x=x_0$,
\begin{align}
\begin{split}
\phi_\alpha(z,x_0,x_0)&=-\theta'_\alpha(z,x_0,x_0)=-\sin(\alpha), \\
\phi'_\alpha(z,x_0,x_0)&=\theta_\alpha(z,x_0,x_0)=\cos(\alpha), \quad
\alpha \in [0,\pi). \lb{2.55}
\end{split}
\end{align}
Again we note that for any fixed $x, x_0\in \bbR$,
$\phi_\alpha(z,x,x_0)$ and $\theta_\alpha(z,x,x_0)$ are entire with
respect to $z$ and that
\begin{equation}
W(\theta_\alpha(z,\cdot,x_0),\phi_\alpha(z,\cdot,x_0))(x)=1, \quad
z\in\bbC.  \lb{2.56}
\end{equation}

Particularly important solutions of \eqref{2.54} are the
{\it Weyl--Titchmarsh solutions} $\psi_{\pm,\alpha}(z,\cdot,x_0)$,
$z\in\bbC\backslash\bbR$, uniquely characterized by
\begin{align}
\begin{split}
&\psi_{\pm,\alpha}(z,\cdot,x_0)\in L^2([x_0,\pm\infty);dx),  \\
&\sin(\alpha)\psi'_{\pm,\alpha}(z,x_0,x_0)
+\cos(\alpha)\psi_{\pm,\alpha}(z,x_0,x_0)=1, \quad
z\in\bbC\backslash\bbR. \lb{2.57}
\end{split}
\end{align}
The crucial condition in \eqref{2.57} is again the $L^2$-property which
uniquely determines $\psi_{\pm,\alpha}(z,\cdot,x_0)$ up to constant
multiples by the limit point hypothesis of $\tau$ at $\pm\infty$. In
particular, for
$\alpha, \beta \in [0,\pi)$,
\begin{align}
&\psi_{\pm,\alpha}(z,\cdot,x_0)=C_\pm(z,\alpha,\beta,x_0)
\psi_{\pm,\beta}(z,\cdot,x_0) \no \\
& \hspace*{.35cm}
\text{ for some coefficients $C_\pm (z,\alpha,\beta,x_0)\in\bbC$.}
\lb{2.58}
\end{align}
The normalization in \eqref{2.57} shows that
$\psi_{\pm,\alpha}(z,\cdot,x_0)$ are of the type
\begin{equation}
\psi_{\pm,\alpha}(z,x,x_0)=\theta_{\alpha}(z,x,x_0)
+m_{\pm,\alpha}(z,x_0)\phi_{\alpha}(z,x,x_0),
\quad  z\in\bbC\backslash\bbR, \; x\in\bbR \lb{2.59}
\end{equation}
for some coefficients $m_{\pm,\alpha}(z,x_0)$, the
{\it Weyl--Titchmarsh $m$-functions} associated with $\tau$, $\alpha$,
and $x_0$.

Again we recall the fundamental identity
\begin{align}
& \int_{x_0}^{\pm\infty} dx\,\psi_{\pm,\alpha}(z_{1},x,x_0)
\psi_{\pm,\alpha}(z_{2},x,x_0) = \frac{m_{\pm,\alpha}(z_{1},x_0)-
m_{\pm, \alpha}(z_{2},x_0)}{z_{1}-z_{2}}, \lb{2.60} \\
&\hspace*{7.75cm}  z_1, z_2 \in\bbC\backslash\bbR, \; z_1 \neq z_2, \no
\end{align}
and as before one concludes
\begin{equation}
\ol{m_{\pm,\alpha}(z,x_0)} = m_{\pm,\alpha}(\ol z,x_0), \quad
z\in\bbC\backslash\bbR.  \lb{2.61}
\end{equation}
Choosing $z_1=z$, $z_2=\ol z$ in \eqref{2.60}, one infers
\begin{equation}
\int_{x_0}^{\pm\infty} dx\,|\psi_{\pm,\alpha}(z,x,x_0)|^2
                               =
\frac{\Im(m_{\pm,\alpha}(z,x_0))}{\Im(z)}, \quad z
\in\bbC\backslash\bbR. \lb{2.62}
\end{equation}
Since $m_{\pm,\alpha}(\cdot,x_0)$ are analytic on $\bbC\backslash\bbR$,
$\pm m_{\pm,\alpha}(\cdot,x_0)$ are Herglotz functions.

The Green's function $G(z,x,x')$ of $H$ then reads
\begin{align}
G(z,x,x') & =
\f{1}{W(\psi_{+,\alpha}(z,\cdot,x_0),\psi_{-,\alpha}(z,\cdot,x_0))} \no \\
& \quad \times \begin{cases}
\psi_{-,\alpha}(z,x,x_0)\psi_{+,\alpha}(z,x',x_0), & x \leq x', \\
\psi_{-,\alpha}(z,x',x_0)\psi_{+,\alpha}(z,x,x_0), & x' \leq x,
\end{cases} \quad z\in\bbC\backslash\bbR,     \lb{2.63}
\end{align}
with
\begin{equation}
W(\psi_{+,\alpha}(z,\cdot,x_0),\psi_{-,\alpha}(z,\cdot,x_0))
=m_{-,\alpha}(z,x_0)-
m_{+,\alpha}(z,x_0), \quad z\in\bbC\backslash\bbR.  \lb{2.64}
\end{equation}
Thus,
\begin{equation}
((H-zI)^{-1}f)(x)
=\int_{\bbR} dx' \, G(z,x,x')f(x'), \quad z\in\bbC\backslash\bbR, \;
x\in\bbR, \; f\in L^{2}(\bbR;dx). \lb{2.65}
\end{equation}

Given $m_\pm(z,x_0)$, we also introduce the $2\times 2$ matrix-valued
Weyl--Titchmarsh function
\begin{align}
M_\alpha(z,x_0)&=\begin{pmatrix}
\f{1}{m_{-,\alpha}(z,x_0)-m_{+,\alpha}(z,x_0)} &
\f{1}{2}\f{m_{-,\alpha}(z,x_0)
+m_{+,\alpha}(z,x_0)}{m_{-,\alpha}(z,x_0)-m_{+,\alpha}(z,x_0)} \\
\f{1}{2}\f{m_{-,\alpha}(z,x_0)
+m_{+,\alpha}(z,x_0)}{m_{-,\alpha}(z,x_0)-m_{+,\alpha}(z,x_0)} &
\f{m_{-,\alpha}(z,x_0)m_{+,\alpha}(z,x_0)}{m_{-,\alpha}(z,x_0)
-m_{+,\alpha}(z,x_0)}\end{pmatrix}, \quad z\in\bbC\backslash\bbR.
\lb{2.71a}
\end{align}
$M_\alpha(z,x_0)$ is a Herglotz matrix with representation
\begin{align}
\begin{split}
& M_\alpha(z,x_0)=C_\alpha(x_0)+\int_{{\mathbb{R}}}
d\Omega_\alpha (\lambda,x_0)\bigg[\frac{1}{\lambda -z}-\frac{\lambda}
{1+\lambda^2}\bigg], \quad z\in\bbC\backslash\bbR, \lb{2.71b} \\
& C_\alpha(x_0)=C_\alpha(x_0)^*, \quad \int_{\bbR}
\f{\|d\Omega_{\alpha}(\lambda,x_0)\|}{1+\lambda^2} <\infty.
\end{split}
\end{align}
The Stieltjes inversion formula for the $2\times 2$ nonnegative
matrix-valued measure $d\Omega_\alpha(\cdot,x_0)$ then reads
\begin{equation}
\Omega_\alpha((\lambda_1,\lambda_2],x_0)
=\pi^{-1} \lim_{\delta\downarrow 0}
\lim_{\varepsilon\downarrow 0} \int^{\lambda_2+\delta}_{\lambda_1+\delta}
d\lambda \, \Im(M_\alpha(\lambda +i\varepsilon,x_0)), \quad \lambda_1,
\lambda_2 \in\bbR, \; \lambda_1<\lambda_2. \lb{2.71c}
\end{equation}
In particular, this implies that the entries
$d\Omega_{\al,\ell,\ell'}$, $\ell,\ell'=0,1$, of the matrix-valued
measure $d\Omega_{\al}$ are
real-valued scalar measures. Moreover, since the diagonal entries
of $M_\al$ are Herglotz functions, the diagonal entries of the
measure $d\Omega_{\al}$ are nonnegative measures. The
off-diagonal entries of the measure $d\Omega_{\al}$ equal a
complex measure which naturally admits a decomposition into a
linear conbination of differences of two nonnegative measures.

We note that in formulas \eqref{2.57}--\eqref{2.71b} one can
replace $z\in\bbC\backslash\bbR$ by $z\in\bbC\backslash\sigma(H)$.

\medskip

Next, we relate the family of spectral projections,
$\{E_H(\lambda)\}_{\lambda\in\bbR}$, of the self-adjoint
operator $H$ and the $2\times 2$ matrix-valued increasing spectral
function $\Omega_{\alpha}(\lambda,x_0)$, $\lambda\in\bbR$, which
generates the matrix-valued measure in the Herglotz representation
\eqref{2.71b} of $M_\alpha(z,x_0)$.

We first note that for $F\in C(\bbR)$,
\begin{align}
&(f,F(H)g)_{L^2(\bbR;dx)}= \int_{\bbR}
d (f,E_H(\lambda)g)_{L^2(\bbR;dx)}\,
F(\lambda), \lb{2.72} \\
& f, g \in\dom(F(H)) =\bigg\{h\in L^2(\bbR;dx) \,\bigg|\,
\int_{\bbR} d \|E_H(\lambda)h\|_{L^2(\bbR;dx)}^2 \, |F(\lambda)|^2
< \infty\bigg\}. \no
\end{align}

Given a $2\times 2$ matrix-valued nonnegative measure
$d\Omega=\big(d\Omega_{\ell,\ell'}\big)_{\ell,\ell'=0,1}$ on $\bbR$
with
\begin{equation}
d\Omega^{\rm tr} = d\Omega_{0,0} + d\Omega_{1,1}
\end{equation}
its trace measure, the density matrix
\begin{equation}
\bigg(\f{d\Omega_{\ell,\ell'}}{d\Omega^{\rm tr}}\bigg)_{\ell,\ell'=0,1}
\end{equation}
is locally integrable on $\bbR$ with respect to $d\Omega^{\rm tr}$. One
then introduces the vector-valued Hilbert space $L^2(\bbR;d\Omega)$ in the
following manner. Consider ordered pairs $f=(f_0,f_1)^\top$ of
$d\Omega^{\rm tr}$-measurable functions such that
\begin{equation}
\sum_{\ell,\ell'=0}^1 \ol{f_{\ell}(\cdot)}\,
\f{d\Omega_{\ell,\ell'}(\cdot)}{d\Omega^{\rm tr}(\cdot)}\,
f_{\ell'}(\cdot)
\end{equation}
is $d\Omega^{\rm tr}$-integrable on $\bbR$ and define
$L^2(\bbR;d\Omega)$ as the set of equivalence classes modulo
$d\Omega$-null functions. Here $g=(g_0,g_1)^{\top}\in
L^2(\bbR;d\Omega)$ is defined to be a $d\Omega$-null function if
\begin{equation}
\int_{\bbR} d\Omega^{\rm tr}(\lambda) \sum_{\ell,\ell'=0}^1
\ol{g_{\ell}(\lambda)} \,
\f{d\Omega_{\ell,\ell'}(\lambda)}{d\Omega^{\rm tr}(\lambda)} \,
g_{\ell'}(\lambda) =0.
\end{equation}
This space is complete with respect to the norm induced by the scalar
product
\begin{equation} \lb{2.72a}
(f,g)_{L^2(\bbR;d\Omega)} = \int_\bbR d\Omega^{\rm tr}(\lambda)
\sum_{\ell,\ell'=0}^1 \ol{f_{\ell}(\lambda)} \,
\f{d\Omega_{\ell,\ell'}(\lambda)}{d\Omega^{\rm tr}(\lambda)} \,
g_{\ell'}(\lambda), \quad f, g \in L^2(\bbR;d\Omega).
\end{equation}
For notational simplicity, expressions of the type \eqref{2.72a} will
usually be abbreviated by
\begin{equation} \lb{2.72b}
(f,g)_{L^2(\bbR;d\Omega)} = \int_\bbR \ol{f(\lambda)^{\top}} \,
d\Omega(\lambda) \, g(\lambda), \quad
f, g \in L^2(\bbR;d\Omega).
\end{equation}
(In this context we refer to \cite[p.\ 1345--1346]{DS88} for some
peculiarities in connection with matrix-valued nonnegative measures.)

\begin{theorem} \lb{t2.9}
Let $\alpha\in [0,\pi)$, $f,g \in C^\infty_0(\bbR)$,
$F\in C(\bbR)$, $x_0\in\bbR$, and $\lambda_1, \lambda_2 \in\bbR$,
$\lambda_1<\lambda_2$. Then,
\begin{align}
&(f,F(H)E_H((\lambda_1,\lambda_2])g)_{L^2(\bbR;dx)} =
\big(\hatt
f_{\alpha}(\cdot,x_0),M_FM_{\chi_{(\lambda_1,\lambda_2]}} \hatt
g_{\alpha}(\cdot,x_0)\big)_{L^2(\bbR;d\Omega_{\alpha}(\cdot,x_0))}
\no  \\
& \quad = \int_{(\lambda_1,\lambda_2]}
\ol{\hatt f_{\alpha}(\lambda,x_0)^\top} \,
d\Omega_{\alpha}(\lambda,x_0) \,
\hatt g_{\alpha}(\lambda,x_0)F(\lambda),  \lb{2.73}
\end{align}
where we introduced the notation
\begin{align}
\begin{split}
&\hatt h_{\alpha,0}(\lambda,x_0)=\int_\bbR dx \,
\theta_\alpha(\lambda,x,x_0) h(x),  \quad
\hatt h_{\alpha,1}(\lambda,x_0)=\int_\bbR dx \,
\phi_\alpha(\lambda,x,x_0) h(x), \lb{2.74} \\
&\hatt h_{\alpha}(\lambda,x_0)=\big(\hatt h_{\alpha,0}(\lambda,x_0),
\hatt h_{\alpha,1}(\lambda,x_0)\big)^\top,  \quad
\lambda \in\bbR, \;  h\in C^\infty_0(\bbR),
\end{split}
\end{align}
and $M_G$ denotes the maximally defined operator of multiplication
by the $d\Omega_{\alpha}^{\rm tr}$-measurable function $G$ in the
Hilbert space $L^2(\bbR;d\Omega_{\alpha}(\cdot,x_0))$,
\begin{align}
\begin{split}
& \big(M_G\hatt h\big)(\lambda)=G(\lambda)\hatt h(\lambda)
=\big(G(\lambda) \hatt h_0(\lambda), G(\lambda) \hatt h_1(\lambda)\big)^\top
\, \text{ for $d\Omega_{\alpha}^{\rm tr}$-a.e.\ $\lambda\in\bbR$}, \lb{2.75} \\
& \hatt h\in\dom(M_G)=\big\{\hatt k \in
L^2(\bbR;d\Omega_{\alpha}(\cdot,x_0)) \,\big|\,
G\hatt k \in L^2(\bbR;d\Omega_{\alpha}(\cdot,x_0))\big\}.
\end{split}
\end{align}
\end{theorem}
\begin{proof}
The point of departure for deriving \eqref{2.73} is again Stone's
formula \eqref{2.26a} applied to $T=H$,
\begin{align}
&(f,F(H)E_H((\lambda_1,\lambda_2]) g)_{L^2(\bbR;dx)}  \no \\
& \quad = \lim_{\delta\downarrow 0}\lim_{\varepsilon\downarrow 0}
\frac{1}{2\pi i} \int_{\lambda_1+\delta}^{\lambda_2+\delta}
d\lambda \, F(\lambda) \big[\big(f,(H-(\lambda+i\varepsilon)
I)^{-1}g\big)_{L^2(\bbR;dx)}  \no \\
& \hspace*{4.9cm} - \big(f,(H-(\lambda-i\varepsilon)
I)^{-1}g\big)_{L^2(\bbR;dx)}\big]. \lb{2.82}
\end{align}
Insertion of \eqref{2.63} and \eqref{2.65} into \eqref{2.82} then
yields the following:
\begin{align}
&(f,F(H)E_H((\lambda_1,\lambda_2]) g)_{L^2(\bbR;dx)} =
\lim_{\delta\downarrow 0}\lim_{\varepsilon\downarrow 0}
\frac{1}{2\pi i}
\int_{\lambda_1+\delta}^{\lambda_2+\delta} d\lambda \, F(\lambda) \no \\
& \quad \times \int_\bbR dx \bigg\{\f{1}{W(\lambda+i\varepsilon)}
\bigg[\ol{f(x)} \psi_{+,\alpha}(\lambda
+i\varepsilon,x,x_0) \int_{-\infty}^x dx'\,
\psi_{-,\alpha}(\lambda+i\varepsilon,x',x_0)g(x') \no \\
& \hspace*{1.6cm} + \ol{f(x)}
\psi_{-,\alpha}(\lambda+i\varepsilon,x,x_0)
\int_x^\infty dx'\, \psi_{+,\alpha}(\lambda+i\varepsilon,x',x_0) g(x')
\bigg]   \no \\
& \hspace*{1.55cm} -\f{1}{W(\lambda-i\varepsilon)}
\bigg[\ol{f(x)}
\psi_{+,\alpha}(\lambda -i\varepsilon,x,x_0) \int_{-\infty}^x dx'\,
\psi_{-,\alpha}(\lambda-i\varepsilon,x',x_0)g(x') \no \\
& \hspace*{1.6cm}  + \ol{f(x)}
\psi_{-,\alpha}(\lambda-i\varepsilon,x,x_0)
\int_x^\infty dx' \,\psi_{+,\alpha}(\lambda-i\varepsilon,x',x_0)
g(x')\bigg]\bigg\},  \lb{2.83}
\end{align}
where we used the abbreviation
\begin{equation}
W(z)=W(\psi_{+,\alpha}(z,\cdot,x_0)),\psi_{-,\alpha}(z,\cdot,x_0)),
\quad z\in\bbC\backslash\bbR. \lb{2.84}
\end{equation}
Freely interchanging the $dx$ and $dx'$ integrals with the limits and the
$d\lambda$ integral (since all integration domains are finite and all
integrands are continuous), and inserting the expressions \eqref{2.59} for
$\psi_{\pm,\alpha}(z,x,x_0)$ into \eqref{2.83}, one obtains
\begin{align}
&\big(f,F(H)E_H((\lambda_1,\lambda_2]) g\big)_{L^2(\bbR;dx)}
=\int_{\bbR} dx\, \ol{f(x)}
\bigg\{\int_{-\infty}^x dx' \, g(x') \no \\
& \quad \times \lim_{\delta\downarrow
0}\lim_{\varepsilon\downarrow 0} \frac{1}{2\pi i}
\int_{\lambda_1+\delta}^{\lambda_2+\delta} d\lambda \, F(\lambda)
\Big[\big[\theta_\alpha(\lambda,x,x_0) +
m_{+,\alpha}(\lambda+i\varepsilon,x_0)\phi_\alpha(\lambda,x,x_0)\big] \no \\
& \hspace*{2cm} \times \big[\theta_\alpha(\lambda,x',x_0) +
m_{-,\alpha}(\lambda+i\varepsilon,x_0)\phi_\alpha(\lambda,x',x_0)\big]
W(\lambda+i\varepsilon)^{-1}  \no \\
& \hspace*{1.7cm} -\big[\theta_\alpha(\lambda,x,x_0) +
m_{+,\alpha}(\lambda-i\varepsilon,x_0)\phi_\alpha(\lambda,x,x_0)\big] \no \\
& \hspace*{2cm} \times \big[\theta_\alpha(\lambda,x',x_0) +
m_{-,\alpha}(\lambda-i\varepsilon,x_0)\phi_\alpha(\lambda,x',x_0)\big]
W(\lambda-i\varepsilon)^{-1}\Big] \no \\
& \quad +\int_x^\infty dx'\, g(x') \lim_{\delta\downarrow 0}
\lim_{\varepsilon\downarrow 0} \frac{1}{2\pi i}
\int_{\lambda_1+\delta}^{\lambda_2+\delta} d\lambda \, F(\lambda)
\lb{2.85} \\
& \hspace*{2cm} \times \Big[\big[\theta_\alpha(\lambda,x,x_0) +
m_{-,\alpha}(\lambda+i\varepsilon,x_0)\phi_\alpha(\lambda,x,x_0)\big]
\no \\
& \hspace*{2cm} \times \big[\theta_\alpha(\lambda,x',x_0) +
m_{+,\alpha}(\lambda+i\varepsilon,x_0)\phi_\alpha(\lambda,x',x_0)\big]
W(\lambda+i\varepsilon)^{-1}   \no \\
& \hspace{1.7cm} -\big[\theta_\alpha(\lambda,x,x_0) +
m_{-,\alpha}(\lambda-i\varepsilon,x_0)\phi_\alpha(\lambda,x,x_0)\big]
\no \\
& \hspace*{2cm} \times\big[\theta_\alpha(\lambda,x',x_0) +
m_{+,\alpha}(\lambda-i\varepsilon,x_0)\phi_\alpha(\lambda,x',x_0)\big]
W(\lambda-i\varepsilon)^{-1}\Big]\bigg\}. \no
\end{align}
Here we employed the fact that for fixed $x\in\bbR$,
$\theta_\alpha(z,x,x_0)$ and $\phi_\alpha(z,x,x_0)$ are entire
with respect to $z$, that $\theta_\alpha(\la,x,x_0)$ and
$\phi_\alpha(\la,x,x_0)$ are real-valued for $\la\in\bbR$, that
$\phi_\alpha(z,\cdot,x_0), \theta_\alpha(z,\cdot,x_0)\in AC_{\loc}(\bbR)$,
and hence that
\begin{align}
\begin{split}
\theta_\alpha(\lambda\pm i\varepsilon,x,x_0)
&\underset{\varepsilon\downarrow 0}{=}
\theta_\alpha(\lambda,x,x_0) \pm
i\varepsilon (d/dz)\theta_\alpha(z,x,x_0)|_{z=\lambda} + \Oh(\ve^2),
\\
\phi_\alpha(\lambda\pm i\varepsilon,x,x_0)
&\underset{\varepsilon\downarrow 0}{=} \phi_\alpha(\lambda,x,x_0)
\pm i\varepsilon (d/dz)\phi_\alpha(z,x,x_0)|_{z=\lambda}
+ \Oh(\ve^2)
\end{split} \lb{2.85a}
\end{align}
with $\Oh(\varepsilon^2)$ being uniform with respect to
$(\lambda,x)$ as long as $\lambda$ and $x$ vary in compact subsets of
$\bbR^2$. Moreover, we used that
\begin{align}
&\varepsilon|M_{\al,\ell,\ell'}(\lambda+i\varepsilon,x_0)|\leq
C(\lambda_1,\lambda_2,\varepsilon_0,x_0), \quad \lambda\in
[\lambda_1,\lambda_2], \; 0<\varepsilon\leq\varepsilon_0, \;
\ell,\ell'=0,1, \no \\
&\varepsilon
|\Re(M_{\al,\ell,\ell'}(\lambda+i\varepsilon,x_0))|
\underset{\varepsilon\downarrow 0}{=}\oh(1), \quad \lambda\in\bbR, \;
\ell,\ell'=0,1, \lb{2.86}
\end{align}
which follows from the properties of Herglotz functions since
$M_{\al,\ell,\ell}$, $\ell=0,1$, are Herglotz and
$M_{\al,0,1}=M_{\al,1,0}$ have Herglotz-type representations by
decomposing the associated complex measure $d\Omega_{\al,0,1}$ into
$d\Omega_{\al,0,1}=d(\omega_1-\omega_2)+id(\omega_3-\omega_4)$, with
$d\omega_k$, $k=1,\dots,4$, nonnegative measures. In particular, utilizing
\eqref{2.61}, \eqref{2.85a},
\eqref{2.86}, and the elementary fact (cf.\ \eqref{2.64})
\begin{align}
\begin{split}
& \Im\bigg[\f{m_{\pm,\alpha}(\lambda+i\varepsilon,
x_0)}{W(\lambda+i\varepsilon)}\bigg]=\f{1}{2}
\Im\bigg[\f{m_{-,\alpha}(\lambda+i\varepsilon,
x_0)+m_{+,\alpha}(\lambda+i\varepsilon,
x_0)}{W(\lambda+i\varepsilon)}\bigg], \\
& \hspace*{8.5cm} \lambda\in\bbR, \; \varepsilon >0,  \lb{2.88}
\end{split}
\end{align}
$\phi_\alpha(\lambda\pm i\varepsilon,x,x_0)$ and
$\theta_\alpha(\lambda\pm i\varepsilon,x,x_0)$ under the $d\lambda$
integrals in \eqref{2.85} have immediately been replaced by
$\phi_\alpha(\lambda,x,x_0)$ and $\theta_\alpha(\lambda,x,x_0)$. Collecting
appropriate terms in
\eqref{2.85} then yields
\begin{align}
&(f,F(H)E_{H}((\lambda_1,\lambda_2]) g)_{L^2(\bbR;dx)} \no \\
& \quad =\int_\bbR dx\, \ol{f(x)}\int_\bbR dx' \, g(x')
\lim_{\delta\downarrow 0}\lim_{\varepsilon\downarrow 0}
\frac{1}{\pi} \int_{\lambda_1+\delta}^{\lambda_2+\delta} d\lambda
\, F(\lambda)
\no \\
& \qquad \times\bigg\{
\theta_\alpha(\lambda,x,x_0)\theta_\alpha(\lambda,x',x_0)
\Im\bigg[\f{1}{m_{-,\alpha}(\lambda+i\varepsilon,x_0)-
m_{+,\alpha}(\lambda+i\varepsilon,x_0)} \bigg] \no \\
& \hspace*{1.3cm}
+[\phi_\alpha(\lambda,x,x_0)\theta_\alpha(\lambda,x',x_0) +
\theta_\alpha(\lambda,x,x_0)\phi_\alpha(\lambda,x',x_0)] \no \\
& \hspace*{1.6cm}
\times \f{1}{2}\Im\bigg[\f{m_{-,\alpha}(\lambda+i\varepsilon,x_0)+
m_{+,\alpha}(\lambda+i\varepsilon,x_0)}{m_{-,\alpha}(\lambda
+i\varepsilon,x_0)-  m_{+,\alpha}(\lambda+i\varepsilon,x_0)}
\bigg] \lb{2.89} \\
& \hspace*{1.3cm}
+\phi_\alpha(\lambda,x,x_0)\phi_\alpha(\lambda,x',x_0)
\Im\bigg[\f{m_{-,\alpha}(\lambda+i\varepsilon,x_0)
m_{+,\alpha}(\lambda+i\varepsilon,x_0)}{m_{-,\alpha}(\lambda
+i\varepsilon,x_0)-  m_{+,\alpha}(\lambda+i\varepsilon,x_0)}
\bigg]\bigg\}. \no
\end{align}
Using the fact that by \eqref{2.71c} ($\ell, \ell'=0,1$)
\begin{align}
\begin{split}
&\int_{(\lambda_1,\lambda_2]} d\Omega_{\alpha,\ell,\ell'}(\lambda,x_0)
= \Omega_{\alpha,\ell,\ell'}((\lambda_1,\lambda_2],x_0) \\
& \quad =
\lim_{\delta\downarrow 0}\lim_{\varepsilon\downarrow 0}
\frac{1}{\pi}\int_{\lambda_1+\delta}^{\lambda_2+\delta} d\lambda \,
\Im(M_{\alpha,\ell,\ell'}(\lambda+i\varepsilon,x_0)), \lb{2.90}
\end{split}
\end{align}
and hence that
\begin{align}
&\int_{\bbR} d\Omega_{\alpha,\ell,\ell'}(\lambda,x_0)\, h(\lambda) =
\lim_{\varepsilon\downarrow 0} \frac{1}{\pi}\int_{\bbR} d\lambda \,
\Im(M_{\alpha,\ell,\ell'}(\lambda+i\varepsilon,x_0))\, h(\lambda),
\quad h\in C_0(\bbR),  \lb{2.91} \\
&\int_{(\lambda_1,\lambda_2]} d\Omega_{\alpha,\ell,\ell'}(\lambda,x_0)
\, k(\lambda) =
\lim_{\delta\downarrow 0} \lim_{\varepsilon\downarrow 0} \frac{1}{\pi}
\int_{\lambda_1+\delta}^{\lambda_2+\delta} d\lambda \,
\Im(M_{\alpha,\ell,\ell'}(\lambda+i\varepsilon,x_0))\, k(\lambda), \no
\\
& \hspace*{9.4cm} k\in C(\bbR), \lb{2.92}
\end{align}
one concludes
\begin{align}
&(f,F(H)E_H((\lambda_1,\lambda_2]) g)_{L^2(\bbR;dx)}
=\int_\bbR dx\, \ol{f(x)}\int_\bbR dx' \, g(x')
\int_{(\lambda_1,\lambda_2]} F(\lambda) \no \\
& \qquad \times
\Big\{\theta_\alpha(\lambda,x,x_0)\theta_\alpha(\lambda,x',x_0)
d\Omega_{\alpha,0,0}(\lambda,x_0)  \no \\
& \hspace*{1.25cm}
+[\phi_\alpha(\lambda,x,x_0)\theta_\alpha(\lambda,x',x_0)
+ \theta_\alpha(\lambda,x,x_0)\phi_\alpha(\lambda,x',x_0)]
d\Omega_{\alpha,0,1}(\lambda,x_0) \no \\
& \hspace*{1.25cm}
+\phi_\alpha(\lambda,x,x_0)\phi_\alpha(\lambda,x',x_0)
d\Omega_{\alpha,1,1}(\lambda,x_0)\Big\} \no \\
& \quad = \int_{(\lambda_1,\lambda_2]} \ol{\hatt
f_{\alpha}(\lambda,x_0)^\top}\, d\Omega_{\alpha}(\lambda,x_0) \,
\hatt g_{\alpha}(\lambda,x_0)F(\lambda), \lb{2.93}
\end{align}
using \eqref{2.74}, $d\Omega_{\alpha,0,1}(\cdot,x_0)=
d\Omega_{\alpha,1,0}(\cdot,x_0)$, and interchanging the $dx$, $dx'$ and
$d\Omega_{\alpha,\ell,\ell'}(\cdot,x_0)$, $\ell,\ell'=0,1$, integrals
once more.
\end{proof}

\begin{remark} \lb{r2.13a}
Again we emphasize that the idea of a straightforward derivation of the
link between the family of spectral projections $E_{H}(\cdot)$ and the
$2\times 2$ matrix-valued spectral function $\Omega_{\alpha}(\cdot)$ of
$H$ in Theorem \ref{t2.9} can already be found in \cite{HS98} as pointed
out in Remark \ref{r2.6a}. It applies equally well to Dirac-type operators
and Hamiltonian systems on $\bbR$ (see the extensive literature cited,
e.g., in \cite{CG02}) and to Jacobi and CMV operators on $\bbZ$
(cf.\ \cite{Be68} and \cite{GZ05}).
\end{remark}

As in the half-line case before, one can improve on Theorem \ref{t2.9}
and remove the compact support restrictions on $f$ and $g$ in the usual
way. To this end one considers the map
\begin{align}
&\widetilde U_{\alpha}(x_0)\colon \begin{cases} C_0^\infty(\bbR)\to
L^2(\bbR; d\Omega_{\alpha}(\cdot,x_0)) \\[1mm]
\hspace*{.9cm} h \mapsto \hatt h_{\alpha}(\cdot,x_0)
=\big(\hatt h_{\alpha,0}(\lambda,x_0),
\hatt h_{\alpha,1}(\lambda,x_0)\big)^\top, \end{cases} \lb{2.95} \\
&  \hatt h_{\alpha,0}(\lambda,x_0)=\int_\bbR dx \,
\theta_\alpha(\lambda,x,x_0) h(x),  \quad
                               \hatt h_{\alpha,1}(\lambda,x_0)=\int_\bbR dx \,
\phi_\alpha(\lambda,x,x_0) h(x).  \no
\end{align}
Taking $f=g$, $F=1$, $\lambda_1\downarrow -\infty$, and
$\lambda_2\uparrow \infty$ in \eqref{2.73} then shows that $\widetilde
U_{\alpha}(x_0)$ is a densely defined isometry in $L^2(\bbR;dx)$,
which extends by continuity to an isometry on $L^2(\bbR;dx)$. The latter
is denoted by $U_{\alpha}(x_0)$ and given by
\begin{align}
&U_{\alpha}(x_0)\colon \begin{cases} L^2 (\bbR;dx)\to
L^2(\bbR; d\Omega_{\alpha}(\cdot,x_0)) \\[1mm]
\hspace*{1.25cm} h \mapsto \hatt h_{\alpha}(\cdot,x_0)
=\big(\hatt h_{\alpha,0}(\cdot,x_0),
\hatt h_{\alpha,1}(\cdot,x_0)\big)^\top, \end{cases} \lb{2.96} \\
& \hatt h_\alpha(\cdot,x_0)=\begin{pmatrix} \hatt h_{\alpha,0}(\cdot,x_0) \\
\hatt h_{\alpha,1}(\cdot,x_0) \end{pmatrix}=
\slimes_{a\downarrow -\infty, b \uparrow\infty} \begin{pmatrix}
\int_{a}^b dx \, \theta_\alpha(\cdot,x,x_0) h(x) \\
\int_{a}^b dx \, \phi_\alpha(\cdot,x,x_0) h(x) \end{pmatrix}, \no
\end{align}
where $\slimes$ refers to the
$L^2(\bbR; d\Omega_{\alpha}(\cdot,x_0))$-limit.

The calculation in \eqref{2.93} also yields
\begin{align}
&(E_{H}((\lambda_1,\lambda_2])g)(x) =\int_{(\lambda_1,\lambda_2]}
(\theta_\alpha(\lambda,x,x_0), \phi_\alpha(\lambda,x,x_0)) \,
d\Omega_{\alpha}(\lambda,x_0)\,
\hatt g_{\alpha}(\lambda,x_0) \no \\
&\quad =\int_{(\lambda_1,\lambda_2]}
\Big\{d\Omega_{\alpha,0,0}(\lambda,x_0)
\, \theta_\alpha(\lambda,x,x_0) \hatt g_{\alpha,0}(\lambda,x_0)
\no \\
&\qquad \quad + d\Omega_{\alpha,0,1}(\lambda,x_0) \,
[\theta_\alpha(\lambda,x,x_0) \hatt g_{\alpha,1}(\lambda,x_0) +
\phi_\alpha(\lambda,x,x_0) \hatt g_{\alpha,0}(\lambda,x_0)] \no \\
& \qquad \quad +d\Omega_{\alpha,1,1}(\lambda,x_0) \,
\phi_\alpha(\lambda,x,x_0) \hatt g_{\alpha,1}(\lambda,x_0)
\Big\}, \quad g\in C_0^\infty(\bbR) \lb{2.97}
\end{align}
and subsequently, \eqref{2.97} extends to all $g\in L^2(\bbR;dx)$
by continuity. Moreover, taking $\lambda_1\downarrow -\infty$ and
$\lambda_2\uparrow \infty$ in \eqref{2.97} and using
\begin{equation}
\slim_{\lambda\downarrow -\infty} E_{H}(\lambda)=0, \quad
\slim_{\lambda\uparrow \infty} E_{H}(\lambda)=I_{L^2(\bbR;dx)},
\lb{2.98}
\end{equation}
where
\begin{equation}
E_{H}(\lambda)=E_{H}((-\infty,\lambda]), \quad \lambda\in\bbR, \lb{2.99}
\end{equation}
then yield
\begin{align}
g(\cdot) &=\slimes_{\mu_1\downarrow -\infty, \mu_2\uparrow \infty}
\int_{(\mu_1,\mu_2]} (\theta_\alpha(\lambda,\cdot,x_0),
\phi_\alpha(\lambda,\cdot,x_0)) \, d\Omega_{\alpha}(\lambda,x_0)\,
\hatt g_{\alpha}(\lambda,x_0) \no \\
&= \slimes_{\mu_1\downarrow -\infty, \mu_2\uparrow \infty}
\int_{\mu_1}^{\mu_2}
\Big\{d\Omega_{\alpha,0,0}(\lambda,x_0)
\, \theta_\alpha(\lambda,\cdot,x_0) \hatt g_{\alpha,0}(\lambda,x_0)
\no \\
&\qquad +
d\Omega_{\alpha,0,1}(\lambda,x_0) \,
[\theta_\alpha(\lambda,\cdot,x_0) \hatt g_{\alpha,1}(\lambda,x_0) +
\phi_\alpha(\lambda,\cdot,x_0) \hatt g_{\alpha,0}(\lambda,x_0)] \no \\
& \qquad +d\Omega_{\alpha,1,1}(\lambda,x_0) \,
\phi_\alpha(\lambda,\cdot,x_0) \hatt g_{\alpha,1}(\lambda,x_0)
\Big\}, \quad g\in L^2(\bbR;dx), \lb{2.100}
\end{align}
where $\slimes$ refers to the $L^2(\bbR;dx)$-limit. In addition, one can show
that the map $U_{\alpha}(x_0)$ in \eqref{2.96} is onto and hence that
$U_{\alpha}(x_0)$ is unitary with
\begin{align}
&U_{\alpha}(x_0)^{-1} \colon \begin{cases}
L^2(\bbR;d\Omega_{\alpha}(\cdot,x_0)) \to  L^2(\bbR;dx)  \\[1mm]
\hspace*{2.4cm} \hatt h \mapsto  h_\alpha, \end{cases} \lb{2.101} \\
& h_\alpha(\cdot)= \slimes_{\mu_1\downarrow -\infty, \mu_2\uparrow
\infty} \int_{\mu_1}^{\mu_2}(\theta_\alpha(\lambda,\cdot,x_0),
\phi_\alpha(\lambda,\cdot,x_0)) \, d\Omega_{\alpha}(\lambda,x_0)\,
\hatt h(\lambda). \no
\end{align}
Indeed, following the argument in \eqref{2.45a}--\eqref{2.45b}, one obtains $U_\al(F(H)f)=FU_\al(f)$ for all $f\in L^2(\bbR;dx)$ and all bounded $F\in C(\bbR)$.
Since $U_\al$ is an isometry, the range of $U_\al$ is closed and hence $U_\al$ is onto if the only function orthogonal to the range of $U_\al$ is the zero function. Let $F\in C_0(\bbR)$ and suppose $\hatt f \in L^2(\bbR;d\Omega_{\alpha}(\, \cdot \,,x_0))$ is orthogonal to the range of $U_\al$ then, in particular, $\hatt f$ is orthogonal to $U_\al(F(H)\chi_{[x_0,y]})=FU_\al(\chi_{[x_0,y]})$ for every $y\in(a,b)$, that is,
\begin{align}
\int_\bbR F(\la) \bigg(\int_{x_0}^y dx\,\theta_\al(\la,x,x_0),
\int_{x_0}^y dx\,\phi_\al(\la,x,x_0)\bigg)d\Omega_{\al}(\la,x_0)\hatt f(\la) = 0. \lb{6.101a}
\end{align}
Differentiating twice with respect to $y$ and taking $y=x_0$ then yields
\begin{align}
&\int_\bbR F(\la) \big(\theta_\al(\la,x_0,x_0),\phi_\al(\la,x_0,x_0)\big) \, d\Omega_{\al}(\la,x_0)\hatt f(\la) \no
\\
&\quad = \int_\bbR F(\la) \big(\cos(\al),-\sin(\al)\big) \,
d\Omega_{\al}(\la,x_0) \, \hatt f(\la) = 0, \lb{6.101b}
\\
&\int_\bbR F(\la) \big(\theta^{[1]}_\al(\la,x_0,x_0),\phi^{[1]}_\al(\la,x_0,x_0)\big) \, d\Omega_{\al}(\la,x_0)\hatt f(\la) \no
\\
&\quad = \int_\bbR F(\la) \big(\sin(\al),\cos(\al)\big) \,
d\Omega_{\al}(\la,x_0) \, \hatt f(\la) = 0. \lb{6.101c}
\end{align}
Taking linear combinations of \eqref{6.101b} and \eqref{6.101c} then implies
\begin{align}
\int_\bbR \big(F(\la),0\big) \, d\Omega_{\al}(\la,x_0) \, \hatt f(\la) =
\int_\bbR \big(0,F(\la)\big) \, d\Omega_{\al}(\la,x_0) \, \hatt f(\la) = 0. \lb{6.101d}
\end{align}
Since the vector functions of the form $\big(F(\la),0\big)^\top$, $\big(0,F(\la)\big)^\top$, $F\in C_0(\bbR)$, are dense in $L^2(\bbR;d\Omega_{\alpha}(\, \cdot \,,x_0))$, \eqref{6.101d} implies that $\hatt f=0$. Thus, $U_\al$ is onto.

We sum up these considerations in a variant of the spectral theorem for
(functions of) $H$.

\begin{theorem} \lb{t2.10}
Let $F\in C(\bbR)$ and $x_0\in\bbR$. Then,
\begin{equation}
U_{\alpha}(x_0) F(H)U_{\alpha}(x_0)^{-1} = M_F \lb{2.102}
\end{equation}
in $L^2(\bbR;d\Omega_{\alpha}(\cdot,x_0))$ $($cf.\ \eqref{2.75}$)$.
Moreover,
\begin{equation}
\sigma(H)=\supp\,(d\Omega_\alpha(\cdot,x_0))=\supp\,
(d\Omega_\alpha^{\rm tr}(\cdot,x_0)).
\end{equation}
\end{theorem}
Here $d\Omega_\alpha^{\rm tr}(\cdot,x_0) =
d\Omega_{\alpha,0,0}(\cdot,x_0) + d\Omega_{\alpha,1,1}(\cdot,x_0)$
denotes the trace measure of $d\Omega_\alpha(\cdot,x_0)$.

We conclude the case of the entire line with an elementary example.

\begin{example} \lb{e2.15}
Let $\alpha=0$, $x_0=0$ and $V(x)=0$ for a.e.\ $x\in\bbR$. Then,
\begin{align}
& \phi_0(\lambda,x,0)=\frac{\sin(\lambda^{1/2}x)}{\lambda^{1/2}}, \quad
\theta_0(\lambda,x,0)=\cos(\lambda^{1/2}x),  \quad
\lambda > 0, \; x\in\bbR, \no \\
& m_{\pm, 0}(z,0)=\pm i z^{1/2}, \quad z\in \bbC\backslash [0,\infty),
\lb{2.124} \\
& d\Omega_{0}(\lambda,0)=\f{1}{2 \pi}\chi_{(0,\infty)}(\lambda)
\begin{pmatrix}\lambda^{-1/2} & 0 \\ 0 & \lambda^{1/2}
\end{pmatrix}d\lambda, \quad \lambda\in\bbR,  \no
\end{align}
and hence,
\begin{align}
& \hatt h_0(\lambda,0)=\begin{pmatrix} \hatt h_{0,0}(\lambda,0) \\
\hatt h_{0,1}(\lambda,0) \end{pmatrix}
=\slimes_{a\downarrow -\infty,b\uparrow\infty}
\begin{pmatrix} \int_a^b dx\, \cos(\lambda^{1/2}x)h(x)  \\
\int_a^b dx\,\lambda^{-1/2}\sin(\lambda^{1/2}x) h(x)\end{pmatrix}, \no \\
& \hspace*{8.2cm}  h\in L^2(\bbR;dx), \no \\
& h(x)=\slimes_{\mu\uparrow\infty} \frac{1}{2 \pi}
\int_{0}^{\mu} \lambda^{1/2}d\lambda \, \bigg[
\f{\cos(\lambda^{1/2}x)}{\lambda} \hatt h_{0,0}(\lambda,0) +
\frac{\sin(\lambda^{1/2}x)}{\lambda^{1/2}} \,
\hatt h_{0,1}(\lambda,0)\bigg],   \no \\
&\hspace*{6.3cm}  \hatt h_0(\cdot,0) \in L^2([0,\infty); d\Omega_0(\cdot,0)).   \lb{2.125}
\end{align}
Introducing the change of variables
\begin{equation}
p=\lambda^{1/2}>0, \quad \hatt
H(p)=\begin{pmatrix} \hatt H_0(p) \\ \hatt H_1(p) \end{pmatrix} =
\f{1}{\pi^{1/2}} \begin{pmatrix}
\hatt h_{0,0} (\lambda,0)  \\ \lambda^{1/2} \hatt h_{0,1}(\lambda,0)
\end{pmatrix},  \lb{2.126}
\end{equation}
the pair of equations in \eqref{2.125} take on the symmetric
form,
\begin{align}
& \hatt H(p)
=\slimes_{a\downarrow -\infty,b\uparrow\infty} \f{1}{\pi^{1/2}}
\begin{pmatrix} \int_a^b dx\, \cos(px)h(x)  \\
\int_a^b dx\,\sin(px) h(x)\end{pmatrix}, \quad  h\in L^2(\bbR;dx), \no \\
& h(x)=\slimes_{q\uparrow\infty} \frac{1}{\pi^{1/2}}
\int_{0}^{q} dp \, \big[\cos(px) \hatt H_{0}(p) +
\sin(px) \, \hatt H_{1}(p)\big],  \lb{2.127} \\
&\hspace*{4.82cm}  \hatt H_\ell \in L^2([0,\infty); dp), \; \ell=0,1.  \no
\end{align}
One verifies that the pair of equations in \eqref{2.127} is equivalent to the
usual Fourier transform
\begin{align}
\begin{split}
& \wti h(p) =\slimes_{y\uparrow\infty}
\f{1}{(2\pi)^{1/2}}
\int_{-y}^{y} dx\, e^{ipx} h(x), \quad h \in L^2 (\bbR; dx), \lb{2.128} \\
& h(x) =\slimes_{q\uparrow\infty} \f{1}{(2\pi)^{1/2}}
\int_{-q}^{q} dp\, e^{-ipx} \wti h(p), \quad \wti h \in L^2(\bbR; dq).
\end{split}
\end{align}
\end{example}

\section{The Case of Strongly Singular Potentials} \lb{s3}

In this section we extend our discussion to a class of strongly
singular potentials $V$ on the half-line $(a,\infty)$ with the singularity
of $V$ being concentrated at the endpoint $a$. We will present and contrast
two approaches to this problem: One in which the reference point
$x_0$ coincides with the singular endpoint $a$ leading to a (scalar)
spectral function, and one in which $x_0$ lies in the interior of the
half-line $(a,\infty)$ and hence is a regular point for the half-line
Schr\"odinger differential expression. The latter case naturally leads to a
$2\times 2$ matrix-valued spectral function which will be shown to be
essentially equivalent to the scalar spectral function obtained from the
former approach. While Herglotz functions still lie at the heart of the
matter of spectral functions (resp., matrices), the direct analog of
half-line Weyl--Titchmarsh coefficients will cease to be Herglotz
functions in the first approach where the reference point $x_0$ coincides
with the endpoint $a$.

\begin{hypothesis} \lb{h4.1}
$(i)$ Let $a\in\bbR$ and assume that
\begin{equation}
V\in L^1_{\loc} ((a,\infty);dx), \quad V \, \text{real-valued.} \lb{4.1}
\end{equation}
$(ii)$ Introducing the differential expression $\tau_+$ given by
\begin{equation}
\tau_+=-\f{d^2}{dx^2} + V(x), \quad x\in(a,\infty), \lb{4.2}
\end{equation}
we assume $\tau_+$ to be in the limit point case at $a$ and at
$+\infty$. \\
$(iii)$ Assume there exists an analytic Weyl--Titchmarsh
solution $\wti\phi(z,\cdot)$ of
\begin{align}
(\tau_+ \psi)(z,x) = z \psi(z,x), \quad x\in(a,\infty), \lb{4.3}
\end{align}
for $z$ in an open neighborhood $\cO$ of $\bbR$ $($containing $\bbR$$)$
in the following sense: \\
\indent $(\alpha)$ For all $x\in(a,\infty)$, $\wti\phi(z,x)$ is
analytic with respect to $z\in\cO$. \\
\indent $(\beta)$ $\wti\phi(z,x)$, $x\in\bbR$, is real-valued for
$z\in\bbR$. \\
\indent $(\gamma)$ $\wti\phi(z,\cdot)$ satisfies an $L^2$-condition near
the end point $a$
\begin{align}
\int_{a}^{b} dx \, \big|\wti\phi(z,x)\big|^2 < \infty \,
\text{ for all $b\in(a,\infty)$}
\end{align}
\indent \quad\;\; for all $z\in\bbC\backslash\bbR$ with $|\Im(z)|$
sufficiently small.
\end{hypothesis}

Without loss of generality we assumed in Hypothesis \ref{h4.1}\,$(iii)$
that the analytic Weyl--Titchmarsh solution satisfies the $L^2$-condition
near the left end point $a$. One can replace this by the analogous
$L^2$-condition at $\infty$.

A class of examples of strongly singular potentials satisfying
Hypothesis \ref{h4.1} will be discussed in Examples \ref{e4.11}
and \ref{e4.14} at the end of this section.

While we focus on strongly singular potentials with $\tau_+$ in the limit
point case at both endpoints $a$ and $\infty$, the case of strongly
singular potentials with $\tau_+$ in the limit circle case at both
endpoints has been studied by Fulton \cite{Fu77}.

Associated with the differential expression $\tau_+$ one
introduces the self-adjoint Schr\"odinger operator $H_+$ in
$L^2([a,\infty);dx)$ by
\begin{align}
&H_+f=\tau_+ f,  \lb{4.5} \\
&f\in \dom(H_+)=\{g\in L^2([a,\infty);dx)\,|\, g, g' \in
AC_{\loc}((a,\infty)); \, \tau_+ g\in L^2([a,\infty); dx)\}.  \no
\end{align}

Next, we introduce the usual fundamental system of solutions
$\phi(z,\cdot,x_0)$ and $\theta(z,\cdot,x_0)$, $z\in\bbC$, of
\eqref{4.3} satisfying the initial conditions at the fixed reference point
$x_0\in(a,\infty)$,
\begin{equation}
\phi(z,x_0,x_0)=\theta'(z,x_0,x_0)=0, \quad
\phi'(z,x_0,x_0)=\theta(z,x_0,x_0)=1. \lb{4.7}
\end{equation}
Thus, for any fixed $x\in(a,\infty)$, the solutions
$\phi(z,x,x_0)$ and $\theta(z,x,x_0)$ are entire with respect to
$z$ and
\begin{equation}
W(\theta(z,\cdot,x_0),\phi(z,\cdot,x_0))(x)=1, \quad z\in\bbC.
\lb{4.8}
\end{equation}

We note, that Hypothesis \ref{h4.1}\,$(iii)$ implies that for fixed
$x\in(a,\infty)$, $\wti\phi'(z,x)$ is also analytic with respect to
$z\in\cO$. This follows from differentiating the identity
\begin{equation}
\wti\phi(z,x)=\wti\phi'(z,x_0)\phi(z,x,x_0)+\wti\phi(z,x_0)\theta(z,x,x_0),
\quad x, x_0 \in (a,\infty)  \lb{4.9}
\end{equation}
for $z\in\cO$. More precisely, one can argue as follows: One considers Volterra integral equations of the type
\begin{align}
& \psi_j(z,x,x_0) = \psi_j(z_0,x,x_0)   \no \\
& \quad + (z-z_0) \int_{x_0}^x dx' \,
\f{\psi_1(z_0,x,x_0) \psi_2(z_0,x',x_0) - \psi_1(z_0,x',x_0) \psi_2(z_0,x,x_0)}{W
(\psi_1(z_0,\cdot,x_0),\psi_2(z_0,\cdot,x_0))} \no \\
& \hspace*{3.1cm} \times \psi_j(z,x',x_0), \quad z, z_0 \in\cO, \; j=1,2,
\end{align}
where $\psi_j(z_0,\cdot,x_0)$, $j=1,2$, is a fundamental system of solutions of
\eqref{4.3} for $z = z_0$ (such as $\phi$ and $\theta$ in \eqref{4.7}). In particular,
\begin{equation}
\psi_j(z_0,x_0,x_0) = \alpha_j \in\bbC, \quad
\psi_j^\prime (z_0,x_0,x_0) = \beta_j \in\bbC, \quad z_0 \in \cO, \; j=1,2,
\end{equation}
implies
\begin{equation}
\psi_j(z,x_0,x_0) = \alpha_j \in\bbC, \quad
\psi_j^\prime (z,x_0,x_0) = \beta_j \in\bbC,  \quad z \in \cO, \; j=1,2,
\end{equation}
Analyticity of $\psi_j^\prime (z,\cdot,x_0)$, $j=1,2$, for $z\in\cO$ then follows from
the equation
\begin{align}
& \psi_j^\prime (z,x,x_0) = \psi_j^\prime (z_0,x,x_0)   \no \\
& \quad + (z-z_0) \int_{x_0}^x dx' \,
\f{\psi_1^\prime (z_0,x,x_0) \psi_2(z_0,x',x_0)
- \psi_1(z_0,x',x_0) \psi_2^\prime (z_0,x,x_0)}{W
(\psi_1(z_0,\cdot,x_0),\psi_2(z_0,\cdot,x_0))} \no \\
& \hspace*{3.1cm} \times \psi_j(z,x',x_0), \quad z, z_0 \in\cO, \; j=1,2.
\end{align}

Next, we also introduce the {\it Weyl--Titchmarsh solutions}
$\psi_{\pm}(z,\cdot,x_0)$, $x_0\in(a,\infty)$, $z\in\bbC\backslash\bbR$
of \eqref{4.3}. Since by Hypothesis \ref{h4.1} $(ii)$, $\tau_+$ is
assumed to be in the limit point case at $a$ and at $\infty$, the
Weyl--Titchmarsh solutions are uniquely characterized (up to constant
multiples) by
\begin{equation}
\psi_{-}(z,\cdot,x_0)\in L^2([a,x_0];dx), \quad
\psi_{+}(z,\cdot,x_0)\in L^2([x_0,\infty);dx), \quad
z\in\bbC\backslash\bbR. \lb{4.10}
\end{equation}
We fix the normalization of $\psi_\pm(z,\cdot,x_0)$ by requiring
$\psi_{\pm}(z,x_0,x_0)=1$ and hence $\psi_\pm(z,\cdot,x_0)$ have the
following structure,
\begin{align}
\psi_\pm(z,x,x_0) = \theta(z,x,x_0) + m_\pm(z,x_0)\phi(z,x,x_0),
\quad x,x_0 \in(a,\infty), \; z\in\bbC\backslash\bbR, \lb{4.11}
\end{align}
where the coefficients $m_{\pm}(z,x_0)$ are given by
\begin{align}
m_\pm(z,x) = \f{\psi'_\pm(z,x,x_0)}{\psi_\pm(z,x,x_0)}, \quad
x,x_0 \in(a,\infty), \; z\in\bbC\backslash\bbR, \lb{4.12}
\end{align}
and are Herglotz and anti-Herglotz functions, respectively.

\begin{lemma} \lb{l4.2}
Assume Hypothesis \ref{h4.1}\,$(i)$ and $(ii)$.\ Then Hypothesis
\ref{h4.1}\,$(iii)$ is equivalent to the assumption that for any fixed
$x\in(a,\infty)$, $m_-(z,x)$ is meromorphic with respect to $z\in\bbC$.
\end{lemma}
\begin{proof}
In the following we fix $x\in (a,\infty)$.
First, assume Hypothesis \ref{h4.1}. By Hypothesis \ref{h4.1}\,$(ii)$, the
Weyl--Titchmarsh solutions are unique up to constant multiples and one
concludes that $\psi_-(z,\cdot,x_0) = c(z,x_0)\wti\phi(z,\cdot)$. Hence
by \eqref{4.12},
\begin{align}
m_-(z,x) = \f{\wti\phi'(z,x)}{\wti\phi(z,x)}, \quad x\in (a,\infty), \;
z\in\bbC\backslash\bbR.  \lb{4.13}
\end{align}
Since by Hypothesis \ref{h4.1}\,$(iii)$, $\wti\phi(z,x)$ and
$\wti\phi'(z,x)$ are analytic with respect to $z\in\cO$ (cf.\ the
paragraph preceding \eqref{4.9}), one concludes that $m_-(z,x$) is
meromorphic in $z\in\cO$ and since $m_-$ is analytic in
$\bbC\backslash\bbR$, $m_-$ is meromorphic on $\bbC$.

Conversely, if $m_-(z,x)$ is meromorphic with respect to $z\in\bbC$,
then it has the following structure,
\begin{align}
m_-(z,x) = \f{\eta_1(z,x)}{\eta_2(z,x)},
\end{align}
where $\eta_1(z,x)$ and $\eta_2(z,x)$ can be chosen to be entire such
that they do not have common zeros. Moreover, since the zeros of
$\eta_j(\cdot,x)$, $j=1,2$, are necessarily all real, the Weierstrass
factorization theorem (cf., e.g., Corollary 2 of Theorem II.10.1 in
\cite[p.\ 284--285]{Ma85}) shows that $\eta_1(z,x)$ and $\eta_2(z,x)$
can be chosen to be real for $z\in\bbR$. Thus, for $x_0\in (a,\infty)$,
\begin{align}
\wti\phi(z,\cdot) = \eta_2(z,x_0)\psi_-(z,\cdot,x_0) =
\eta_2(z,x_0)\theta(z,\cdot,x_0)+\eta_1(z,x_0)\phi(z,\cdot,x_0)
\end{align}
is entire in $z$, and moreover, it is a Weyl--Titchmarsh solution of
\eqref{4.3} that satisfies Hypothesis \ref{h4.1}\,$(iii)$.
\end{proof}

\begin{lemma} \lb{l4.4}
Assume Hypothesis \ref{h4.1} $(iii)$. Then, there is an open
neighborhood $\cO'$ of $\bbR$ $($containing $\bbR$$)$,
$\cO'\subseteq\cO$, and a solution $\wti\theta(z,\cdot)$ of \eqref{4.3},
which, for each $x\in(a,\infty)$, is analytic with respect to
$z\in\cO'$, real-valued for $z\in\bbR$, such that,
\begin{align} \lb{4.18}
W(\wti\theta(z,\cdot),\wti\phi(z,\cdot))(x)=1, \quad z\in\cO'.
\end{align}
\end{lemma}
\begin{proof}
Let $x_0\in(a,\infty)$ and consider the following solution of \eqref{4.3},
\begin{align}
&\wti\theta(z,x) =
\f{\wti\phi'(z,x_0)}{\wti\phi(z,x_0)^2+\wti\phi'(z,x_0)^2}\theta(z,x,x_0)
-
\f{\wti\phi(z,x_0)}{\wti\phi(z,x_0)^2+\wti\phi'(z,x_0)^2}\phi(z,x,x_0),
\no \\
&\hspace*{9cm} x\in(a,\infty) \lb{4.19}
\end{align}
for $z$ in a sufficiently small neighborhood of $\bbR$. Since for
$x,x_0\in(a,\infty)$, $\wti\phi(z,x)$, $\theta(z,x,x_0)$, and
$\phi(z,x,x_0)$ are analytic with respect to $z\in\cO$, real-valued for
$z\in\bbR$, and $\wti\phi(z,x_0)$, $\wti\phi'(z,x_0)$ are not both zero
for all $z$ in a sufficiently small neighborhood of $\bbR$,
$\wti\theta(z,x)$ in \eqref{4.19} is analytic with respect to
$z\in\cO'$, $\cO'\subseteq\cO$, for fixed $x\in (a,\infty)$ and
real-valued for $z\in\bbR$. Moreover, $\wti\theta(z,x)$ satisfies
\eqref{4.18} since for $z\in\cO'$,
\begin{align}
W&(\wti\theta(z,\cdot),\wti\phi(z,\cdot))(x) =
W(\wti\theta(z,\cdot),\wti\phi(z,\cdot))(x_0) \\
&=
\f{\wti\phi'(z,x_0)}{\wti\phi(z,x_0)^2+\wti\phi'(z,x_0)^2}\wti\phi'(z,x_0)
+ \f{\wti\phi(z,x_0)}{\wti\phi(z,x_0)^2+\wti\phi'(z,x_0)^2}\wti\phi(z,x_0)
= 1.
\end{align}
\end{proof}

Having a system of two linearly independent solutions
$\wti\phi(z,x)$ and $\wti\theta(z,x)$ we introduce a function $\wti
m_+(z)$ such that the following solution of \eqref{4.3}
\begin{align} \lb{4.20}
\wti \psi_+(z,x) = \wti\theta(z,x) + \wti m_+(z)\wti\phi(z,x), \quad
x\in(a,\infty),
\end{align}
satisfies
\begin{align} \lb{4.21}
\wti\psi_{+}(z,\cdot)\in L^2([b,\infty);dx) \, \text{ for all
$b\in(a,\infty)$,}
\end{align}
for $z\in\cO'\backslash\bbR$. By Hypothesis \ref{h4.1}\,$(ii)$,
the solution $\wti\psi_+(z,\cdot)$ is proportional to
$\psi_+(z,\cdot,x_0)$. Hence, using \eqref{4.12} and \eqref{4.13}, one computes,
\begin{align}
m_+(z,x) &= \f{\wti\theta'(z,x)+\wti m_+(z)\wti\phi'(z,x)}
{\wti\theta(z,x)+\wti m_+(z)\wti\phi(z,x)}, \lb{4.28}
\\
\wti m_+(z) &= \f{\wti\theta(z,x)m_+(z,x)-\wti\theta'(z,x)}
{\wti\phi'(z,x)-\wti\phi(z,x)m_+(z,x)} =
\f{W(\wti\theta(z,\cdot),\psi_+(z,\cdot,x_0))}
{W(\psi_+(z,\cdot,x_0),\wti\phi(z,\cdot))} \lb{4.29}
\\
&= \f{\wti\theta(z,x)}{\wti\phi(z,x)}\f{m_+(z,x)}{m_-(z,x) - m_+(z,x)}-
\f{\wti\theta'(z,x)}{\wti\phi(z,x)}\f{1}{m_-(z,x) - m_+(z,x)}. \lb{4.30}
\end{align}
By \eqref{4.29}, $\wti m_+$ is independent of $x\in(a,\infty)$.

Having in mind the fact that $m_\pm(\cdot,x)$ are Herglotz and
anti-Herglotz functions, that $\wti\phi(z,x)\neq 0$ for
$z\in\CmR$, $|\Im(z)|$ sufficiently small, and that
$\wti\theta(z,x)$ and $\wti\theta'(z,x)$ are analytic with respect
to $z\in\cO'$, one concludes from
\eqref{4.30} that $\wti m_+$ is analytic in $\cO'\backslash\bbR$. In
contrast to $m_+$, the function
$\wti m_+$, in general, is not a Herglotz function. Nevertheless,
$\wti m_+$ shares  some properties with Herglotz functions which
are crucial for the proof of our main result, Theorem \ref{t4.6}.
Before we derive these properties we mention that by using
Hypothesis \ref{h4.1}\,$(iii)$, \eqref{4.18}, and \eqref{4.21},
a computation of the Green's function $G_+ (z,x,x')$ of $H_{+}$ yields
\begin{equation}
G_+ (z,x,x') = \begin{cases}
\wti\phi(z,x)\wti \psi_{+}(z,x'), & a < x \leq x', \\
\wti\phi(z,x')\wti \psi_{+}(z,x), & a < x' \leq x
\end{cases} \lb{4.34}
\end{equation}
and thus,
\begin{align}
\begin{split}
& ((H_{+}-zI)^{-1}f)(x) =\int^{\infty}_{a} dx' \, G_{+}(z,x,x')f(x'),
\lb{4.35} \\
& \hspace*{2.3cm} x\in (a,\infty), \quad f\in L^{2}([a,\infty);dx)
\end{split}
\end{align}
for $z\in \cO'\backslash\bbR$.

The basic properties of $\wti m_+$ then read as follows:

\begin{lemma} \lb{l4.5}
Assume Hypothesis \ref{h4.1}. Then the function $\wti m_+$
introduced in \eqref{4.20} satisfies the following properties:
\begin{enumerate}[$($i$)$]
\item $\wti m_+(z) = \ol{\wti m_+(\ol{z})}, \quad z\in\bbC_+, \;
|\Im(z)| \text{ sufficiently small.}$
\item
$\varepsilon|\wti m_{+}(\lambda+i\varepsilon)|\leq
C(\lambda_1,\lambda_2,\varepsilon_0) \, \text{ for } \, \lambda\in
[\lambda_1,\lambda_2], \; 0<\varepsilon\leq\varepsilon_0$.
\item
$\varepsilon |\Re(\wti m_{+}(\lambda+i\varepsilon))|
\underset{\varepsilon\to 0}{=}\oh(1) \, \text{ for } \,
\lambda\in [\lambda_1,\lambda_2], \; 0<\varepsilon\leq\varepsilon_0$.
\item
$-i \lim_{\varepsilon\downarrow 0} \varepsilon \wti
m_+(\lambda+i\varepsilon)=
\lim_{\varepsilon\downarrow 0} \varepsilon \Im(\wti
m_+(\lambda+i\varepsilon))$ exists for all $\lambda\in\bbR$ and is
nonnegative.
\item
$\wti m_+(\lambda+i0) =\lim_{\varepsilon\downarrow 0}\wti
m_+(\lambda+i\varepsilon) \, \text{ exists for a.e.\
$\lambda\in [\lambda_1,\lambda_2]$ and
} \\ \Im(\wti m_+(\lambda+i0))\geq 0 \, \text{ for a.e.\ $\lambda\in
[\lambda_1,\lambda_2]$}$.
\end{enumerate}
Here $0<\varepsilon_0=\varepsilon(\lambda_1,\lambda_2)$ is assumed to be
sufficiently small. Moreover, one can introduce a nonnegative measure
$d\wti\rho_+$ associated with $\wti m_+$ in a manner similar to the
Herglotz situation
\eqref{A.4} by
\begin{align} \lb{4.33}
\int_{(\lambda_1,\lambda_2]} d\wti\rho_+(\lambda) =
\wti\rho_+((\lambda_1,\lambda_2]) = \lim_{\delta\downarrow
0}\lim_{\varepsilon\downarrow 0}
\frac{1}{\pi}\int_{\lambda_1+\delta}^{\lambda_2+\delta} d\lambda
\, \Im(\wti m_+(\lambda+i\varepsilon)).
\end{align}
\end{lemma}
\begin{proof}
Since $\wti\phi(\lambda,x)$ and $\wti\theta(\lambda,x)$ are real-valued for
$(\lambda,x)\in\bbR\times (a,\infty)$, and analytic for $\lambda\in\cO'$ for fixed $x\in (a,\infty)$, an application of the Schwarz reflection principle yields
\begin{equation}
\wti\phi(z,x) = \ol{\wti\phi(\ol{z},x)}, \quad
\wti\theta(z,x) = \ol{\wti\theta(\ol{z},x)}, \quad x\in(a,\infty), \; z\in\cO'.  \lb{4.33a}
\end{equation}
Thus, picking real numbers $c$ and $d$ such that $a\leq c<d<\infty$,
\eqref{4.34} and \eqref{4.35} imply for the analog of \eqref{2.20} in the
present context of $H_+$,
\begin{align}
& \int_{\sigma(H_{+})} \f{d\big\|E_{H_+} (\lambda)\chi_{[c,d]}
\big\|_{L^2([a,\infty);dx)}^2}{\lambda-z} = \big(\chi_{[c,d]},(H_+
-zI)^{-1}
\chi_{[c,d]}\big)_{L^2([a,\infty);dx)} \no \\
& \quad = \int_c^d dx \int_c^x dx'\, \wti\theta (z,x) \wti\phi
(z,x') + \int_c^d dx\int_x^d dx'\,
\wti\phi (z,x)\wti\theta (z,x')  \lb{4.33b} \\
& \qquad + \wti m_+ (z) \bigg[\int_c^d dx \,
\wti\phi (z,x)\bigg]^2, \quad
z\in\bbC\backslash\sigma(H_+).  \no
\end{align}
Choosing $c(z_0),d(z_0)\in [a,\infty)$ such that
\begin{equation}
\int_{c(z_0)}^{d(z_0)} dx\, \wti\phi(z,x) \neq 0 \lb{4.33c}
\end{equation}
for $z$ in an open neighborhood $\cN(z_0)$ of $z_0\in
\bbC\backslash\sigma(H_+)$ with $\Im(z_0)$ sufficiently small (cf.\ the
proof of Lemma \ref{l2.3}), items $(i)$--$(v)$ follow from \eqref{4.33a} and \eqref{4.33b}
since the left-hand side in \eqref{4.33b},
\begin{equation}
\int_{\sigma(H_{+})} \f{d\big\|E_{H_+} (\lambda)\chi_{[c,d]}
\big\|_{L^2([a,\infty);dx)}^2} {\lambda-z}, \quad
z\in\bbC\backslash\sigma(H_+),
\end{equation}
is a Herglotz function and $\wti\phi(z,x)$, $\wti\theta(z,x)$ are
analytic with respect to $z\in\cO'$, where $\cO'\subseteq\cO$ is
an open neighborhood of $\bbR$. In addition, $\wti\phi(z,x)$ and
$\wti\theta(z,x)$ are real-valued for $(z,x)\in\bbR\times
(a,\infty)$. Next, we pick $\la_1,\la_2 \in \bbR$,
$\lambda_1<\lambda_2$, such that for some $c_0, d_0 \in
[a,\infty)$,
\begin{equation}
\int_{c_0}^{d_0} dx \, \wti\phi(z,x) \neq 0
\end{equation}
for all $z$ in a complex neighborhood of the interval
$(\lambda_1,\lambda_2)$. Then \eqref{4.33b} applied to
$z=\lambda+i\varepsilon$, for real-valued $\lambda$ in a neighborhood of
$(\lambda_1,\lambda_2)$, $0 < \varepsilon \leq \varepsilon_0$, implies
that $\wti\rho_+$ defined in \eqref{4.33} satisfies
\begin{align}
&\wti\rho_+((\lambda_1,\lambda_2]) =
\lim_{\delta\downarrow 0}\lim_{\varepsilon\downarrow 0} \f{1}{\pi}
\int_{\lambda_1+\delta}^{\lambda_2+\delta} d\lambda \,
\Im(\wti m_+(\lambda+i\varepsilon)) \no \\
& \quad = \lim_{\delta\downarrow 0}\lim_{\varepsilon\downarrow 0}
\f{1}{\pi} \int_{\lambda_1+\delta}^{\lambda_2+\delta} d\lambda \,
\Im \bigg\{\int_{\sigma(H_+)}\f{d\big\|E_{H_+}(\lambda')\chi_{[c_0,
d_0]}\big\|_{L^2([a,\infty);dx)}^2}{\lambda'-\lambda-i\varepsilon} \no \\
& \qquad \quad \times \bigg[\bigg(\int_{c_0}^{d_0} dx \,
\wti\phi(\lambda,x)\bigg)^2+2i\varepsilon  \bigg(\int_{c_0}^{d_0} dx
\, (d/dz)\wti\phi(z,x)|_{z=\lambda}\bigg)+\Oh(\varepsilon^2)\bigg]^{-1}
\no \\
& \qquad + \Oh(\varepsilon) \bigg\} \no \\
& \quad = \int_{(\lambda_1,\lambda_2]} d\big\|E_{H_+}
(\la)\chi_{[c_0,d_0]} \big\|_{L^2([a,\infty);dx)}^2
\left[\int_{c_0}^{d_0} dx\, \wti\phi(\la,x)\right]^{-2},
\end{align}
using item $(ii)$, item $(iii)$, the dominated convergence theorem, and
the analog of \eqref{2.43} applied to the present context. Hence,
$\wti\rho_+$ generates the nonnegative measure $d\wti\rho_+$.
\end{proof}

Next, we relate the family of spectral projections,
$\{E_{H_{+}}(\lambda)\}_{\lambda\in\bbR}$, of the self-adjoint
operator $H_{+}$ and the spectral function
$\wti\rho_{+}(\lambda)$, $\lambda\in\bbR$, defined in
\eqref{4.33}.

We first note that for $F\in C(\bbR)$,
\begin{align}
&\big(f,F(H_+)g\big)_{L^2([a,\infty);dx)}= \int_{\bbR}
d\big(f,E_{H_+}(\lambda)g\big)_{L^2([a,\infty);dx)}\,
F(\lambda), \\
& f, g \in\dom(F(H_+)) \no \\
& \quad =\bigg\{h\in L^2([a,\infty);dx)
\,\bigg|\,
\int_{\bbR} d\big\|E_{H_+}(\lambda)h\big\|_{L^2([a,\infty);dx)}^2
\, |F(\lambda)|^2 < \infty\bigg\}. \no
\end{align}

\begin{theorem} \lb{t4.6}
Let $f,g \in C^\infty_0((a,\infty))$, $F\in C(\bbR)$, and
$\lambda_1, \lambda_2 \in\bbR$, $\lambda_1<\lambda_2$. Then,
\begin{equation}
\big(f,F(H_{+})E_{H_{+}}((\lambda_1,\lambda_2])
g\big)_{L^2([a,\infty);dx)} =  \big(\hatt
f_{+},M_FM_{\chi_{(\lambda_1,\lambda_2]}} \hatt
g_{+}\big)_{L^2(\bbR;d\wti\rho_{+})}, \lb{4.37}
\end{equation}
where we introduced the notation
\begin{equation}
\hatt h_{+}(\lambda)=\int_a^\infty dx \, \wti\phi(\lambda,x)h(x),
\quad \lambda \in\bbR, \; h\in C^\infty_0((a,\infty)), \lb{4.38}
\end{equation}
and $M_G$ denotes again the maximally defined operator of multiplication
by the $d\wti\rho_{+}$-measurable function $G$ in the Hilbert
space $L^2(\bbR;d\wti\rho_{+})$,
\begin{align}
\begin{split}
& (M_G\hatt h)(\lambda)=G(\lambda)\hatt h(\lambda)
\, \text{ for a.e.\ $\lambda\in\bbR$},  \\
& \hatt h\in\dom(M_G)=\big\{\hatt k \in L^2(\bbR;d\wti\rho_{+}) \,|\,
G\hatt k \in L^2(\bbR;d\wti\rho_{+})\big\}.
\end{split} \lb{4.38a}
\end{align}
\end{theorem}
\begin{proof}
The point of departure for deriving \eqref{4.37} is again Stone's
formula \eqref{2.26a} applied to $T=H_{+}$,
\begin{align}
&\big(f,F(H_{+})E_{H_{+}}((\lambda_1,\lambda_2])
g\big)_{L^2([a,\infty);dx)} \no \\ & \quad =
\lim_{\delta\downarrow 0}\lim_{\varepsilon\downarrow 0}
\frac{1}{2\pi i} \int_{\lambda_1+\delta}^{\lambda_2+\delta}
d\lambda \, F(\lambda) \big[\big(f,(H_{+}-(\lambda+i\varepsilon)
I)^{-1}g\big)_{L^2([a,\infty);dx)}
\no \\
& \hspace*{4.9cm} - \big(f,(H_{+}-(\lambda-i\varepsilon)
I)^{-1}g\big)_{L^2([a,\infty);dx)}\big]. \lb{4.39}
\end{align}
Insertion of \eqref{4.34} and \eqref{4.35} into \eqref{4.39} then
yields the following:
\begin{align}
&\big(f,F(H_{+})E_{H_{+}}((\lambda_1,\lambda_2])
g\big)_{L^2([a,\infty);dx)} = \lim_{\delta\downarrow
0}\lim_{\varepsilon\downarrow 0} \frac{1}{2\pi i}
\int_{\lambda_1+\delta}^{\lambda_2+\delta} d\lambda \, F(\lambda) \no \\
& \quad \times\int_a^\infty dx
\bigg\{\bigg[\ol{f(x)} \wti\psi_{+}(\lambda
+i\varepsilon,x) \int_a^x dx'\,
\wti\phi(\lambda+i\varepsilon,x')g(x') \no \\
& \hspace*{2.4cm} + \ol{f(x)} \wti\phi(\lambda+i\varepsilon,x)
\int_x^\infty dx'\, \wti\psi_{+}(\lambda+i\varepsilon,x') g(x')
\bigg]  \no \\
& \hspace*{1.9cm} -\bigg[\ol{f(x)}
\wti\psi_{+}(\lambda -i\varepsilon,x) \int_a^x dx'\,
\wti\phi(\lambda-i\varepsilon,x')g(x') \no \\
& \hspace*{2.4cm}  + \ol{f(x)} \wti\phi(\lambda-i\varepsilon,x)
\int_x^\infty dx' \,\wti\psi_{+}(\lambda-i\varepsilon,x')
g(x')\bigg]\bigg\}.  \lb{4.40}
\end{align}
Freely interchanging the $dx$ and $dx'$
integrals with the limits and the $d\lambda$ integral
(since all integration domains are finite and all
integrands are continuous), and inserting expression \eqref{4.20}
for $\wti\psi_{+}(z,x)$ into \eqref{4.40}, one obtains
\begin{align}
&\big(f,F(H_{+})E_{H_{+}}((\lambda_1,\lambda_2])
g\big)_{L^2([a,\infty);dx)}
=\int_a^\infty dx\, \ol{f(x)}\bigg\{\int_a^x dx' \, g(x') \no \\
& \qquad \times \lim_{\delta\downarrow
0}\lim_{\varepsilon\downarrow 0} \frac{1}{2\pi i}
\int_{\lambda_1+\delta}^{\lambda_2+\delta} d\lambda \, F(\lambda)
\Big[\big[\wti\theta(\lambda,x) + \wti
m_{+}(\lambda+i\varepsilon)\wti\phi(\lambda,x)\big]
\wti\phi(\lambda,x')   \no \\
& \hspace*{4cm} -\big[\wti\theta(\lambda,x) + \wti
m_{+}(\lambda-i\varepsilon)\wti\phi(\lambda,x)\big]
\wti\phi(\lambda,x')\Big] \no \\
& \quad +\int_x^\infty dx'\, g(x') \lim_{\delta\downarrow 0}
\lim_{\varepsilon\downarrow 0} \frac{1}{2\pi i}
\int_{\lambda_1+\delta}^{\lambda_2+\delta} d\lambda \, F(\lambda)
\lb{4.41} \\
& \hspace*{3.5cm} \times \Big[ \wti\phi(\lambda,x)
\big[\wti\theta(\lambda,x') + \wti
m_{+}(\lambda+i\varepsilon)\wti\phi(\lambda,x')\big]   \no \\
& \hspace{4cm}  -\wti\phi(\lambda,x) \big[\wti\theta(\lambda,x') +
\wti
m_{+}(\lambda-i\varepsilon)\wti\phi(\lambda,x')\big]\Big]\bigg\}.
\no
\end{align}
Here we employed the fact that for fixed $x\in (a,\infty)$,
$\wti\phi(z,x)$, $\wti\theta(z,x)$ are analytic with respect to
$z\in\cO$ and real-valued for
$z\in\bbR$, the fact that $\wti\phi(z,\cdot), \wti\theta(z,\cdot) \in
AC_{\loc}((a,\infty))$, and hence that
\begin{align}
\begin{split}
\wti\phi(\lambda\pm i\varepsilon,x)
&\underset{\varepsilon\downarrow 0}{=} \wti\phi(\lambda,x) \pm
i\varepsilon (d/dz)\wti\phi(z,x)|_{z=\lambda} + \Oh(\ve^2),  \lb{4.42}
\\
\wti\theta(\lambda\pm i\varepsilon,x)
&\underset{\varepsilon\downarrow 0}{=} \wti\theta(\lambda,x) \pm
i\varepsilon (d/dz) \wti\theta(z,x)|_{z=\lambda} + \Oh(\ve^2),
\end{split}
\end{align}
with $\Oh(\varepsilon^2)$ being uniform with respect to
$(\lambda,x)$ as long as $\lambda$ and $x$ vary in compact
subsets of $\bbR\times (a,\infty)$. (Here real-valuedness of $\wti\phi(z,x)$ and
 $\wti\theta(z,x)$ for $z\in\bbR$, $x\in (a,\infty)$ yields a purely imaginary $\Oh(\varepsilon)$-term in \eqref{4.42}.) Moreover, we used items $(ii)$ and
$(iii)$ of Lemma \ref{l4.5} to replace $\wti\phi(\lambda\pm
i\varepsilon,x)$ and $\wti\theta(\lambda\pm i\varepsilon,x)$ by
$\wti\phi(\lambda,x)$ and $\wti\theta(\lambda,x)$ under the $d\lambda$
integrals in \eqref{4.41}. Cancelling appropriate terms in \eqref{4.41},
simplifying the remaining terms, and using item $(i)$ of Lemma \ref{l4.5}
then yield
\begin{align}
&\big(f,F(H_{+})E_{H_{+}}((\lambda_1,\lambda_2])
g\big)_{L^2([a,\infty);dx)}
=\int_a^\infty dx\, \ol{f(x)}\int_a^\infty dx' \, g(x') \no \\
& \quad \times \lim_{\delta\downarrow
0}\lim_{\varepsilon\downarrow 0} \frac{1}{\pi}
\int_{\lambda_1+\delta}^{\lambda_2+\delta} d\lambda \, F(\lambda)
\wti\phi(\lambda,x)\wti\phi(\lambda,x') \Im(\wti
m_{+}(\lambda+i\varepsilon)). \lb{4.43}
\end{align}
Using \eqref{4.33},
\begin{align}
\int_{\bbR} d\wti\rho_{+}(\lambda)\, h(\lambda) &=
\lim_{\varepsilon\downarrow 0} \frac{1}{\pi}\int_{\bbR} d\lambda
\, \Im(\wti m_{+}(\lambda+i\varepsilon))\, h(\lambda), \quad
h\in C_0(\bbR), \lb{4.44} \\
\int_{(\lambda_1,\lambda_2]} d\wti\rho_{+}(\lambda)\, k(\lambda)
&= \lim_{\delta\downarrow 0} \lim_{\varepsilon\downarrow 0}
\frac{1}{\pi} \int_{\lambda_1+\delta}^{\lambda_2+\delta} d\lambda
\, \Im(\wti m_{+}(\lambda+i\varepsilon))\, k(\lambda), \quad k\in
C(\bbR), \lb{4.45}
\end{align}
and hence
\begin{align}
&\big(f,F(H_{+})E_{H_{+}}((\lambda_1,\lambda_2])
g\big)_{L^2([a,\infty);dx)} \no \\ & \quad =\int_a^\infty dx\,
\ol{f(x)}\int_a^\infty dx' \, g(x') \int_{(\lambda_1,\lambda_2]}
d\wti\rho_{+}(\lambda) F(\lambda)
\wti\phi(\lambda,x)\wti\phi(\lambda,x') \no \\
& \quad = \int_{(\lambda_1,\lambda_2]} d\wti\rho_{+}(\lambda)
F(\lambda) \, \ol{\hatt f_{+}(\lambda)}\, \hatt g_{+}(\lambda),
\lb{4.46}
\end{align}
using \eqref{4.38} and interchanging the $dx$, $dx'$ and
$d\wti\rho_{+}$ integrals once more.
\end{proof}

Again one can improve on Theorem \ref{t4.6} and remove
the compact support restrictions on $f$ and $g$ in the usual way.
To this end we consider the map
\begin{equation}
\widetilde U_{+}\colon \begin{cases} C_0^\infty((a,\infty))\to
L^2(\bbR; d\wti\rho_{+}) \\[1mm]
\hspace*{1.6cm} h \mapsto \hatt h_{+}(\cdot)= \int_a^\infty dx\,
\wti\phi(\cdot,x) h(x). \end{cases} \lb{4.46a}
\end{equation}
Taking $f=g$, $F=1$, $\lambda_1\downarrow -\infty$, and
$\lambda_2\uparrow \infty$ in \eqref{4.37} then shows that $\widetilde
U_{+}$ is a densely defined isometry in $L^2([a,\infty);dx)$, which
extends by continuity to an isometry on $L^2([a,\infty);dx)$. The latter
is denoted by $U_{+}$ and given by
\begin{equation}
U_{+}\colon \begin{cases}L^2([a,\infty);dx)\to L^2(\bbR;
d\wti\rho_{+}) \\[1mm]
\hspace*{1.95cm}  h \mapsto \hatt h_{+}(\cdot)=
\slimes_{b\uparrow\infty}\int_a^b dx\, \wti\phi(\cdot,x) h(x),
\end{cases}  \lb{4.46b}
\end{equation}
where $\slimes$ refers to the $L^2(\bbR; d\wti\rho_+)$-limit.

The calculation in \eqref{4.46} also yields
\begin{equation}
(E_{H_{+}}((\lambda_1,\lambda_2])g)(\cdot)
=\int_{(\lambda_1,\lambda_2]} d\wti\rho_{+}(\lambda)\,
\wti\phi(\lambda,\cdot) \hatt g_{+}(\lambda), \quad g\in
C_0^\infty((a,\infty)) \lb{4.46c}
\end{equation}
and subsequently, \eqref{4.46c} extends to all $g\in
L^2([a,\infty);dx)$ by continuity. Moreover, taking
$\lambda_1\downarrow -\infty$ and $\lambda_2\uparrow \infty$ in
\eqref{4.46c} and using
\begin{equation}
\slim_{\lambda\downarrow -\infty} E_{H_{+}}(\lambda)=0, \quad
\slim_{\lambda\uparrow \infty}
E_{H_{+}}(\lambda)=I_{L^2([a,\infty);dx)},
\end{equation}
where
\begin{equation}
E_{H_{+}}(\lambda)=E_{H_{+}}((-\infty,\lambda]),  \quad \lambda\in\bbR,
\end{equation}
then yields
\begin{equation}
g(\cdot)=\slimes_{\mu_1\downarrow -\infty, \mu_2\uparrow\infty}
\int_{\mu_1}^{\mu_2} d\wti\rho_{+}(\lambda)\,
\wti\phi(\lambda,\cdot) \hatt g_{+}(\lambda), \quad g\in
L^2([a,\infty);dx), \lb{4.46d}
\end{equation}
where $\slimes$ refers to the $L^2([a,\infty); dx)$-limit.

In addition, one can show that the map $U_{+}$ in \eqref{4.46b} is
onto and hence that $U_{+}$ is unitary (i.e., $U_{+}$ and
$U_{+}^{-1}$ are isometric isomorphisms between
$L^2([a,\infty);dx)$ and $L^2(\bbR;d\wti\rho_{+})$) with
\begin{equation}
U_{+}^{-1} \colon \begin{cases} L^2(\bbR;d\wti\rho_{+}) \to
L^2([a,\infty);dx)  \\[1mm]
\hspace*{1.5cm} \hatt h \mapsto \slimes_{\mu_1\downarrow -\infty,
\mu_2\uparrow\infty} \int_{\mu_1}^{\mu_2} d\wti\rho_{+}(\lambda)\,
\wti\phi(\lambda,\cdot) \hatt h(\lambda). \end{cases} \lb{4.46e}
\end{equation}
To show this one can follow the corresponding proof of unitarity of $U_{+,\alpha}$ in \eqref{2.45a}--\eqref{2.45b} line by line.

We sum up these considerations in a variant of the spectral theorem for
(functions of) $H_{+}$.

\begin{theorem}
Let $F\in C(\bbR)$, Then,
\begin{equation}
U_{+} F(H_{+})U_{+}^{-1} = M_F \lb{4.46l}
\end{equation}
in $L^2(\bbR;d\wti\rho_{+})$ $($cf.\ \eqref{4.38a}$)$. Moreover,
\begin{align}
& \sigma(F(H_{+}))= \essran_{d\wti\rho_{+}}(F), \lb{4.46m} \\
& \sigma(H_{+})=\supp\,(d\wti\rho_{+}),  \lb{4.46n}
\end{align}
and the spectrum of $H_{+}$ is simple.
\end{theorem}

Simplicity of the spectrum of $H_+$ is consistent with the observation that
\begin{align}\begin{split}
&\det\left(\Im(M(\lambda +i0,x_0))\right) \\
&\quad =\det\left(\begin{pmatrix} \f{\Im(m_+(\la+i0,x_0))}
{\abs{m_-(\la+i0,x_0)-m_+(\la+i0,x_0)}^2} &
\f{m_-(\la+i0,x_0)\Im(m_+(\la+i0,x_0))}
{\abs{m_-(\la+i0,x_0)-m_+(\la+i0,x_0)}^2}
\\
\f{m_-(\la+i0,x_0)\Im(m_+(\la+i0,x_0))}
{\abs{m_-(\la+i0,x_0)-m_+(\la+i0,x_0)}^2} &
\f{\abs{m_-(\la+i0,x_0)}^2\Im(m_+(\la+i0,x_0))}
{\abs{m_-(\la+i0,x_0)-m_+(\la+i0,x_0)}^2}
\end{pmatrix}\right) \\
&\quad = 0 \, \text{ for a.e.\ $\lambda\in\bbR$}
\end{split}\end{align}
since by Lemma \ref{l4.2}, $m_-(z,x_0)$ is meromorphic and
real-valued for $z\in\bbR$. In this context we also refer to
\cite{Gi89}, \cite{Gi98}, \cite{Ka62}, \cite{Ka63}, where necessary and
sufficient conditions for simplicity of the spectrum in terms of properties
of $m_\pm(\cdot,x_0)$ can be found.

Next, we consider the alternative way of deriving the
(matrix-valued) spectral function corresponding to a reference point
$x_0\in(a,\infty)$ and subsequently compare the two approaches.

As in the half-line context in Section \ref{s2} we introduce the usual
fundamental system of solutions $\phi (z,\cdot,x_0)$ and
$\theta (z,\cdot,x_0)$, $z\in\bbC$, of
\begin{equation}
(\tau_+ \psi)(z,x) = z \psi(z,x), \quad x\in (a,\infty) \lb{4.46o}
\end{equation}
with respect to a fixed reference point $x_0\in(a,\infty)$, satisfying the
initial conditions at the point $x=x_0$,
\begin{equation}
\phi(z,x_0,x_0)=\theta'(z,x_0,x_0)=0, \quad
\phi'(z,x_0,x_0)=\theta(z,x_0,x_0)=1. \lb{4.46p}
\end{equation}
Again we note that for any fixed $x, x_0\in (a,\infty)$,
$\phi(z,x,x_0)$ and $\theta(z,x,x_0)$ are entire with respect to $z$ and
that
\begin{equation}
W(\theta(z,\cdot,x_0),\phi(z,\cdot,x_0))(x)=1, \quad
z\in\bbC.  \lb{4.46q}
\end{equation}
The {\it Weyl--Titchmarsh solutions} $\psi_{\pm,\alpha}(z,\cdot,x_0)$,
$z\in\bbC\backslash\bbR$, of \eqref{4.46o} are uniquely characterized by
\begin{align}
\begin{split}
&\psi_{-}(z,\cdot,x_0)\in L^2([a,x_0];dx), \;
\psi_{+}(z,\cdot,x_0)\in L^2([x_0,\infty);dx), \;
z\in\bbC\backslash\bbR, \\
&\psi_{\pm}(z,x_0,x_0)=1. \lb{4.46r}
\end{split}
\end{align}
The normalization in \eqref{4.46r} shows that
$\psi_{\pm}(z,\cdot,x_0)$ are of the type
\begin{equation}
\psi_{\pm}(z,x,x_0)=\theta(z,x,x_0)
+m_{\pm}(z,x_0)\phi(z,x,x_0),
\quad  z\in\bbC\backslash\bbR, \; x\in\bbR \lb{4.46s}
\end{equation}
for some coefficients $m_{\pm}(z,x_0)$, the half-line
{\it Weyl--Titchmarsh $m$-functions} associated with $\tau_+$ and $x_0$.
Again we recall the fundamental identity
\begin{align}
&\int^{x_0}_{a} dx\,\psi_{-}(z_{1},x,x_0)
\psi_{-}(z_{2},x,x_0) = -\frac{m_{-}(z_{1},x_0)-
m_{-}(z_{2},x_0)}{z_{1}-z_{2}}, \lb{4.46t} \\
&\int_{x_0}^{\infty} dx\,\psi_{+}(z_{1},x,x_0)
\psi_{+}(z_{2},x,x_0) = \frac{m_{+}(z_{1},x_0)-
m_{+}(z_{2},x_0)}{z_{1}-z_{2}}, \lb{4.46ta} \\
&\hspace*{6.1cm}  z_1, z_2 \in\bbC\backslash\bbR, \; z_1 \neq z_2, \no
\end{align}
and as before one concludes
\begin{equation}
\ol{m_{\pm}(z,x_0)} = m_{\pm}(\ol z,x_0), \quad
z\in\bbC\backslash\bbR.  \lb{4.46tb}
\end{equation}
Choosing $z_1=z$, $z_2=\ol z$ in \eqref{4.46t}, \eqref{4.46ta} one infers
\begin{align}
\begin{split}
\int^{x_0}_{a} dx\,|\psi_{-}(z,x,x_0)|^2
                               & = -\frac{\Im(m_{-}(z,x_0))}{\Im(z)},
\lb{4.46tc} \\
\int_{x_0}^{\infty} dx\,|\psi_{+}(z,x,x_0)|^2
                               & = \frac{\Im(m_{+}(z,x_0))}{\Im(z)}, \quad z
\in\bbC\backslash\bbR.
\end{split}
\end{align}
Since $m_{\pm}(\cdot,x_0)$ are analytic on $\bbC\backslash\bbR$,
$\pm m_{\pm}(\cdot,x_0)$ are Herglotz functions.

The Green's function $G_+(z,x,x')$ of $H_+$ then admits the
alternative representation (cf.\ also \eqref{4.34}, \eqref{4.35})
\begin{align}
& G_+(z,x,x') =
\f{1}{W(\psi_+(z,\cdot,x_0),\psi_-(z,\cdot,x_0))}\begin{cases}
\psi_{-}(z,x,x_0)\psi_{+}(z,x',x_0), & x \leq x', \\
\psi_{-}(z,x',x_0)\psi_{+}(z,x,x_0), & x' \leq x,
\end{cases} \no \\
&  \hspace*{9.7cm} z\in\bbC\backslash\bbR \lb{4.46u}
\end{align}
with
\begin{equation}
W(\psi_+(z,\cdot,x_0),\psi_-(z,\cdot,x_0))=m_{-}(z,x_0)-
m_{+}(z,x_0), \quad z\in\bbC\backslash\bbR.  \lb{4.46v}
\end{equation}
Thus,
\begin{align}
\begin{split}
&((H_+-zI)^{-1}f)(x)
=\int_{a}^\infty dx' \, G_+(z,x,x')f(x'), \\
& \hspace*{1.1cm} z\in\bbC\backslash\bbR, \;
x\in[a,\infty), \; f\in L^{2}([a,\infty);dx). \lb{4.46w}
\end{split}
\end{align}

Given $m_\pm(z,x_0)$, we also introduce the $2\times 2$ matrix-valued
Weyl--Titchmarsh function
\begin{align}
M(z,x_0)&=\begin{pmatrix}
\f{1}{m_{-}(z,x_0)-m_{+}(z,x_0)} &
\f{1}{2}\f{m_{-}(z,x_0)
+m_{+}(z,x_0)}{m_{-}(z,x_0)-m_{+}(z,x_0)} \\
\f{1}{2}\f{m_{-}(z,x_0)
+m_{+}(z,x_0)}{m_{-}(z,x_0)-m_{+}(z,x_0)} &
\f{m_{-}(z,x_0)m_{+}(z,x_0)}{m_{-}(z,x_0)
-m_{+}(z,x_0)}\end{pmatrix}, \quad z\in\bbC\backslash\bbR. \lb{4.46x}
\end{align}
$M(z,x_0)$ is a Herglotz matrix with representation
\begin{align}
\begin{split}
& M(z,x_0)=C(x_0)+\int_{{\mathbb{R}}}
d\Omega (\lambda,x_0)\bigg[\frac{1}{\lambda -z}-\frac{\lambda}
{1+\lambda^2}\bigg], \quad z\in\bbC\backslash\bbR, \lb{4.46y} \\
& C(x_0)=C(x_0)^*, \quad \int_{\bbR}
\f{\|d\Omega(\lambda,x_0)\|}{1+\lambda^2} <\infty,
\end{split}
\end{align}
where
\begin{align}
\Omega((\lambda_1,\lambda_2],x_0) &=
\f{1}{\pi}\lim_{\delta\downarrow 0}\lim_{\varepsilon\downarrow
0}\int^{\lambda_2+\delta}_{\lambda_1+\delta}d\lambda \,
\Im(M(\lambda +i\varepsilon,x_0)), \quad \lambda_1, \lambda_2
\in\bbR, \; \lambda_1<\lambda_2. \lb{4.46z}
\end{align}

Again one can of course replace $z\in\bbC\backslash\bbR$ by
$z\in\bbC\backslash\sigma(H_+)$ in formulas \eqref{4.46r}--\eqref{4.46y}.

Next, we relate once more the family of spectral projections,
$\{E_{H_+}(\lambda)\}_{\lambda\in\bbR}$, of the self-adjoint
operator $H_+$ and the $2\times 2$ matrix-valued nondecreasing spectral
function $\Omega(\lambda,x_0)$, $\lambda\in\bbR$, which
generates the matrix-valued measure in the Herglotz representation
\eqref{4.46y} of $M(z,x_0)$.

\begin{theorem} \lb{t4.7}
Let $f,g \in C^\infty_0((a,\infty))$, $F\in C(\bbR)$,
$x_0\in(a,\infty)$, and $\lambda_1, \lambda_2 \in\bbR$,
$\lambda_1<\lambda_2$. Then,
\begin{align}
&\big(f,F(H_+)E_{H_+}((\lambda_1,\lambda_2])g\big)_{L^2([a,\infty);dx)}
= \big(\hatt f(\cdot,x_0),M_FM_{\chi_{(\lambda_1,\lambda_2]}}
\hatt g(\cdot,x_0)\big)_{L^2(\bbR;d\Omega(\cdot,x_0))}
\no  \\
& \quad = \int_{(\lambda_1,\lambda_2]} \ol{\hatt
f(\lambda,x_0)^\top} \, d\Omega(\lambda,x_0) \, \hatt
g(\lambda,x_0)F(\lambda), \lb{4.47}
\end{align}
where we introduced the notation
\begin{align}
\begin{split}
&\hatt h_{0}(\lambda,x_0)=\int_a^\infty dx \,
\theta(\lambda,x,x_0) h(x), \quad \hatt
h_{1}(\lambda,x_0)=\int_a^\infty dx \,
\phi(\lambda,x,x_0) h(x),  \\
&\hatt h(\lambda,x_0)=(\hatt h_{0}(\lambda,x_0), \hatt
h_{1}(\lambda,x_0))^\top,  \quad \lambda \in\bbR, \; h\in
C^\infty_0((a,\infty)),
\end{split} \lb{4.48}
\end{align}
and $M_G$ denotes the maximally defined operator of multiplication
by the $d\Omega^{\rm tr}(\cdot,x_0)$-measurable function $G$ in the Hilbert
space
$L^2(\bbR;d\Omega(\cdot,x_0))$,
\begin{align}
\begin{split}
& (M_G\hatt h)(\lambda)=G(\lambda)\hatt h(\lambda)
=\big(G(\lambda) \hatt h_0(\lambda), G(\lambda) \hatt h_1(\lambda)\big)^\top
\, \text{ for a.e.\ $\lambda\in\bbR$},  \\
& \hatt h\in\dom(M_G)=\big\{\hatt k \in L^2(\bbR;d\Omega(\cdot,x_0))
\,|\, G \hatt k \in L^2(\bbR;d\Omega(\cdot,x_0))\big\}.
\end{split} \lb{4.49}
\end{align}
\end{theorem}

We omit the proof of Theorem \ref{t4.7} since it parallels that of
Theorem \ref{t2.9}.

Repeating the proof of Theorem \ref{t2.10} one also obtains the
following result.

\begin{theorem}  \lb{t4.9}
Let $F\in C(\bbR)$, $x_0\in(a,\infty)$,
\begin{align}
& U(x_0)\colon \begin{cases} L^2 ([a,\infty);dx) \to
L^2(\bbR; d\Omega(\cdot,x_0)) \\[1mm]
                  \hspace*{1.98cm} h \mapsto \hatt h(\cdot,x_0) =\big(\hatt
h_{0}(\cdot,x_0), \hatt h_{1}(\cdot,x_0)\big)^\top, \end{cases}  \lb{4.49a} \\
& \hatt h(\cdot,x_0)=\begin{pmatrix} \hatt h_{0}(\cdot,x_0) \\
\hatt h_{1}(\cdot,x_0) \end{pmatrix} =
\slimes_{b\downarrow a, c \uparrow\infty} \begin{pmatrix}
\int_{b}^c dx \, \theta(\cdot,x,x_0) h(x) \\
\int_{b}^c dx \, \phi(\cdot,x,x_0) h(x) \end{pmatrix}, \no
\end{align}
where $\slimes$ refers to the $L^2(\bbR; d\Omega(\cdot,x_0))$-limit and
\begin{align}
&U(x_0)^{-1} \colon \begin{cases}
L^2(\bbR;d\Omega(\cdot,x_0)) \to  L^2([a,\infty);dx)  \\[1mm]
\hspace*{2.2cm} \hatt h \mapsto  h, \end{cases} \\
& h(\cdot)= \slimes_{\mu_1\downarrow -\infty, \mu_2\uparrow
\infty} \int_{\mu_1}^{\mu_2} (\theta(\lambda,\cdot,x_0),
\phi(\lambda,\cdot,x_0)) \, d\Omega(\lambda,x_0)\, \hatt
h(\lambda,x_0), \no
\end{align}
where $\slimes$ refers to the $L^2([a,\infty); dx)$-limit. Then,
\begin{equation}
U(x_0) F(H_+)U(x_0)^{-1} = M_F
\end{equation}
in $L^2(\bbR;d\Omega(\cdot,x_0))$ $($cf.\ \eqref{4.49}$)$.
Moreover,
\begin{equation}
\sigma(H_+)=\supp\,(d\Omega(\cdot,x_0))=\supp\,
(d\Omega^{\rm tr}(\cdot,x_0)).
\end{equation}
\end{theorem}

\begin{corollary} \lb{c4.10}
The expansions in \eqref{4.38} and \eqref{4.49a} are related by,
\begin{align}
\hatt h_+(\la) = \wti\phi(\la,x_0) \hatt h_0(\la,x_0) +
\wti\phi'(\la,x_0) \hatt h_1(\la,x_0), \quad \lambda\in\sigma(H_+).
\lb{4.50}
\end{align}
The measures $d\wti\rho_+$ and $d\Omega(\cdot,x_0)$ are related
by,
\begin{align}
d\wti\rho_+(\la) &=
\f{\wti\theta(\la,x_0)}{\wti\phi(\la,x_0)}d\Omega_{0,1}(\la,x_0) -
\f{\wti\theta'(\la,x_0)}{\wti\phi(\la,x_0)}d\Omega_{0,0}(\la,x_0)
\lb{4.51} \\
&=
\f{\wti\phi'(\la,x_0)}{\wti\phi(\la,x_0)}
\f{1}{\wti\phi(\la,x_0)^2+\wti\phi'(\la,x_0)^2}d\Omega_{0,1}(\la,x_0)
\no \\
& \quad +
\f{1}{\wti\phi(\la,x_0)^2+\wti\phi'(\la,x_0)^2}d\Omega_{0,0}(\la,x_0),
\quad \lambda\in\sigma(H_+). \lb{4.52}
\end{align}
\end{corollary}
\begin{proof}
\eqref{4.50} follows from \eqref{4.9}, \eqref{4.38}, and
\eqref{4.48}. \eqref{4.51} and \eqref{4.52} follow from \eqref{4.7}, \eqref{4.19},
\eqref{4.30}, \eqref{4.33}, \eqref{4.46x}, and \eqref{4.46z}.
\end{proof}

Finally, we illustrate the applicability of our approach to strongly singular
potentials by verifying Hypothesis \ref{h4.1} under very general
circumstances.

We start with a simple example first.

\begin{example} \lb{e4.11}
The class of potentials $V$ of the form
\begin{equation}
V(x) = \f{\ga^2-1/4}{x^2} + \wti V(x), \quad \ga \in[1,\infty), \;
x\in(0,\infty),  \lb{4.53}
\end{equation}
where $\wti V$ is a real-valued measurable function on $[0,\infty)$ such
that 
\begin{equation}
\wti V\in L^1([0,b];x\,dx) \, \text{ for all $b>0$,}
\lb{4.54}
\end{equation}
assuming that $\tau_+=-d^2/dx^2+[\gamma^2-(1/4)]x^{-2}+\wti V(x)$
is in the limit point case at $\infty$, satisfies Hypothesis
\ref{h4.1}.
\end{example}

To verify that the potential $V$ in \eqref{4.53} indeed satisfies
Hypothesis \ref{h4.1} we first state the following result. As kindly
pointed out to us by Don Hinton, this is a special case of his Theorem\ 1
in \cite{Hi74}. For convenience of the reader we include the following
elementary and short proof we found independently (and which differs from
the proof in \cite{Hi74}).

\begin{lemma} \rm{(\cite{Hi74}.)} \lb{l4.12}
Let $b\in (0,\infty)$. Then the differential expression $\tau_+$ given by
\begin{equation}
\tau_+=-\f{d^2}{dx^2}+\f{\gamma^2-(1/4)}{x^{2}}+\wti V(x), \quad x\in
(0,b),  \;  \gamma\in [1,\infty),  \lb{4.55}
\end{equation}
with $\wti V$ a real-valued and measurable function on $[0,b]$ satisfying
\begin{equation}
\wti V\in L^1([0,b];x\,dx),  \lb{4.56}
\end{equation}
is in the limit point case at $x=0$.
\end{lemma}
\begin{proof}
Consider a solution $\theta$ of
\begin{align} \lb{4.89}
\begin{split}
& (\tau_+ \theta)(x) = 0, \quad x\in(0,b),
\\
&\theta(x_0) = x_0^{1/2-\ga}, \; \theta'(x_0) =
(1/2-\ga)x_0^{-1/2-\ga} \quad \text{for some } x_0\in(0,b).
\end{split}
\end{align}
By the ``variation of constants'' formula, $\theta$ satisfies
\begin{align} \lb{4.90}
\theta(x) = x^{1/2-\ga} + \f{1}{2\ga}\int_{x_0}^x dt
\big[x^{1/2+\ga}t^{1/2-\ga}-x^{1/2-\ga}t^{1/2+\ga}\big]
\wti V(t)\theta(t).
\end{align}
Introducing
\begin{align} \lb{4.91}
&\theta_0 (x) = x^{1/2-\ga}, \no \\
&\theta_k (x) =
\f{1}{2\ga}\int_{x_0}^x dt
\big[x^{1/2+\ga}t^{1/2-\ga}-x^{1/2-\ga}t^{1/2+\ga}\big]
\wti V(t)\theta_{k-1}(t), \quad k\in\bbN,
\end{align}
and estimating $\theta_k$ by
\begin{align} \lb{4.92}
\abs{\theta_k (x)} \leq x^{1/2-\ga}
\f{1}{k!}\left(\f{1}{2\ga}\int_0^{x_0}dt\, t \big|\wti V(t)\big|\right)^k,
\quad x\in(0,x_0), \; k\geq 0,
\end{align}
then imply
\begin{align} \lb{4.93}
\theta(x) = \sum_{k=0}^{\infty} \theta_k (x),
\end{align}
where the sum converges absolutely and uniformly on any compact subset of
$(0,x_0)$. In addition,
\begin{align} \lb{4.94}
\abs{\theta(x)} \leq \sum_{k=0}^\infty\abs{\theta_k (x)}
\leq x^{1/2-\ga} \exp\bigg(\f{1}{2\ga}\int_0^{x_0}dt\, t \big|\wti
V(t)\big|\bigg),
\quad x\in(0,x_0).
\end{align}
Since $\wti V\in L^1([0,b];x\,dx)$, there exists $x_0\in(0,b)$
such that
\begin{align} \lb{4.95}
\f{1}{2\ga}\int_0^{x_0}dt \, t\big|\wti V(t)\big| \leq \ln(3/2),
\end{align}
and hence by \eqref{4.91}, \eqref{4.93}, \eqref{4.94}, and
\eqref{4.95},
\begin{align} \lb{4.96}
\theta(x) \geq 2\theta_0 - \sum_{k=0}^{\infty}
\big|\theta_k (x)\big| \geq x^{1/2-\ga} \Big(2-e^{\ln(3/2)}\Big) \geq
\f{1}{2}x^{1/2-\ga}, \quad x\in(0,x_0).
\end{align}
Thus, $\theta \notin L^2((0,x_0); dx)$ and hence $\tau_+$ is in the
limit point case at $x=0$.
\end{proof}

Moreover, by the ``variation of constants'' formula, the Weyl--Titchmarsh
solution $\wti\phi(z,\cdot)$ of
\begin{align}
&-\psi''(z,x) + V(x)\psi(z,x) = z \psi(z,x),
\quad  x\in(0,\infty),   \lb{4.67} \\
&\; \psi(z,\cdot) \in L^2((0,b);dx) \, \text{ for some
$b\in(0,\infty)$, $z\in\bbC$} \lb{4.68}
\end{align}
satisfies the Volterra integral equation
\begin{align} \lb{4.69}
\wti\phi(z,x) = x^{1/2+\ga} + \f{1}{2\ga} \int_0^x dt \,
\big[x^{1/2+\ga}t^{1/2-\ga} - t^{1/2+\ga}x^{1/2-\ga}\big]
U(z,t)\wti\phi(z,t),
\end{align}
where
\begin{align} \lb{4.70}
U(z,x) = \wti V(x) - z.
\end{align}
To verify this claim one iterates \eqref{4.69} to obtain a solution
$\wti\phi(z,x)$ of \eqref{4.67} in the form
\begin{align} \lb{4.71}
\wti\phi(z,x) = \sum_{k=0}^\infty \wti\phi_k (z,x), \quad
z\in\bbC, \,\; x\in(0,\infty),
\end{align}
where
\begin{align}
\wti\phi_0 (z,x) &= x^{1/2+\ga},  \no \\
\wti\phi_k (z,x) &= \f{1}{2\ga} \int_0^x dx'
\big[x^{1/2+\ga}(x')^{1/2-\ga} - (x')^{1/2+\ga}x^{1/2-\ga}\big]
U(z,x')\wti\phi_{k-1}(z,x'),  \no \\
& \hspace{5.5cm} k\in\bbN, \; z\in\bbC, \; x\in(0,\infty).  \lb{4.72}
\end{align}
Since $\wti\phi_k (z,x)$, $k\in\bbN$, is continuous in
$(z,x)\in\bbC\times(0,\infty)$, entire with respect to $z$ for all fixed
$x\in(0,\infty)$, and since
\begin{align} \lb{4.74}
\abs{\wti\phi_k (z,x)} &\leq
\f{x^{1/2+\ga}}{k!}\left(\f{1}{\ga}\int_0^xdx' x'
\abs{U(z,x')}\right)^k, \quad (z,x)\in K,
\end{align}
where $K$ is any compact subs      et of $\bbC\times (0,\infty)$, the
series in \eqref{4.71} converges absolutely and uniformly on $K$,
and hence $\wti\phi(z,x)$ is continuous in
$(z,x)\in\bbC\times(0,\infty)$ and entire in $z$ for all fixed
$x\in(0,\infty)$. Moreover, it follows from \eqref{4.71} and
\eqref{4.74} that
\begin{align}
\abs{\wti\phi(z,x)} &\leq x^{1/2+\ga} \exp\bigg(\f{1}{\ga}\int_0^xdx' x'
\abs{U(z,x')}\bigg), \quad (z,x)\in K,
\end{align}
and hence, $\wti\phi(z,\cdot)$ satisfies \eqref{4.68}. Summarizing these
considerations, $\wti\phi(z,\cdot)$ satisfies Hypotheses
\ref{h4.1}\,$(iii)$\,$(\alpha)$--$(\gamma)$.

While this represents just an elementary example, we now turn to a vast
class of singular potentials.

We first state the following auxiliary result.

\begin{lemma} \lb{l4.13}
Let $b\in(0,\infty)$ and $f, f' \in AC_{\loc}((0,b))$, $f$
real-valued, and $f(x)\neq 0$ for all $x\in (0,b)$. \\
$(i)$  Introduce
\begin{equation}
\eta_{\pm} (x)= 2^{-1/2}f(x) \exp\bigg(\pm \int_x^{x_0} dx' \, f(x')^{-2}
\bigg),  \quad x, x_0 \in (0,b).  \lb{4.105}
\end{equation}
Then $\eta_\pm$ represent a fundamental system of solutions of
\begin{equation}
-\psi''(x) + \bigg[\f{f''(x)}{f(x)} + \f{1}{f(x)^4} \bigg]\psi(x)=0, \quad
x\in (0,b)
\lb{4.106}
\end{equation}
and
\begin{equation}
W(\eta_+,\eta_-)(x) = 1.  \lb{4.107}
\end{equation}
$(ii)$ Assume in addition that $f\in L^2([0,b'];dx)$ for some
$b'\in(0,b)$ and $\wti V\in L^1([0,c];f^2dx)$ for all $c\in(0,b)$. Then
there exists an entire Weyl--Titchmarsh solution $\wti \phi(z,\cdot)$ of
\begin{equation}
-\phi''(z,x) +\bigg[\f{f''(x)}{f(x)} +\f{1}{f(x)^4} + \wti V(x)
\bigg]\phi(z,x)=z\phi(z,x),
\quad z\in\bbC, \; x\in(0,b)
\end{equation}
in the following sense: \\
\indent $(\alpha)$ For all $x\in(0,b)$, $\wti\phi(\cdot,x)$ is
entire. \\
\indent $(\beta)$ $\wti\phi(z,x)$, $x\in(0,b)$, is real-valued for
$z\in\bbR$. \\
\indent $(\gamma)$ $\wti\phi(z,\cdot)$ satisfies the $L^2$-condition near
the end point $0$ and hence
\begin{equation}
\wti\phi(z,\cdot)\in L^2([0,c];dx) \, \text{ for all
$z\in\bbC$ and all $c\in(0,b)$.}
\end{equation}
\end{lemma}
\begin{proof}
Verifying item $(i)$ is a straightforward computation. To verify item
$(ii)$, consider the Volterra integral equation
\begin{align}
& \wti \phi(z,x)=\eta_-(x)+\int_0^x dx' \,
[\eta_+(x')\eta_-(x)-\eta_+(x)\eta_-(x')] \big[\wti V(x')-z
\big]\wti\phi(z,x'), \no  \\
& \hspace*{8cm}  z\in\bbC, \; x\in(0,b). \lb{4.110}
\end{align}
Again, iterating \eqref{4.110} then yields
\begin{align}
&\wti\phi(z,x)=\sum_{k=0}^\infty \wti\phi_k(z,x), \quad
\wti\phi_0(z,x)=\eta_-(x), \\
& \wti \phi_k(z,x)=\int_0^x dx' \,
[\eta_+(x')\eta_-(x)-\eta_+(x)\eta_-(x')] \big[\wti V(x')-z
\big]\wti\phi_{k-1}(z,x'), \quad k\in\bbN.
\end{align}
The elementary estimate
\begin{equation}
\bigg|\f{\eta_+(x)}{\eta_-(x)}\f{\eta_-(x')}{\eta_+(x')}\bigg| \leq
\exp\bigg(-\int_{x'}^x dy\, f(y)^{-2}\bigg)\leq 1, \quad 0\leq x' \leq x<b
\end{equation}
then yields
\begin{align}
\big|\wti\phi_1(z,x) \big|& \leq |\eta_-(x)| \int_0^x dx' \,
|\eta_+(x')\eta_-(x')|\bigg|1+\f{\eta_+(x)}{\eta_-(x)}
\f{\eta_-(x')}{\eta_+(x')}\bigg| \big|\wti V(x')-z \big|  \no \\
& \leq  |\eta_-(x)| \int_0^x dx' f(x')^2|\wti V(x') -z|
\end{align}
and hence
\begin{equation}
\big|\wti\phi_k(z,x) \big| \leq |\eta_-(x)|\f{1}{k!} \bigg(
 \int_0^x dx' f(x')^2 \big|\wti V(x)-z \big|\bigg)^k,
\quad  k\in\bbN, \; z\in\bbC, \; x\in(0,b).
\end{equation}
Thus,
\begin{equation}
\big|\wti\phi(z,x) \big| \leq |\eta_-(x)|\exp\bigg(
 \int_0^x dx' f(x')^2 \big|\wti V(x)-z \big|\bigg),
\quad  k\in\bbN, \; z\in\bbC, \; x\in(0,b).
\end{equation}
This proves items $(ii)\, (\alpha)$ and $(ii)\, (\beta)$. Since by
hypothesis, $f\in L^2([0,b'];dx)$ for some $b'\in(0,b)$ and hence
$\eta_-\in L^2([0,c];dx)$ for all $c\in(0,b)$, item $(ii)\, (\gamma)$
holds as well.
\end{proof}

A general class of examples of strongly singular potentials satisfying
Hypothesis \ref{h4.1}\,$(iii)$ is then described in the following example.

\begin{example} \lb{e4.14}
Let $b\in(0,\infty)$. Then the class of potentials $V$ such that
\begin{align}
& V, V' \in AC_{\loc}((0,b)), \, V \in L^1_\loc((0,\infty);dx), \,
\text{ $V$ real-valued},  \lb{4.117}
\\ & V(x)>0, \quad x\in (0,b),  \lb{4.118} \\
& V^{-1/2} \in L^1 ([0,b];dx),  \lb{4.119}  \\
& V' V^{-5/4} \in L^2([0,b];dx), \lb{4.120}  \\
& \text{either $V^{-3/2}V'' \in L^1([0,b];dx)$, or else,}  \lb{4.120a}
\\
& \text{$V''>0$ a.e.\ on $(0,b)$
        and $\lim_{x\downarrow 0} V'(x)V(x)^{-3/2}$ exists and is
finite,}   \lb{4.121}
\end{align}
satisfies Hypothesis \ref{h4.1}\,$(iii)$ $(\alpha)$--$(\gamma)$ in the
following sense: There exists an entire Weyl--Titchmarsh solution $\wti
\phi(z,\cdot)$ of
\begin{equation}
-\psi''(z,x) + V(x) \psi(z,x)=z\psi(z,x), \quad z\in\bbC, \;
x\in(0,\infty)
\end{equation}
satisfying the following conditions $(\alpha)$--$(\beta)$: \\
\indent $(\alpha)$ For all $x\in(0,\infty)$, $\wti\phi(\cdot,x)$ is
entire. \\
\indent $(\beta)$ $\wti\phi(z,x)$, $x\in(0,\infty)$, is real-valued for
$z\in\bbR$. \\
\indent $(\gamma)$ $\wti\phi(z,\cdot)$ satisfies the $L^2$-condition near
the end point $0$ and hence
\begin{equation}
\wti\phi(z,\cdot)\in L^2([0,c];dx) \, \text{ for all $z\in\bbC$ and all
$c\in(0,\infty)$.}
\end{equation}
\end{example}

Since $V$ is strongly singular at most at $x=0$, it suffices to
discuss this  example for $x\in (0,b)$ only. Moreover, for simplicity, we
focus only on sufficient conditions for Hypotheses \ref{h4.1}\,$(iii)$
$(\alpha)$--$(\gamma)$ to hold. The additional limit point assumptions on
$V$ at zero and at infinity can easily be supplied (cf.\
\cite[Sects.\ XIII.6, XIII.9, XIII.10]{DS88}). Moreover, we made no
efforts to optimize the conditions on $V$. The point of the example is
just to show the wide applicability of our approach based on Hypothesis
\ref{h4.1}.

In order to reduce Example \ref{e4.14} to Lemma \ref{l4.13}, one can
argue as follows: Introduce
\begin{align}
f(x)&=V(x)^{-1/4},  \\
\widetilde V(x) &= -f''(x)/f(x).
\end{align}
Then $f, f'\in AC_{\loc}((0,b))$, $f\neq 0$ on $(0,b)$, and $f\in
L^2([0,c];dx)$ for all $c\in(0,b)$. Moreover, since
\begin{equation}
f^2\wti V =-f f''=-\f{5}{16} \big[V^{-5/4}V' \big]^2+\f{1}{4}V^{-3/2} V'',
\end{equation}
$\wti V \in L^1([0,c];f^2 dx)$ for some $c\in (0,b)$. (This is clear from
\eqref{4.120} if condition in \eqref{4.120a} is assumed. In case
\eqref{4.121} is assumed, a straightforward
integration by parts, using \eqref{4.120}, yields
$\wti V \in L^1([0,c];f^2 dx)$ for some $c\in (0,b)$.) Thus, Lemma
\ref{l4.13} applies to
\begin{equation}
V=f^{-4}=[(f''/f)+f^{-4}]+\wti V.
\end{equation}

\begin{remark} \lb{r4.15}
We focused on the strongly singular case where $\tau_+$ is in the limit
point case at the singular endpoint $x=a$. The singular case, where $V$ is
not integrable at the endpoint $a$ and $\tau_+$ is in the limit circle
case at $a$ is similar to the regular case (associated with a
Weyl--Titchmarsh coefficient having the Herglotz property) considered in
Section \ref{s2}. For pertinent references to this case see \cite{Fu77},
\cite{Fu04}.
\end{remark}

\section{An Illustrative Example} \lb{s4}

In this section we provide a detailed treatment of the following
well-known singular potential example (which fits into Lemma \ref{l4.13} with
$f(x)=(x/\gamma)^{1/2}$, $x>0$, $\gamma \in [1,\infty)$, and $\wti V=0$),
\begin{equation}
V(x,\gamma)=\f{\gamma^2-(1/4)}{x^2}, \quad x\in(0,\infty), \;
\gamma\in [1,\infty)  \lb{3.1}
\end{equation}
with associated differential expression
\begin{equation}
\tau_+(\gamma) = -\f{d^2}{dx^2} + V(x,\gamma), \quad
x\in(0,\infty), \; \gamma\in [1,\infty). \lb{3.2}
\end{equation}
Numerous references have been devoted to this example, we refer, for
instance, to \cite{DLS04}, \cite{DS00}, \cite[p.\ 1532--1536]{DS88},
\cite{EK05}, \cite{Fu04}, \cite{FP98}, \cite{Me64},
\cite[p.\ 142--144]{Na68}, \cite{Na74}, \cite[p.\ 87--90]{Ti62}, and the
literature therein. The corresponding maximally defined self-adjoint
Schr\"odinger operator $H_+(\gamma)$ in $L^2([0,\infty);dx)$ is then defined
by
\begin{align}
&H_+(\gamma)f=\tau_+(\gamma) f,  \no  \\
&f\in\dom(H_+(\gamma))=\{g\in L^2([0,\infty);dx)\,|\, g, g'\in
AC_{\loc}((0,\infty));  \lb{3.3} \\
& \hspace*{5.8cm} \tau_+(\gamma) g \in L^2([0,\infty);dx)\}.  \no
\end{align}
The potential $V(\cdot,\gamma)$ in \eqref{3.1} is so strongly singular
at the finite end point $x=0$ that $H_+(\gamma)$ (in stark contrast to
cases regular at $x=0$, cf.\ \eqref{2.3}) is self-adjoint in
$L^2([0,\infty);dx)$ without imposing any boundary condition at $x=0$.
Equivalently, the corresponding minimal Schr\"odinger operator
$\widetilde H_+(\gamma)$, defined by
\begin{align}
&\widetilde H_+(\gamma)f=\tau_+(\gamma) f,  \no  \\
&f\in\dom(\widetilde H_+(\gamma))=\{g\in L^2([0,\infty);dx)\,|\, g,
g'\in AC_{\loc}((0,\infty));  \lb{3.3a} \\
& \hspace*{1.55cm} \text{${\rm supp}\,(g)\subset (0,\infty)$ compact;}
\; \tau_+(\gamma) g \in L^2([0,\infty);dx)\},  \no
\end{align}
is essentially self-adjoint in $L^2([0,\infty);dx)$.

A fundamental system of solutions of
\begin{equation}
(\tau_+(\gamma)\psi)(z,x)=z \psi(z,x), \quad x\in(0,\infty) \lb{3.4}
\end{equation}
is given by
\begin{equation}
x^{1/2} J_\gamma(z^{1/2}x), \; x^{1/2} Y_\gamma(z^{1/2}x),
\quad z\in\bbC\backslash\{0\}, \; x\in(0,\infty), \; \gamma \in [1,\infty)
\lb{3.5}
\end{equation}
with $J_\gamma(\cdot)$ and $Y_\gamma(\cdot)$ the usual Bessel functions
of order $\gamma$ (cf.\ \cite[Ch.\ 9]{AS72}). We first treat the
case where
\begin{equation}
\gamma\in(1,\infty), \quad \gamma \notin \bbN, \lb{3.5a}
\end{equation}
in which case
\begin{equation}
x^{1/2} J_\gamma(z^{1/2}x), \; x^{1/2} J_{-\gamma}(z^{1/2}x),
\quad z\in\bbC\backslash\{0\}, \;
x\in(0,\infty), \; \gamma\in(1,\infty)\backslash\bbN  \lb{3.6}
\end{equation}
is a fundamental system of solutions of \eqref{3.4}. Since the system
of solutions in \eqref{3.6} exhibits the branch cut $[0,\infty)$ with
respect to $z$, we slightly change it into the following system,
\begin{align}
\begin{split}
\phi(z,x,\gamma)&=C^{-1}\pi [2\sin(\pi\gamma)]^{-1} z^{-\gamma/2}
x^{1/2} J_\gamma(z^{1/2}x), \\
\theta(z,x,\gamma)&=Cz^{\gamma/2} x^{1/2} J_{-\gamma}(z^{1/2}x),
\quad z\in\bbC, \;
x\in(0,\infty), \; \gamma\in(1,\infty)\backslash\bbN, \lb{3.7}
\end{split}
\end{align}
which for each $x\in(0,\infty)$ represents entire functions with respect to
$z$. Here $C\in\bbR\backslash\{0\}$ is a normalization constant to be
discussed in Remark \ref{r4.4}. One verifies that (cf.\ \cite[p.\
360]{AS72})
\begin{equation}
W(\theta(z,\cdot,\gamma), \phi(z,\cdot,\gamma)) = 1, \quad
z\in\bbC, \; \gamma\in(1,\infty)\backslash\bbN \lb{3.7a}
\end{equation}
and that (cf.\ \cite[p.\ 360]{AS72})
\begin{align}
\begin{split}
& z^{\mp\gamma/2} x^{1/2} J_{\pm\gamma} (z^{1/2}x)=2^{-\gamma}
x^{(1/2)\pm\gamma} \sum_{k=0}^\infty
\f{(-zx^2/4)^k}{k!\Gamma(k+1 \pm\gamma)}, \lb{3.8}  \\
& \hspace*{3.45cm} z\in\bbC, \; x\in(0,\infty), \;
\gamma\in(1,\infty)\backslash\bbN.
\end{split}
\end{align}
Hence the fundamental system $\phi(z,\cdot,\gamma),
\theta(z,\cdot,\gamma)$ in \eqref{3.7} of solutions of \eqref{3.4} is
entire with respect to $z$ and real-valued for $z\in\bbR$.

The corresponding solution of \eqref{3.4}, square integrable in a
neighborhood of infinity, is given by
\begin{align}
\begin{split}
& x^{1/2} H^{(1)}_\gamma (z^{1/2}x)=\f{i}{\sin(\pi\gamma)}x^{1/2}
\big[e^{-i\pi\gamma}J_\gamma(z^{1/2}x)-J_{-\gamma}(z^{1/2}x)\big], \\
& \hspace*{3.45cm} z\in\bbC\backslash [0,\infty), \; x\in(0,\infty), \;
\gamma\in(1,\infty)\backslash\bbN \lb{3.9}
\end{split}
\end{align}
with $H_\gamma^{(1)}(\cdot)$ the usual Hankel function of order
$\gamma$ (cf.\ \cite[Ch.\ 9]{AS72}). In order to be compatible
with our modified system $\phi, \theta$ of solutions of
\eqref{3.4}, we replace it by
\begin{align}
\psi_+(z,x,\gamma)&= Cz^{\gamma/2}x^{1/2}J_{-\gamma}(z^{1/2}x)
-C^{2} e^{-i\pi\gamma} z^{\gamma} C^{-1}
z^{-\gamma/2}x^{1/2}J_\gamma(z^{1/2}x) \no \\
&= \theta(z,x,\gamma) + m_+(z,\gamma) \phi(z,x,\gamma), \lb{3.10} \\
& \hspace*{.5cm} z\in\bbC\backslash [0,\infty), \; x\in(0,\infty), \;
\gamma\in(1,\infty)\backslash\bbN,  \no
\end{align}
where
\begin{equation}
m_+(z,\gamma)=-C^2 (2/\pi)\sin(\pi\gamma) e^{-i\pi\gamma}z^\gamma,
\quad z\in \bbC\backslash [0,\infty), \;
\gamma\in(1,\infty)\backslash\bbN  \lb{3.11}
\end{equation}
and
\begin{equation}
\ol{m_+(z,\gamma)}=m_+(\ol z,\gamma), \quad z\in \bbC\backslash
[0,\infty).
\lb{3.11a}
\end{equation}

Next, we consider the case,
\begin{align} \lb{3.11b}
\gamma = n \in \bbN,
\end{align}
in which
\begin{equation}
x^{1/2}J_n(z^{1/2}x), \; x^{1/2}Y_n(z^{1/2}x), \quad
z\in\bbC\backslash\{0\}, \; x\in(0,\infty), \; n\in\bbN,  \lb{3.11ba}
\end{equation}
is a fundamental system of solutions of \eqref{3.4}. As before, we
slightly change it into the following system,
\begin{align}
\phi(z,x,n) &= C^{-1} (\pi/2) z^{-n/2} x^{1/2}
J_n(z^{1/2}x), \no \\
\theta(z,x,n) &= C z^{n/2} x^{1/2} \big[ -Y_n(z^{1/2}x) +
\pi^{-1}\ln(z)J_n(z^{1/2}x) \big],  \lb{3.11c} \\
& \hspace*{3.65cm} z\in\bbC, \; x\in(0,\infty), \; n\in\bbN. \no
\end{align}
Here $C\in\bbR\backslash\{0\}$ is a normalization constant to be
discussed in Remark \ref{r4.4}. One verifies that (cf.\ \cite[p.\
360]{AS72})
\begin{equation}
W(\theta(z,\cdot,n), \phi(z,\cdot,n))(x) = 1, \quad z\in\bbC, \; n\in
\bbN, \lb{3.11d}
\end{equation}
and that the fundamental system of solutions of \eqref{3.4},
$\phi(z,\cdot,n), \theta(z,\cdot,n)$ in \eqref{3.11c}, is entire
with respect to $z$ and real-valued for $z\in\bbR$.

The corresponding solution of \eqref{3.4}, square integrable in a
neighborhood of infinity, is given by
\begin{align}
\begin{split}
& x^{1/2} H^{(1)}_n (z^{1/2}x)= x^{1/2} \big[J_n(z^{1/2}x) +
iY_n(z^{1/2}x)\big],  \\
& \hspace*{2.5cm} z\in\bbC\backslash [0,\infty),x\in(0,\infty), \;
n\in \bbN \lb{3.11e}
\end{split}
\end{align}
with $H_n^{(1)}(\cdot)$ the usual Hankel function of order $n$
(cf.\ \cite[Ch.\ 9]{AS72}). In order to be compatible with our
modified system $\phi, \theta$ of solutions of \eqref{3.4}, we
replace it by
\begin{align}
\psi_+(z,x,n)&= Cz^{n/2}x^{1/2}iH_n(z^{1/2}x) = C z^{1/2} x^{1/2}
\big[-Y_n(z^{1/2}x) + iJ_n(z^{1/2}x)\big] \no \\
&= \theta(z,x,n) + m_+(z,n) \phi(z,x,n), \lb{3.11f} \\
& \; z\in\bbC\backslash [0,\infty), \; x\in(0,\infty), \; n\in \bbN, \no
\end{align}
where
\begin{equation}
m_+(z,n)=C^2 (2/\pi) z^n\big[i-(1/\pi)\ln(z)\big], \quad z\in
\bbC\backslash [0,\infty), \; n\in\bbN \lb{3.11g}
\end{equation}
and
\begin{equation}
\ol{m_+(z,n)}=m_+(\ol z,n), \quad z\in \bbC\backslash [0,\infty).
\lb{3.11h}
\end{equation}

\begin{remark} \lb{r3.1} $(i)$
We emphasize that in stark contrast to the case of regular
half-line Schr\"odinger operators in Section \ref{s2},
$m_+(\cdot,\gamma)$ in \eqref{3.11} and \eqref{3.11g} is not a
Herglotz function for $\gamma\in[1,\infty)$. \\
$(ii)$ After finishing our
paper, we received a manuscript by Everitt and Kalf \cite{EK05} in which
the Friedrichs extension and the associated Hankel eigenfunction transform
are treated in detail for the case $\gamma\in [0,1)$ in \eqref{3.1}. In
this case the corresponding Weyl--Titchmarsh coefficient turns out to be a
Herglotz function.
\end{remark}

Since $\tau_+(\gamma)$ is in the limit point case at $x=0$ and at
$x=\infty$, \eqref{3.4} has a unique solution (up to constant
multiples) that is $L^2$ near $0$ and $L^2$ near $\infty$. Indeed,
that unique $L^2$-solution near $0$ (up to normalization) is
precisely $\phi(z,\cdot,\gamma)$; similarly, the unique
$L^2$-solution near $\infty$ (up to normalization) is
$\psi_+(z,\cdot,\gamma)$.

By \eqref{3.7a} and \eqref{3.11d}, a computation of the Green's function
$G_+(z,x,x',\gamma)$ of $H_+(\gamma)$ yields
\begin{align}
G_+(z,x,x',\gamma)&=\f{i\pi}{2}\begin{cases} x^{1/2}J_\gamma(z^{1/2}x)
x'^{1/2}H^{(1)}_\gamma(z^{1/2}x'), & 0<x\leq x', \\
x'^{1/2}J_\gamma(z^{1/2}x')
x^{1/2}H^{(1)}_\gamma(z^{1/2}x), & 0<x'\leq x, \end{cases} \lb{3.13} \\
&=\begin{cases} \phi(z,x,\gamma) \psi_+(z,x',\gamma), & 0<x\leq x', \\
\phi(z,x',\gamma)\psi_+(z,x,\gamma), & 0<x'\leq x, \end{cases}
\lb{3.14} \\
& \hspace*{2.17cm} z\in\bbC\backslash [0,\infty), \; \gamma\in[1,\infty).
\no
\end{align}
Thus,
\begin{align}
\begin{split}
&((H_{+}(\gamma)-zI)^{-1}f)(x)
=\int^{\infty}_{0} dx' \, G_{+}(z,x,x',\gamma)f(x'), \lb{3.15} \\
& z\in\bbC\backslash [0,\infty), \; x\in(0,\infty), \;
f\in L^{2}([0,\infty);dx), \; \gamma\in[1,\infty).
\end{split}
\end{align}

Given $m_+(z,\gamma)$ in \eqref{3.11}, we define the associated measure
$\rho_+(\cdot,\gamma)$ by
\begin{align}
\rho_+((\lambda_1,\lambda_2],\gamma)&=
\pi^{-1}\lim_{\delta\downarrow 0} \lim_{\varepsilon\downarrow 0}
\int_{\lambda_1+\delta}^{\lambda_2+\delta} d\lambda \,
\Im(m_+(\lambda+i\varepsilon,\gamma)) \lb{3.16} \\
&=C^2 \f{\lambda_2^{\gamma+1} -\lambda_1^{\gamma+1}}{\gamma+1}
\f{2}{\pi^2}\begin{cases}\sin^2(\pi\gamma), & \ga\notin\bbN, \\ 1,
& \ga\in\bbN,\end{cases} \lb{3.17} \\
& \hspace*{2.65cm} 0\leq\lambda_1<\lambda_2, \; \gamma
\in [1,\infty), \no
\end{align}
generated by the function
\begin{align}
\rho_+(\lambda,\gamma)= C^2 \chi_{[0,\infty)}(\lambda)
\f{\lambda^{\gamma+1}}{\gamma+1} \f{2}{\pi^2} \begin{cases}
\sin^2(\pi\gamma), & \ga \notin \bbN, \\ 1, & \ga \in \bbN,
\end{cases} \quad \lambda\in\bbR, \; \gamma \in [1,\infty).
\lb{3.18}
\end{align}

Even though $m_+(\cdot,\gamma)$ is not a Herglotz function for
$\gamma\in[1,\infty)$, $d\rho_+(\cdot,\gamma)$ is defined as in
\eqref{4.33}, in analogy to the case of Herglotz functions
discussed in Appendix \ref{sA} (cf.\ \eqref{A.4}).

Next, we introduce the family of spectral projections,
$\{E_{H_{+}(\gamma)}(\lambda)\}_{\lambda\in\bbR}$, of the
self-adjoint operator $H_{+}(\gamma)$ and note that for $F\in C(\bbR)$,
\begin{align}
&\big(f,F(H_{+}(\gamma))g\big)_{L^2([0,\infty);dx)}= \int_{\bbR}
d\big(f,E_{H_{+}(\gamma)}(\lambda)g\big)_{L^2([0,\infty);dx)}\,
F(\lambda), \no \\ & f, g \in\dom(F(H_{+}(\gamma))) \lb{3.19} \\
& \quad =\bigg\{h\in L^2([0,\infty);dx)
\,\bigg|\,
\int_{\bbR} d\big\|E_{H_{+}(\gamma)}(\lambda)h\big\|_{L^2([0,\infty);dx)}^2
\, |F(\lambda)|^2 < \infty\bigg\}. \no
\end{align}

The connection between $\{E_{H_{+}(\gamma)}(\lambda)\}_{\lambda\in\bbR}$
and $\rho_+(\lambda,\gamma)$, $\lambda\geq0$, is described in the next
result.

\begin{lemma} \lb{l3.2}
Let $\gamma \in [1,\infty)$, $f,g \in C^\infty_0((0,\infty))$, $F\in
C(\bbR)$, and $\lambda_1, \lambda_2 \in [0,\infty)$,
$\lambda_1<\lambda_2$. Then,
\begin{align}
\begin{split}
&\big(f,F(H_{+}(\gamma))E_{H_{+}(\gamma)}((\lambda_1,\lambda_2])
g\big)_{L^2([0,\infty);dx)}   \\
& \quad = \big(\hatt
f_{+}(\gamma),M_FM_{\chi_{(\lambda_1,\lambda_2]}} \hatt
g_{+}(\gamma)\big)_{L^2(\bbR;d\rho_{+}(\cdot,\gamma))}, \lb{3.20}
\end{split}
\end{align}
where
\begin{equation}
\hatt h_{+}(\gamma)(\lambda)=\int_0^\infty dx \,
\phi(\lambda,x,\gamma)h(x), \quad \lambda \in [0,\infty), \;
h\in C^\infty_0((0,\infty)),  \lb{3.21}
\end{equation}
and $M_G$ denotes the operator of multiplication by the
$d\rho_+(\cdot,\gamma)$-measurable function $G$ in the Hilbert space
$L^2(\bbR;d\rho_{+}(\cdot,\gamma))$.
\end{lemma}

The proof of Lemma \ref{l3.2} is a special case of that of Theorem
\ref{t4.6} and hence omitted.

As in Section \ref{s3} one can remove the compact support restrictions on
$f$ and $g$ in Lemma \ref{l3.2}. To this end one considers the map
\begin{equation}
U_{+}(\gamma)\colon \begin{cases}L^2([0,\infty);dx)\to L^2(\bbR;
d\rho_{+}(\cdot,\gamma)) \\[1mm]
\hspace*{1.95cm}  h \mapsto \hatt h_{+}(\cdot,\gamma)=
\slimes_{b\uparrow\infty}\int_0^b dx\, \phi(\cdot,x,\gamma) h(x),
\end{cases}  \lb{3.36}
\end{equation}
where $\slimes$ refers to the $L^2(\bbR;d\rho_+(\cdot,\gamma))$-limit.

In addition, it is of course known (cf., e.g., \cite[p.\ 1535]{DS88})
that  the Bessel transform $U_{+}(\gamma)$ in \eqref{3.36} is onto
and hence that $U_{+}(\gamma)$ is unitary with
\begin{equation}
U_{+}(\gamma)^{-1} \colon \begin{cases}
L^2(\bbR;d\rho_{+}(\cdot,\gamma))
\to  L^2([0,\infty);dx)  \\[1mm]
\hspace*{2.2cm}
\hatt h \mapsto \slimes_{\mu_1\downarrow -\infty,
\mu_2\uparrow\infty}
\int_{\mu_1}^{\mu_2} d\rho_{+}(\lambda,\gamma)\,
\phi(\lambda,\cdot,\gamma) \hatt h(\lambda), \end{cases} \lb{3.41}
\end{equation}
where $\slimes$ refers to the $L^2([0,\infty);dx)$-limit.

Again we sum up these considerations in a variant of the spectral
theorem for (functions of) $H_{+}(\gamma)$.

\begin{theorem} \lb{t3.3}
Let $\gamma \in [1,\infty)$, $F\in C(\bbR)$. Then,
\begin{equation}
U_{+}(\gamma) F(H_{+}(\gamma))U_{+}(\gamma)^{-1} = M_F \lb{3.42}
\end{equation}
in $L^2(\bbR;d\rho_{+}(\cdot,\gamma))$. Moreover,
\begin{align}
& \sigma(F(H_{+}(\gamma)))= \essran_{d\rho_{+}(\cdot,\gamma)}(F),
\lb{3.43m} \\
& \sigma(H_{+}(\gamma))=\supp\,(d\rho_{+}(\cdot,\gamma)), \lb{3.44}
\end{align}
and the spectrum of $H_{+}(\gamma)$ is simple.
\end{theorem}

Next, we reconsider spectral theory for $H_+(\gamma)$ by choosing a
reference point $x_0\in(0,\infty)$ away from the singularity of
$V(\cdot,\gamma)$ at $x=0$.

Consider a system $\phi(z,\cdot,x_0,\gamma)$, $\theta(z,\cdot,x_0,\gamma)$
of solutions of \eqref{3.4} with the following initial conditions at the
reference point $x_0\in(0,\infty)$,
\begin{align*}
\phi(z,x_0,x_0,\gamma) = \theta'(z,x_0,x_0,\gamma) = 0, \quad
\phi'(z,x_0,x_0,\gamma) = \theta(z,x_0,x_0,\gamma) = 1.
\end{align*}
Denote by $m_\pm(z,x_0,\gamma)$ two Weyl--Titchmarsh $m$-functions
corresponding to the restriction of our problem to the intervals
$(0,x_0]$ and $[x_0,\infty)$, respectively. Then the
Weyl--Titchmarsh  solutions
$\psi_{\pm}(z,\cdot,x_0,\gamma)$ and the $2\times 2$ matrix-valued
Weyl--Titchmarsh $M$-function $M(z,x_0,\gamma)$ are given by
\begin{align}
\psi_\pm(z,x,x_0,\gamma) &= \theta(z,x,x_0,\gamma) + m_\pm(z,x_0,\gamma)
\phi(z,x,x_0,\gamma),  \\
M(z,x_0,\gamma) &=\begin{pmatrix}
\f{1}{m_{-}(z,x_0,\gamma)-m_{+}(z,x_0,\gamma)} &
\f{1}{2}\f{m_{-}(z,x_0,\gamma)
+m_{+}(z,x_0,\gamma)}{m_{-}(z,x_0,\gamma)-m_{+}(z,x_0,\gamma)} \\
\f{1}{2}\f{m_{-}(z,x_0,\gamma)
+m_{+}(z,x_0,\gamma)}{m_{-}(z,x_0,\gamma)-m_{+}(z,x_0,\gamma)} &
\f{m_{-}(z,x_0,\gamma)m_{+}(z,x_0,\gamma)}{m_{-}(z,x_0,\gamma)
-m_{+}(z,x_0,\gamma)}\end{pmatrix}.
\end{align}
Since any $L^2$-solution near $0$ and near $\infty$ (i.e., any
Weyl--Titchmarsh solution) is necessarily proportional to
$x^{1/2}J_\ga(z^{1/2}x)$ and $x^{1/2}H_\ga^{(1)}(z^{1/2}x)$,
respectively, one explicitly computes for $m_\pm(z,x_0,\gamma)$,
\begin{align}
m_-(z,x_0,\gamma) &= \f{1}{2x_0} +
z^{1/2}\f{J_\ga'(z^{1/2}x_0)}{J_\ga(z^{1/2}x_0)}, \\
m_+(z,x_0,\gamma) &= \f{1}{2x_0} +
z^{1/2}\f{{H^{(1)}_\ga}'(z^{1/2}x_0)}{H^{(1)}_\ga(z^{1/2}x_0)},
\end{align}
and for $M(z,x_0,\gamma)$,
\begin{align}
M_{0,0}(z,x_0,\gamma) &= \f{i\pi x_0}{2}
J_\ga(z^{1/2}x_0)H^{(1)}_\ga(z^{1/2}x_0),
\\
M_{0,1}(z,x_0,\gamma) &= M_{1,0}(z,x_0,\gamma) = \f{i\pi}{4}
\Bigg[J_\ga(z^{1/2}x_0)H^{(1)}_\ga(z^{1/2}x_0) + x_0z^{1/2}
\\
& \qquad \times \Big(J_\ga(z^{1/2}x_0){H^{(1)}_\ga}'(z^{1/2}x_0) +
J'_\ga(z^{1/2}x_0)H^{(1)}_\ga(z^{1/2}x_0)\Big) \Bigg],  \no  \\
M_{1,1}(z,x_0,\gamma) &= \f{i\pi}{8x_0}
\Bigg[J_\ga(z^{1/2}x_0)H^{(1)}_\ga(z^{1/2}x_0) + 2x_0z^{1/2}  \no
\\
& \qquad \times \Big(J_\ga(z^{1/2}x_0){H^{(1)}_\ga}'(z^{1/2}x_0) +
J'_\ga(z^{1/2}x_0)H^{(1)}_\ga(z^{1/2}x_0)\Big)
\\
& \quad + 4x_0^2z J'_\ga(z^{1/2}x_0){H^{(1)}_\ga}'(z^{1/2}x_0)
\Bigg].  \no
\end{align}
Using \eqref{3.9}, \eqref{3.11e}, and the calculation above, one can
also compute the $2\times 2$ spectral measure
$d\Omega(\cdot,x_0,\gamma)$ and its density
$d\Omega(\cdot,x_0,\gamma)/d\lambda$,
\begin{align}
\f{d\Omega(\la,x_0,\gamma)}{d\lambda} &= \f{1}{\pi}\lim_{\ve\downarrow 0}
\Im(M(\la+i\ve,x_0, \ga)),  \quad \lambda\in\bbR, \\
\f{d\Omega_{0,0}(\la,x_0,\gamma)}{d\lambda} &= \begin{cases}\f{x_0}{2}
J_\ga(\la^{1/2}x_0)^2, & \la>0, \\ 0, & \la \leq 0, \end{cases}
\lb{3.60} \\
\f{d\Omega_{0,1}(\la,x_0,\gamma)}{d\lambda} &=
\f{d\Omega_{1,0}(\la,x_0,\gamma)}{d\la}
\no \\
&=
\begin{cases}
\f{1}{4}\left[J_\ga(\la^{1/2}x_0)^2 +
2x_0\la^{1/2}J_\ga(\la^{1/2}x_0)J_\ga'(\la^{1/2}x_0)\right], &
\la>0,
\\ 0, & \la \leq 0, \end{cases} \lb{3.61}
\\
\f{d\Omega_{1,1}(\la,x_0,\gamma)}{d\lambda} &=
\begin{cases} \f{1}{8x_0}\left[J_\ga(\la^{1/2}x_0) +
2x_0\la^{1/2}J_\ga'(\la^{1/2}x_0)\right]^2, & \la>0, \\ 0, & \la
\leq 0.
\end{cases}
\end{align}
Moreover, one verifies that,
\begin{align}
\rank\bigg(\f{d\Omega(\la,x_0,\gamma)}{d\lambda}\bigg)
= \begin{cases} 1, & \la > 0, \\ 0, & \la \leq 0. \end{cases}
\end{align}

Finally, we will show that the results of Section \ref{s3} which
let one obtain a scalar spectral measure $d\wti\rho_+(\la,\gamma)$ from
the $2 \times 2$ spectral measure $d\Omega(\la,x_0,\gamma)$ lead to the
measure equivalent to $d\rho_+(\la,\ga)$ obtained in the first part
of this section.

Let
\begin{align}
\wti\phi(z,x,\gamma) = z^{-\gamma/2}x^{1/2}J_\ga(z^{1/2}x) \lb{3.64}
\end{align}
be the Weyl--Titchmarsh solution satisfying Hypothesis
\ref{h4.1}\,$(iii)$. Inserting \eqref{3.60}, \eqref{3.61}, and
\eqref{3.64} into \eqref{4.52} then yields
\begin{align}
\f{d\wti\rho_+(\la,\gamma)}{d\lambda} = \begin{cases} \f{1}{2} \la^\ga, &
\la > 0,
\\ 0, & \la \leq 0, \end{cases}
\end{align}
which, up to a constant multiple, is the same as
$d\rho_+(\la,\ga)/d\lambda$ in \eqref{3.18}.

Of course the analogs of Theorem \ref{t4.7}, Theorem \ref{t4.9}, and
Corollary \ref{c4.10} all hold in the present context of the potential
\eqref{3.1}; we omit the details.

\begin{remark} \lb{r4.4}
We explicitly introduced the normalization constant
$C\in\bbR\backslash\{0\}$ in \eqref{3.7} and \eqref{3.11c} to determine
its effect on (the analog of) the Weyl--Titchmarsh coefficient $m_+$
(cf.\ \eqref{3.11} and \eqref{3.11g}) and the associated spectral
function $\rho_+$ (cf.\ \eqref{3.18}). As $C$ enters quadratically in
$m_+$ and $\rho_+$, it clearly has an effect on their asymptotic
behavior as $|z|\to\infty$, respectively, $|\lambda|\to\infty$. The same
observation applies of course in the regular half-line case considered
in the first half of Section \ref{s2}. It just so happens that in this
case the standard normalization of the fundamental system of solutions
$\phi_{\alpha}$ and $\theta_{\alpha}$ of \eqref{2.4} in \eqref{2.5}
represents a canonical choice and the normalization dependence can safely
be ignored. In the strongly singular case in Sections \ref{s3} and
\ref{s4} no such canonical choice of normalization exists. Of course, the
actual spectral properties of the corresponding half-line Schr\"odinger
operator are independent of such a choice of normalization.
\end{remark}

\appendix
\section{Basic Facts on Herglotz Functions} \lb{sA}
\setcounter{theorem}{0}
\setcounter{equation}{0}

In this appendix we recall the definition and basic properties of
Herglotz functions.

\begin{definition} \lb{dA.1}
Let $\bbC_\pm=\{z\in\bbC \mid \Im(z)\gtrless 0 \}$.
$m:{\mathbb{C_+}}\to {\mathbb{C}}$ is called a {\it Herglotz}
function (or {\it Nevanlinna} or {\it Pick} function) if $m$ is analytic
on
${\mathbb{C}}_+$ and
$m({\mathbb{C}}_+)\subseteq {\mathbb{C}}_+$.
\end{definition}

\smallskip

\noindent One then extends $m$ to $\bbC_-$ by reflection, that is, one
defines
\begin{equation}
m(z)=\overline{m(\overline z)},
\quad z\in{\mathbb{C}}_-. \lb{A.1}
\end{equation}
Of course, generally, \eqref{A.1} does not represent the analytic
continuation of $m\big|_{\bbC_+}$ into $\bbC_-$.

The fundamental result on Herglotz functions and their representations
as Borel transforms, in part due to Fatou, Herglotz, Luzin,
Nevanlinna, Plessner, Privalov, de la Vall{\'e}e Poussin, Riesz, and
others, then reads as follows.

\begin{theorem} \rm {(\cite{AG81}, Sect.\ 69, \cite{AD56},
\cite{Do74}, Chs.~II, IV,
\cite{KK74}, \cite{Ko80}, Ch.~6, \cite{Pr56}, Chs.~II,
IV, \cite{RR94}, Ch.~5.)} \label{tA.2}
Let $m$ be a Herglotz function. Then, \\
$(i)$ $m(z)$ has finite normal limits $m(\lambda
\pm i0)=\lim_{\varepsilon
\downarrow 0} m(\lambda \pm i\varepsilon)$ for
a.e.~$\lambda \in {\mathbb{R}}$. \\
$(ii)$ Suppose $m(z)$ has a zero normal limit on a
subset of
${\mathbb{R}}$ having positive Lebesgue measure. Then
$m\equiv 0$. \\
$(iii)$ There exists a Borel measure $d\omega$ on
${\mathbb{R}}$
satisfying
\begin{equation} \lb{A.2}
\int_{{\mathbb{R}}} \frac{d\omega (\lambda )}{1+\lambda^2}<\infty
\end{equation}
such that the Nevanlinna, respectively, Riesz-Herglotz
representation
\begin{align}
&m(z)=c+dz+\int_{{\mathbb{R}}}
d\omega (\lambda)\bigg[\frac{1}{\lambda -z}-\frac{\lambda}
{1+\lambda^2}\bigg], \quad z\in\bbC_+, \lb{A.3} \\[2mm]
&c=\Re(m(i)),\quad d=\lim_{\eta \uparrow
\infty}m(i\eta )/(i\eta )
\geq 0 \no
\end{align}
holds. Conversely, any function $m$ of the type \eqref{A.3} is a
Herglotz function. \\
$(iv)$ Let $\lambda _1,\lambda_2\in\bbR$, $\lambda_1<\lambda_2$. Then
the Stieltjes inversion formula for $d\omega$ reads
\begin{equation} \lb{A.4}
\omega((\lambda_1,\lambda_2])=\pi^{-1}\lim_{\delta\downarrow 0}
\lim_{\varepsilon\downarrow
0}\int^{\lambda_2+\delta}_{\lambda_1+\delta}d\lambda \, \Im(m(\lambda
+i\varepsilon)).
\end{equation}
$(v)$ The absolutely continuous $({\it ac})$ part $d\omega_{ac}$ of
$d\omega$ with respect to Lebesgue measure $d\lambda$ on
${\mathbb{R}}$ is
given by
\begin{equation}\lb{A.5}
d\omega_{ac}(\lambda)=\pi^{-1}\Im(m(\lambda+i0))\,d\lambda.
\end{equation}
$(vi)$ Local singularities of $m$ and $m^{-1}$ are necessarily
real and at most of first order in the sense that
\begin{align}
& \lim_{\epsilon\downarrow0}
(-i\epsilon)\, m(\lambda+i\epsilon) \geq 0, \quad
\lambda\in\bbR, \lb{A.6} \\
& \lim_{\epsilon\downarrow0}
(i\epsilon)\, m(\lambda+i\epsilon)^{-1}
\geq 0, \quad \lambda\in\bbR. \lb{A.7}
\end{align}
\end{theorem}

Further properties of Herglotz functions are collected in the
following theorem. We denote by
\begin{equation}
d\omega =d\omega_{\rm ac}+d\omega_{\rm sc}
+d\omega_{\rm pp} \lb{A.8}
\end{equation}
the decomposition of $d\omega$ into its absolutely continuous
$({\it ac})$, singularly continuous $({\it sc})$, and pure point $({\it
pp})$ parts with respect to Lebesgue measure on ${\mathbb{R}}$.

\begin{theorem} \rm{(\cite{AD56}, \cite{GT00}, \cite{KK74}, \cite{Si95},
\cite{Si95a}.)} \label{tA.3} Let $m$ be a Herglotz function with
representation \eqref{A.3}. Then, \\
$(i)$
\begin{align}
& d=0 \text{ and } \int_{{\mathbb{R}}}d\omega
(\lambda)(1+|\lambda|^s)^{-1}<\infty \text{ for some }
s\in (0,2)
\no \\
& \text{if and only if }
\int^\infty_1d\eta \, \eta^{-s} \, \Im (m(i\eta ))<\infty.
\lb{A.9}
\end{align}
$(ii)$ Let $(\lambda_1,\lambda_2)\subset {\mathbb{R}}$,
$\eta_1>0$.
Then there is a constant
$C(\lambda_1,\lambda_2,\eta_1)>0$ such that
\begin{equation}\lb{A.10}
\eta|m(\lambda+i\eta)|\leq
C(\lambda_1,\lambda_2,\eta_1),\quad (\lambda,\eta)\in
[\lambda_1,\lambda_2]\times (0,\eta_1).\end{equation}
$(iii)$
\begin{equation} \lb{A.11}
\sup_{\eta
>0}\eta |m(i\eta)|<\infty \text{ if and only if }
m(z)=\int_{{\mathbb{R}}}d\omega (\lambda)(\lambda-z)^{-1}
\text{ and }
\int_{{\mathbb{R}}}d\omega(\lambda)<\infty.
\end{equation}
In this case,
\begin{equation}\lb{A.12}
\int_{{\mathbb{R}}}d\omega
(\lambda)=\sup_{\eta >0}\eta |m(i\eta )|=-i\lim_{\eta \uparrow
\infty}\eta m(i\eta).
\end{equation}
$(iv)$ For all $\lambda \in {\mathbb{R}}$,
\begin{align}\lb{A.13}
&\lim_{\varepsilon\downarrow 0}\varepsilon
\Re(m(\lambda+i\varepsilon ))=0,\\
&\omega(\{\lambda\})=\lim_{\varepsilon\downarrow 0}\varepsilon
\Im
(m(\lambda+i\varepsilon))=-i\lim_{\varepsilon \downarrow 0}
\varepsilon
m(\lambda+i\varepsilon ).\lb{A.14}
\end{align}
$(v)$ Let $L>0$ and suppose $0\leq \Im (m(z))\leq L$
for all $z\in
{\mathbb{C}}_+$. Then $d=0$, $d\omega$ is purely absolutely
continuous, $d\omega = d\omega_{ac}$, and
\begin{equation} \lb{A.15}
0\leq
\frac{d\omega(\lambda)}{d\lambda}=\pi^{-1}\lim_{\varepsilon
\downarrow 0}\Im
(m(\lambda+i\varepsilon))\leq \pi^{-1}L \text{ for a.e. }
\lambda\in
{\mathbb{R}}.
\end{equation}
$(vi)$ Let $p\in (1,\infty)$, $[\lambda_3,\lambda_4]
\subset
(\lambda_1,\lambda_2)$, $[\lambda_1,\lambda_2]\subset
(\lambda_5,\lambda_6)$. If
\begin{equation} \lb{A.16}
\sup_{0<\varepsilon <1}\int^{\lambda_2}_{\lambda_1}d\lambda
\, |\Im (m(\lambda
+i\varepsilon))|^p < \infty,
\end{equation}
then $d\omega=d\omega _{ac}$ is purely
absolutely continuous on $(\lambda_1,\lambda_2)$,
$\frac{d\omega_{ac}}{d\lambda}\in L^p((\lambda_1,
\lambda_2);d\lambda)$, and
\begin{equation} \lb{A.17}
\lim_{\varepsilon\downarrow 0}\bigg\|\pi^{-1}\Im (m(\cdot
+i\varepsilon )) -
\frac{d\omega_{ac}}{d\lambda}\bigg\|_{L^p((\lambda_3,\lambda_4);d\lambda)}
= 0.
\end{equation}
Conversely, if $d\omega$ is purely absolutely continuous on
$(\lambda_5,\lambda_6)$, and if
$\frac{d\omega_{ac}}{d\lambda}\in$ \linebreak
$L^p((\lambda_5,\lambda_6);d\lambda)$,
then \eqref{A.16} holds. \\
$(vii)$ Let $(\lambda_1,\lambda_2)
\subset{\mathbb{R}}$.
Then a local version
of Wiener's theorem reads for $p\in (1,\infty)$,
\begin{align}
&\lim_{\varepsilon\downarrow 0}\varepsilon^{p-1}
\int^{\lambda_2}_{\lambda_1}d\lambda \, |\Im
(m(\lambda+i\varepsilon ))|^p \no \\
& \quad =\frac{\Gamma (\frac{1}{2})\Gamma (p-\frac{1}{2})}
{\Gamma (p)} \bigg[ \frac{1}{2}\omega (\{\lambda_1\})^p
+\frac{1}{2}\omega (\{\lambda_2\})^p+
\sum_{\lambda\in (\lambda_1,\lambda_2)}
\omega (\{\lambda \})^p
\bigg]. \lb{A.18}
\end{align}
Moreover, for $0<p<1$,
\begin{equation} \lb{A.19}
\lim_{\varepsilon \downarrow
0}\int^{\lambda_2}_{\lambda_1}d\lambda \, |\pi^{-1}\Im
(m(\lambda + i\varepsilon))|^p
=\int^{\lambda_2}_{\lambda_1}d\lambda
\left|\frac{d\omega_{ac}(\lambda)}{d\lambda}\right|^p.
\end{equation}
\end{theorem}

It should be stressed that Theorems \ref{tA.2} and \ref{tA.3} record
only the tip of an iceberg of results in this area. A substantial number
of additional references relevant in this context can be found in
\cite{GT00}.

\medskip

\noindent {\bf Acknowledgments.}
We are indebted to Norrie Everitt, 
Don Hinton, Hubert Kalf, Konstantin Makarov, Kwang Shin, and Gerald Teschl
for most helpful comments, and especially to Don Hinton for kindly pointing
out to us the origin of the use of Stone's formula in deriving spectral
functions (cf.\ the proofs of Theorems \ref{t2.5} and \ref{t2.9}) and the
origin of (an extension of) Lemma \ref{l4.12}. We also gratefully acknowledge the kind remarks and constructive comments by three anonymous referees.


\end{document}